\newcommand{\dsp}{\renewcommand{\baselinestretch}{1.5}}

\newcommand{\lr}{\longrightarrow}
\newcommand{\cd}{\cdots}

\newcommand{\N}{\mathbb{N}}
\newcommand{\C}{\mathbb{C}}
\newcommand{\R}{\mathbb{R}}

\newcommand{\Z}{\mathbb{Z}}
\newcommand{\T}{\mathbb{T}}

\documentclass[pdflatex,fleqn]{article}
\usepackage{}
\usepackage{xcolor}
\usepackage{mathrsfs}
\usepackage{amsfonts}
\usepackage{amssymb}
\usepackage{latexsym}
\usepackage{amsmath}
\usepackage{amscd}
\usepackage{graphicx}
\usepackage{hyperref}
\usepackage{hypernat}
\usepackage{cite}
\usepackage{mathtools}

\usepackage{geometry,amsthm,graphics,amssymb,amsmath,enumerate,
latexsym,tabularx,shapepar}
\usepackage[all]{xy}

\usepackage{extarrows}
\date{}

\dsp

\begin{document}

%\title{\bf A classification of inductive limit $C^{*}$-algebras  with ideal property}

\centerline{\large \bf A classification of inductive limit $C^{*}$-algebras  with ideal property}

\vspace{3mm}

\centerline{Guihua Gong, Chunlan Jiang, and Liangqing Li}

%\author{Guihua Gong, Chunlan Jiang, and Liangqing Li}

%\begin{document}

%\large

%\maketitle

%\vspace{-0.3in}

%\noindent{\underline {\noindent{ }\hspace {157mm}}}

\noindent{\bf Abstract}
 Let $A$ be an $AH$ algebra
%that is $A$ is the inductive limit $C^{*}$-algebra
% of
%$$A_{1}\xrightarrow{\phi_{1,2}}A_{2}\xrightarrow{\phi_{2,3}}A_{3}\longrightarrow\cdots\longrightarrow A_{n}\longrightarrow\cdots$$ with
 $A=\lim\limits_{n\to \infty}(A_{n}=\bigoplus\limits_{i=1}\limits^{t_{n}}P_{n,i}M_{[n,i]}(C(X_{n,i}))P_{n,i}, \phi_{n,m})$, where $X_{n,i}$ are compact metric spaces, $t_{n}$ and $[n,i]$ are
positive integers, and $P_{n,i}\in M_{[n,i]}(C(X_{n,i}))$ are projections. Suppose that $A$ has the ideal property: each closed two-sided ideal of $A$ is
generated by the projections inside the ideal, as a closed two sided ideal. In this article, we will classify all $AH$ algebras with ideal property of no
dimension growth---that is, $sup_{n,i}dim(X_{n,i})<+\infty$. This result generalizes and unifies the classification of $AH$ algebras of real rank zero in
[EG] and [DG] and the classification of simple $AH$ algebras in [G5] and [EGL1]. This completes one of two
important possible generalizations of [EGL1] suggested in the introduction of [EGL1]. The invariants for the classification include the scaled ordered total
$K$-group $(\underline{K}(A), \underline{K}(A)_{+},\Sigma A)$ (as already used in real rank zero case in [DG]), for each $[p]\in\Sigma A$, the tracial
state space $T(pAp)$ of cut down algebra $pAp$ with a certain compatibility, (which is used by [Stev] and [Ji-Jiang] for $AI$ algebras with the ideal property), and
a new ingredient, the invariant $U(pAp)/\overline{DU(pAp)}$ with a certain compatibility condition, where $\overline{DU(pAp)}$ is the closure of
commutator subgroup $DU(pAp)$ of the unitary group $U(pAp)$ of the cut down algebra $pAp$. In [GJL] a counterexample is presented to show that this new ingredient must be included in the invariant. The discovery of this new invariant is  analogous to that of the order structure on the total K-theory when one advances  from  the classification of simple real rank zero $C^*$-algebras to that  of non simple real rank zero $C^*$-algebras in [G2], [Ei],  [DL] and [DG] (see Introduction below).

%The notion of simple $TAI$ (and $TAF$) algebras was introduced by  Lin to axiomatize  the decomposition theorem of  [G5] (and of [EG2] respectively) and ispartially inspired by Popa's work [Popa].

%Keywords:  $C^*$-algebra, AH algebra, ideal property, Elliott invariant.

%\vspace{4mm}

% \noindent \emph{Keywords}: $C^*$-algebra, AH algebra, ideal property, Elliott invariant\\
% \noindent \emph{AMS subject classification}: Primary:  46L35, 46L80.

%\newpage

%\noindent{\underline {\noindent{ }\hspace {157mm}}}

%$$\noindent\textbf{Content}$$

\noindent\textbf{\S0. Introduction}
\vspace{-3mm}

\noindent\textbf{\S1. Notation and terminology}
\vspace{-3mm}

\noindent\textbf{\S2. The invariant}
\vspace{-3mm}

\noindent\textbf{\S3. Shape equivalence}
\vspace{-3mm}

\noindent\textbf{\S4. The spectral distribution property and the decomposition theorems}
%\vspace{-3mm}

%\noindent\textbf{\S5. Additional decomposition theorems}
\vspace{-3mm}

\noindent\textbf{\S5. The existence theorem }
\vspace{-3mm}

\noindent\textbf{\S6. The uniqueness theorem }
\vspace{-3mm}

\noindent\textbf{\S7. Proof of the main theorem}\\
\vspace{-3mm}

%\noindent\textbf{\S9. Appendix A: Further discussion of the  invariant}
%\vspace{-3mm}

%\noindent\textbf{\S10. Appendix B:  ATAI algebras}

%\noindent{\underline {\noindent{ }\hspace {157mm}}}
%\vspace{3mm}

\noindent\textbf{\S0. Introduction}

%\vspace{3mm}

A $C^*$ algebra $A$ is called AH algebras if it is inductive limit of the $C^*$ algebras of form \\
$(A_ {n}=\bigoplus\limits_{i=1}\limits^{t_{n}}P_{n,i}M_{[n,i]}(C(X_{n,i}))P_{n,i}, \phi_{n,m})$, where $X_{n,i}$ are compact metric spaces, $t_{n}$ and $[n,i]$ are
positive integers, and $P_{n,i}\in M_{[n,i]}(C(X_{n,i}))$ are projections.  More general, if  $A_n\subset \bigoplus\limits_{i=1}\limits^{t_{n}}M_{[n,i]}(C(X_{n,i}))$, then $A$ is called an ASH algebra.  It seems AH algebras and ASH algebras are quite special class of $C^*$, but these classes of $C^*$-algebras play important role in the Elliott classification program. It is a conjecture that all simple separable stable finite nuclear $C^*$-algebras are ASH algebras. Recent work of Gong-Lin-Niu and Elliott-Gong-Lin-Niu  have verified the conjecture for all simple separable $C^*$-algebras of finite decomposition rank with UCT. Combine with the work of Tikusis-White-Winter, one knows that the conjecture is also true for all simple separable stably finite $C^*$-algebras of finite nuclear dimension. The successful classification has been obtained for all simple separable $C^*$-algebras of finite nuclear dimension which satisfies universal coefficient theorem (see [KP] [P] [GLN1-2], [EGLN1-2] and [TWW]). So now it is timed to advance the classification project to non simple $C^*$ algebras.

Classification theorems have been obtained for $AH$ algebras---the inductive limits of cut downs of matrix algebras over compact metric spaces by projections---and $AD$ algebras---the inductive limits of Elliott dimension drop interval algebras in two special cases long ago:

 $1$. Real rank zero case: all such $AH$ algebras
with no dimension growth and such $AD$ algebras (See [Ell1], [EGLP], [EG1-2], [EGS], [D], [G1-4], [Ei], and [DG]) and

$2$. Simple case: all
such $AH$ algebras with no dimension growth (which includes all simple $AD$ algebras by [EGJS]) (See [Ell2-3], [NT], [Thm2-3], [Li1-4], [G5], and
[EGL1]).

In [EGL1], the authors pointed out two important possible next steps after the complete classification of simple $AH$ algebras (with no dimension
growth). One of these is the classification of simple $ASH$ algebras---the simple inductive limits of sub-homogeneous algebras (with no dimension growth), which has been done in [EGLN1]. The
other is to generalize and  unify the above-mentioned classification theorems for simple $AH$ algebras and real rank zero $AH$ algebras by classifying
$AH$ algebras with the ideal property. In this article, we have achieved the second of these goals  by classifying all $AH$ algebras (of no dimension growth) and
$AD$ algebras, provided the algebras have the ideal property, as explained below.

As in [EG2], let $T_{\uppercase\expandafter{\romannumeral2},k}$ be the $2$-dimensional connected simplicial complex with $H^{1}(T_{\uppercase\expandafter{\romannumeral2},k})=0$ and $H^{2}(T_{\uppercase\expandafter{\romannumeral2},k})=\mathbb{Z}/k\mathbb{Z}$, and
let $I_{k}$ be the subalgebra of $M_{k}(C[0,1])=C([0,1],M_{k}(\mathbb{C}))$ consisting of all functions $f$ with the properties $f(0)\in \mathbb{C}\cdot \textbf{1}_{k}$ and $f(1)\in \mathbb{C}\cdot \textbf{1}_{k}$
(this algebra is called an Elliott dimension drop interval algebra). Denote  $\mathcal{HD}$ the class of algebras consisting of direct sums of building
blocks of the forms $M_{l}(I_{k})$ and $PM_{n}(C(X))P$, with $X$ being one of the spaces $\{pt\}$, $[0,1]$, $S^{1}$, and $T_{\uppercase\expandafter{\romannumeral2},k}$, and with $P\in M _{n}(C(X))$ being a
projection. (In [DG], this class is denoted by $SH(2)$, and in [Jiang1], this class is denoted by $\mathcal{B}$). We will call a $C^{*}$-algebra an $A\mathcal{HD}$
algebra, if it is an inductive limit of algebras in $\mathcal{HD}$. In this article, we will classify all $A\mathcal{HD}$ algebras with the ideal
property. Evidently, all simple $AH$ algebras and all real rank zero $C^{*}$-algebras have the ideal property (see Abstract). In [GJLP1-2], [Li4],  and [Jiang2], it is proved that
all $AH$ algebras with the ideal property of no dimension growth are inductive limits of algebras in the class $\mathcal{HD}$---that is, they are $A\mathcal{HD}$ algebras. By this reduction theorem, then we obtained the classification of $C^{*}$-algebras with the ideal property which are either
$AH$ algebras of no dimension growth or $AD$ algebras.

The ideal property is a property of structural interest for a $C^*$-algebra. Many interesting and important $C^*$-algebras have the ideal property. It was proved by Cuntz-Echterhoff-Li that semigroup $C^*$-algebras of $ax + b$-semigroups over Dedekind domains have the ideal property ([CEL]). As pointed out in [GJLP1], there are many $C^{*}$-algebras naturally arising from $C^{*}$-dynamical systems which have the ideal property. For example, if
$(A,G,\alpha)$ is a $C^{*}$-dynamical system such that $G$ is a discrete amenable group and the action of $G$ on the spectrum $\widehat{A}$ of $A$ is essentially free, that is, for
every $G$-invariant closed subset $F\subset\widehat{A}$, the subset $\{x\in F: gx\neq x$ for all $g\in G\backslash \{1\}\}$ is dense in $F$, then the cross product $C^*$-algebra $A\rtimes_{\alpha}G$
has the ideal property provided that $A$ has the ideal property (see [Sierakowski, Theorem 1.16]). Many of these $C^{*}$-algebras are neither simple nor of real
rank zero. Even if we assume, in the above class of examples, that $A$ is the commutative $C^{*}$-algebra $C(X)$ with $dim(X)=0$ and $G=\mathbb{Z}$, it is not
known whether $A\rtimes_{\alpha}G$ is of real rank zero (it is known that $A\rtimes_{\alpha}G$ is not simple if $\alpha$ is not minimal). However it
follows from the above result or from [Pa2] that $C(X)\rtimes_{\alpha}\mathbb{Z}$ has the ideal property. Hence it is important and natural to extend
the classification of simple $C^{*}$-algebras and the classification of real rank zero $C^{*}$-algebras to $C^{*}$-algebras with the ideal property.
It is proved by Elliott-Niu [EN] that the simple $C^*$-algebra $C(X)\rtimes_{\alpha}\Z$ is an $ASH$ algebra   and are classifiable, provided that $\alpha$ is minimum (this is equivalent to the condition that the cross product is simple)
and mean dimension of the dynamical system is zero (the special case that $dim(X)<\infty$  is due to Lin [Lin3]).
The proof of this important theorem (see [EN]) involves many techniques, in particular, the decomposition theorems,
from the classification of simple $AH$ and $ASH$ algebras. In this article we generalize some of these techniques to the case of $C^*$-algebras with ideal property. The classification of the extended class of $C^*$-algebras---the $C^*$ algebras with ideal property (e.g., the result in this paper) combining with the techniques in the proof may
be very useful to the study of $AH$ or $ASH$ decompositions of the possible non-simple cross product $C^*$-algebras $C(X)\rtimes_{\alpha}\Z$, and in turn, in  the above mentioned case of an essential free action  $\alpha$ on a zero dimensional  space $X$, the real rank zero property of the cross product $C^*$-algebras $C(X)\rtimes_{\alpha}\Z$ could be deduced from its $AH$ decomposition.
Let us point out that C. Pasnicu studied the $C^{*}$-algebras with the ideal property intensively [Pa1--5].

It is proved in [CJL]  that all inductive limits of direct sums of simple $TAI$ algebras (with $UCT$),  which are called $ATAI$ algebras (these strictly include the $ATAF$ algebras studied by Fang in [Fa]), are in the above class---that is, they are $A\mathcal{HD}$ algebras with the ideal property. In [Jiang1], the second named author classified $ATAI$ algebras $A$ under the  extra condition that $K_{1}(A)$ is a torsion group.  In this classification, Jiang used the scaled ordered total $K$-group (from [DG]) and the tracial state
spaces $T(pAp)$ of cut-down algebras $pAp$, with certain compatibility conditions (from [Stev] and [Ji-Jiang]) as the invariant---we will call it $Inv^0(A)$ (see 2.18 of [Jiang1]). In our
classification, we will use the above invariant  $Inv^0(A)$ together  with the extra new invariant $U(pAp)/\overline{DU(pAp)}$ to deal with the torsion-free part of $K_1(A)$ (see Abstract---we will use a simplified version
of this group), with certain compatibility conditions---we will call this $Inv(A)$. In [GJL], we will  prove that our invariant reduces to Jiang's invariant in the case that $K_{1}(A)$ is a torsion group (see Proposition 2.38 in [GJL]). That is, Jiang's classification in [Jiang1] is a special case of our classification theorem (see [CJL]). Let us point out that
$U(pAp)/\overline{DU(pAp)}$ is completely determined, non-canonically,   by $T(pAp)$ and $K_{1}(pAp)$. But in [GJL], we will present an example to show that the compatibility of
$U(pAp)/\overline{DU(pAp)}$ (between different projections $p<q$) is not completely decided by the weaker invariant $Inv^0(A)$ used in Jiang's paper. That is, we will present two unital ${\cal Z}$-stable (where ${\cal Z}$ is the Jiang-Su algebra of [JS]) A$\T$ algebras $A$ and $B$ with the ideal property such that $Inv^0(A)\cong Inv^0(B)$, but $Inv(A)\not\cong Inv(B)$ (see Proposition 3.4 and 3.10 in [GJL]). So the new
invariant is indispensable. Furthermore, for these $C^*$-algebras $A$ and $B$, we have that $Cu(A)\cong Cu(B)$ and $Cu(A\otimes C(S^1))\cong Cu(B\otimes C(S^1))$. That is, our new invariant can not be be detected by Cuntz semigroup.

%%%%%%%%

It is worthwhile to compare the present introduction of the new invariant $U(pAp)/\overline{DU(pAp)}$ with its compatibility data to the introduction of the invariant consisting of the order structure on total K-theory for the case of non-simple real rank zero $C^{*}$-algebras (in [DG] and [DL]).  In the simple real rank zero case, it was sufficient to look at the K-theory alone including its order structure (introduced in [Ell1]).  In the non-simple real rank zero case, in [DG], it was necessary to include the total K-theory with its order structure. As a matter of fact, the total K-theory by itself, is completely decided by the K-theory alone; it is the order structure of the total K-theory which reflects certain compatibility information associated with  the ideals.  In the general simple case, for the isomorphic invariant, it is enough to look at the ordering of the K-theory,  now together with the simplex of traces (as that is no longer determined). In the general non-simple case with the ideal property, with general $K_1$-group, both the total K-theory (with its order) and also the Hausdorffifized  algebraic $K_1$-group
$U(pAp)/\overline{DU(pAp)}$ of the cut-down algebra $pAp$ with compatibility are needed in the  isomorphism theorem.  As pointed out in [NT] (see [Ell3] also) the Hausdorffifized  algebraic $K_1$-group $U(A)/\overline{DU(A)}$ is completely decided by the Elliott invariant consisting of $K_*(A)$, AffT$A$, and the map $\rho: ~ K_0(A) \to \mbox{AffT}A$---it is the compatibility of $U(pAp)/\overline{DU(pAp)}$ and $U(qAq)/\overline{DU(qAq)}$ (for $[p]<[q]$) which makes it a new invariant, analogous to the order structure of the total K-theory in the  non-simple real rank zero case.  Let us point out that, even in the simple case, for uniqueness theorems for homomorphisms, both the total K-theory  (which was reflected by KK  classes already in the real rank zero case in [EG2] and [G2]), and the Hausdorffifized  algebraic $K_1$-group $U(A)/\overline{DU(A)}$ (introduced in [NT], see also [Ell3])  are needed.

For the simple case, the additional invariants of total K-theory or of Hausdorffifized  algebraic $K_1$, are completely decided by the other parts of the invariant, non-canonically.  By ``non-canonically" we mean that an isomorphism between K-theories may allow  different choices of isomorphisms between the total K-theories, while an isomorphisms between Elliott invariants may allow different choices of isomorphisms between the Hausdorffifized  algebraic $K_1$-groups.  These choices make the compatibility become a problem for the non-simple case (either in the case of real rank zero, or in the case with ideal property), which does not appear  in the simple case.

Robert proved in [Rob] that an inductive limit of 1-dimensional non commutative CW complexes can be classified by a modified version $Cu^{\sim}(A)$ of Cuntz semigroup (note that if $A$ is unital then $Cu^{\sim}(A)$ is completed decided by $Cu(A)$). On the other, as we already mentioned, our new invariant can not be recover from the Cuntz semigroups of $A$ and $A\otimes C(S^1)$.
We believe this new addition of the invariant (together with Cuntz semi group [Cu]) may be useful for the classification of the general non simple inductive limits of 1-dimensional non commutative CW complexes, beyond the case with ideal property.

In [W], Wang proves that for the $C^*$algebras with the  ideal property, the invariant $Inv^0(A)$ is equivalent to an extended version of the Elliott invariant. It seems there is no natural  way to extend the Elliott invariant to recover the invariant $Inv(A)$.  Let us also pointed out an interesting phenomenon: if the   $C^*$-algebra $A$  (with the ideal property) has only  finitely  many ideals, then $Inv(A)$ will be completely determined by $Inv^0(A)$---the proof will be presented in a forthcoming paper joint with K. Wang.

%%%%%%%%%%

Beyond the techniques used in [EG2], [DG], [Thm1-4], [NT], [G5], [EGL1-2], [Li4], [Ji-Jiang], [Jiang1-2],  and [GJLP1-2], some important new
techniques are introduced in this paper. Some of them are introduced to deal with the new invariant both in the existence theorem and in the uniqueness theorem. Another
major difficulty is as follows. The local uniqueness theorem requires the homomorphisms involved to satisfy a certain spectral distribution property, called the  $sdp$
property (more specifically, $sdp(\eta,\delta)$ property for certain $\eta>0$ and $\delta>0$ introduced in [G5] and [EGL1]). This property automatically
 holds for the homomorphisms $\phi_{n,m}$ (provided that $m$ is large enough) giving rise to  a simple inductive limit procedure. But for the case of general inductive limit $C^{*}$-algebras
with the ideal property, considered in our paper, to obtain this $sdp$ property, we must pass to certain good quotient algebras which   correspond to
simplicial sub-complexes of the  original spaces; a uniform uniqueness theorem, that does not depend on the choice of simplicial sub-complexes involved, is
required. For the case of an interval, whose simplicial sub-complexes are finite unions of subintervals and points, such a uniform uniqueness theorem is proved
in [Ji-Jiang] (see [Li2] and [Ell2] also). But for the general case, no such theorem is  true. In this paper, we carefully manage to avoid the (lack of) uniqueness theorem for the general case
involving arbitrary finite subsets of $M_{n}(C(T_{\uppercase\expandafter{\romannumeral2},k}))$ (or $M_{l}(I_{k})$), and, instead,  only use  one involving special finite subsets of $M_{n}(C(T_{\uppercase\expandafter{\romannumeral2},k}))$ (or $M_{l}(I_{k})$),
those  are approximately constant to within a small number $\varepsilon$---for which we can prove the required uniform uniqueness theorem.
Of course, we will use the uniform uniqueness theorem for intervals in [Ji-Jiang].
Also,  some technical
problems should be dealt with involving the building blocks $M_{l}(I_{k})$. The decomposition theorems  between
a building block of this kind and a homogeneous building block are proved in [JLW], which will be used here.

To summarize, we have obtained a rather general classification, for inductive limits of direct sum of building blocks which are either homogeneous with no dimension growth or dimension drop interval algebras, assuming that ideals are generated by projections (the ideal property).  This last assumption is needed so far even in the case of A$\T$ algebras. The range of the invariant has not been calculated yet, but might very well  include everything that arises from stably finite separable ${\cal Z}$-stable  $C^*$-algebras with ideal property, with $K_*$ to be Riesz group and the map $T(pAp) \to SK_0(pAp)$ (for any projection $p\in A$), takes the extreme points to extreme points (see [Vi]), where $SK_0(pAp)$ is  the space of homomorphisms $f: K_0(pAp) \to \R$ such that $f(K_0(pAp)_+)\subset \R_+$ and $f([p])=1$.

%\newpage

\vspace{3mm}

\noindent\textbf{\S1. Notation and terminology}

\vspace{-2.4mm}

In this section, we will introduce some notation and terminology.

\noindent\textbf{1.1.}~~Let $A$ and $B$ be two $C^{*}$-algebras. We use $Map(A,B)$ to denote the {\bf space of all completely positive  $*$-contractions} from $A$ to $B$. If
both $A$ and $B$ are unital, then $Map(A,B)_{1}$ will denote the subset of $Map(A,B)$ consisting of all such unital maps. By word ``map", we mean linear, completely
positive $*$-contraction between $C^{*}$-algebras, or else we shall mean continuous map between topological spaces, which one will be clear from the context.

\noindent\textbf{Definition 1.2.}~~ Let $G\subset A$ be a finite set and $\delta>0$. We shall say that $\phi\in Map(A,B)$ is {\bf $G-\delta$ multiplicative} if
$\|\phi(ab)-\phi(a)\phi(b)\|<\delta$ for all $a,b\in G$.

We also use $Map_{G-\delta}(A,B)$ to denote all $G-\delta$ multiplicative maps.

\noindent\textbf{1.3.}~~In the notation for an inductive limit system $\lim (A_{n},\phi_{n,m})$, we understand that $$\phi_{n,m}=\phi_{m-1,m}\circ\phi_{m-2,m-1}\circ\cdots\circ\phi_{n,n+1},$$
where all $\phi_{n,m}:A_{n}\rightarrow A_{m}$ are homomorphisms.

We shall assume that, for any summand $A^{i} _{n}$ in the direct sum $A_{n}=\bigoplus^{t_{n}}_{i=1}A^{i} _{n}$, necessarily,
$\phi_{n,n+1}(\textbf{1}_{A^{i} _{n}})\neq0$, since, otherwise, we could simply delete $A^{i} _{n}$ from $A_{n}$, without changing the limit algebra.

If $A_{n}=\bigoplus_{i}A^{i} _{n}$, $A_{m}=\bigoplus_{j}A^{j} _{m}$, we use $\phi^{i,j}_{n,m}$ to denote the partial map of $\phi_{n,m}$ from the
$i$-th block $A^{i} _{n}$ of $A_{n}$ to the $j$-th block $A^{j} _{m}$ of $A_{m}$. Also, we use $\phi_{n,m}^{-,j}$ to denote the partial
map of $\phi_{n,m}$ from $A_{n}$ to $A_{m}^{j}$. That is,  $\phi_{n,m}^{-,j}=\bigoplus\limits_{i}\phi_{n,m}^{i,j}=\pi_{j}\phi_{n,m}$,
where $\pi_{j}:A_{m}\rightarrow A_{m}^{j}$ is the canonical projection. Some times, we also use $\phi_{n,m}^{i,-}$ to denote $\phi_{n,m}|_{A_n^i}:~ A_n^i \to A_m$.

\noindent\textbf{1.4.}~~As mentioned in the introduction, we use $\mathcal{HD}$ to denote all $C^{*}$-algebras $C=\bigoplus C^{i}$, where each $C^{i}$ is of the forms
$M_{l}(I_{k})$ or $PM_{n}C(X)P$ with $X$ being one of the spaces $\{pt\}$, $[0,1]$, $[a,b]$ $(-1< a<b\leq 1)$, $S^{1}$, $T_{\uppercase\expandafter{\romannumeral2},k}$. (Here the interval $[a,b]$ is included  in this list of spaces for convenience.) Each block $C^{i}$ will be called a basic $\mathcal{HD}$ block or a basic building block. Note
that in this notation, a simplicial sub-complex of interval $[0,1]$ is a finite union $\bigcup\limits_{i}[a_{i},b_{i}]$ with $0\leq a_{i}\leq b_{i}\leq 1$ (including degenerated intervals which are single point sets), and a proper  sub-complex $X$ of $S^{1}$ (that is $\emptyset \neq X\subsetneqq S^{1}$), is
a finite union $X=\bigcup_{i}[a_{i},b_{i}]$, with $-1< a_{i}<b_{i}<a_{i+1}\leq 1$, where each $[a_{i},b_{i}]$ is identified with $\{e^{2\pi it},a_{i}\leq t\leq b_{i}\}$.
Note that if $X=\{pt\},[0,1]$ or $S^{1}$ (but not $T_{\uppercase\expandafter{\romannumeral2},k}$) then $PM_{n}C(X)P\cong M_{n_{1}}(C(X))$ with $n_{1}=rank(P)$, since all the projections on those spaces
are trivial. It is easy to see that for any projection $Q\in M_{l}(I_{k})$, the algebra $QM_{l}(I_{k})Q$ is isomorphic to a $C^{*}$-algebra of the form $M_{l_{1}}(I_{k})$ with $l_{1}\leq l$.

\noindent\textbf{1.5.}~~By $A\mathcal{HD}$ algebra, we mean the inductive limit of
$$A_{1}\xrightarrow{\phi_{1,2}}A_{2}\xrightarrow{\phi_{2,3}}A_{3}\longrightarrow\cdots\longrightarrow\cdots,$$
where $A_{n}\in \mathcal{HD}$ for each $n$. Obviously, in such an expression, we do not
need $[a,b]$ other than $[0,1]$.

For  an $A\mathcal{HD}$ inductive limit $A\!=\!\lim (A_n, \phi_{nm})$, we write $A_n\!=\!\oplus_{i=1}^{t_n}A_n^i$, where
$A_n^i=P_{n,i}M_{[n,i]}(C(X_{n,i}))P_{n,i}$ or $A_n^i=M_{[n,i]}(I_{k_{n,i}})$. For convenience, even for a block $A_n^i=M_{[n,i]}(I_{k_{n,i}})$, we still use  $X_{n,i}$ for $Sp(A_n^i)=[0,1]$---that is, $A_n^i$ is regarded as a homogeneous algebra or a  sub-homogeneous algebra over $X_{n,i}$.
%the spectrum $X_{n,i}$.

\noindent\textbf{1.6.}~~In [GJLP1-2], joint with Cornel Pasnicu, the authors proved the reduction theorem for $AH$ algebras with ideal property provided that the inductive limit
have no dimension growth. That is, if $A$ is an inductive limit of $A_{n}=\bigoplus A^{i}_{n}=\bigoplus P_{n,i}M_{[n,i]}C(X_{n,i})P_{n,i}$ with $sup_{n,i}dim(X_{n,i})<+\infty$, and if we further assume that $A$ has the ideal property, then $A$ can be rewritten as an  inductive limit of $B_{n}=\bigoplus B^{j}_{n}=\bigoplus Q_{n,j}M_{\{n,j\}}C(Y_{n,i})Q_{n,j}$, with $Y_{n,i}$ being one of $\{pt\}$, $[0,1]$, $S^{1}$, $T_{\uppercase\expandafter{\romannumeral2},k}$, $T_{\uppercase\expandafter{\romannumeral3},k}$, $S^{2}$. In turn, the
 second author proved in [Jiang2](also see [Li4]), that the above inductive limit can be rewritten as the inductive limit of direct sums of homogeneous algebras
over $\{pt\}$, $[0,1]$, $S^{1}$, $T_{\uppercase\expandafter{\romannumeral2},k}$ and $M_{l}(I_{k})$. Combining these two results, we know that all $AH$ algebras of no dimension growth with the ideal
property are $A\mathcal{HD}$ algebras. Let us point out that,  as proved in [DG], there are real rank zero $A\mathcal{HD}$ algebras which are not $AH$ algebras.

\noindent\textbf{1.7.}~~
(a)~~ we use $\sharp(.)$ to denote the cardinal number of a set. Very often, the sets under consideration will be sets with multiplicity, in which case we shall
also count multiplicity when we use the notation $\sharp$.\\
(b)~~ We shall use $a^{\sim k}$ to denote $\underbrace{a,a,\cdots a}\limits_{k}$. For example $\{a^{\sim3},b^{\sim2}\}=\{a,a,a,b,b\}$.\\
(c)~~ For any metric space $X$, any $x_{0}\in X$ and $c>0$, let $B_{c}(x_{0})\triangleq\{x\in X|d(x,x_{0})<c\}$,  the open ball with radius $c$ and center $x_{0}$.\\
(d)~~ Suppose that $A$ is a $C^{*}$-algebra, $B\subset A$ a subset (often a subalgebra), $F\subset A$ is a finite subset and $\varepsilon>0$. If for each element $f\in F$, there is an element $g\in B$ such
that $\|f-g\|<\varepsilon$, then we shall say that $F$ is approximately contained in $B$ to within $\varepsilon$, and denote this by $F\subset_{\varepsilon}B$.\\
(e)~~ Let $X$ be a compact metric space. For any $\delta>0$, a finite set $\{x_{1},x_{2},\cdots,x_{n}\}$ is said to be $\delta$-dense in $X$ if for
any $x\in X$, there is $x_{i}\in \{x_{1},x_{2},...,x_{n}\}$ such that $dist(x,x_{i})<\delta$.\\
(f)~~ we shall use $\bullet$ or $\bullet\bullet$ to denote any possible positive integers.\\
(g)~~ For any two projections $p,q\in A$, by $[p]\leq[q]$ we mean that $p$ is unitarily equivalent to a sub-projection of $q$. And we
use $p\sim q$ to denote that $p$ is unitarily equivalent to $q$.

\noindent\textbf{Definition 1.8.}([EG2], [DG]) ~~ Let $X$ be a compact connected space and let $P\in M_{N}(C(X))$ be a projection of rank $n$. The {\bf weak variation} of a finite set $F\subset PM_{N}C(X)P$ is defined by $$\omega(F)=\sup\limits_{\pi,\sigma}\inf\limits_{u\in U(n)}\|u\pi(a)u^{*}-\sigma(a)\|,$$ where
$\pi,\sigma$ run through the set of irreducible representations of $PM_{N}C(X)P$ into $M_{n}(\mathbb{C})$.

For $F\subset M_{r}(I_{k})$, we define $\omega(F)=\omega(\imath(F))$,
where $\imath: M_{r}(I_{k})\longrightarrow M_{rk}(C[0,1])$ is the canonical embedding and $\imath(F)$ is regarded as a finite subset of $M_{rk}(C[0,1])$. Let
$A$ be  a basic $ \mathcal{HD}$ block, a finite set $F\subset A$ is said to be {\bf weakly approximately constant to within} $\varepsilon$ if $\omega(F)<\varepsilon$.

Evidently, the following lemma is true. \\
\noindent\textbf{Lemma 1.9.}~~Let $A,B$ be basic $ \mathcal{HD}$ blocks. and $\phi: A\longrightarrow B$, a homomorphism. For any finite
set $F\subset A$, we have $\omega(\phi(F))\leq\omega(F)$. Consequently,  if $F$ is weakly approximately constant to within $\varepsilon$, then so is $\phi(F)$.

The following dilation lemma is proved in [EG2].\\
 \noindent\textbf{Proposition 1.10.}~~(see Lemma 2.13 of [EG2]) Let $X$ and $Y$ be any connected finite $CW$ complexes. If $\phi: QM_{k}(C(X))Q\longrightarrow
 PM_{n}(C(Y))P$ is a unital homomorphism, then there are:  an $n_{1}\in \mathbb{Z}_{+}$, a projection $P_{1}\in M_{n_{1}}(C(Y))$, and a unital
 homomorphism $\widetilde{\phi}: M_{k}(C(X))\longrightarrow P_{1}(M_{n_{1}}C(Y))P_{1}$ with the property that $QM_{k}(C(X))Q$ and $PM_{n}(C(Y))P$ can be identified
as corner subalgebras of $M_{k}(C(X))$ and $P_{1}M_{n_{1}}(C(Y))P_{1}$,  respectively (i,e $Q$ and $P$ can be considered to be sub-projections of $\textbf{1}_{k}$ and $P_{1}$
respectively),  and,  in such a way that $\phi$ is the restriction of $\widetilde{\phi}$.

If $\phi_{t}: QM_{k}(C(X))Q\longrightarrow PM_{n}(C(Y))P$ is a path of unital homomorphisms, then there exist  $P_{1}M_{n_{1}}(C(Y))P_{1}$ (as above) and a path
of unital homomorphisms $\widetilde{\phi}_{t}: M_{k}(C(X))\longrightarrow P_{1}M_{n_{1}}(C(Y))P_{n_{1}}$ such that $QM_{k}(C(X))Q$ and $PM_{n}(C(Y))P$ are corner subalgebras of $M_{k}(C(X))$ and $P_{1}M_{n_{1}}(C(Y))P_{1}$,  respectively, and $\phi_{t}$ is the restriction of $\widetilde{\phi}_{t}$.

\noindent\textbf{Lemma 1.11.}~~Let $X$ be a connected finite $CW$ complex. If $\phi: QM_{k}(C(X))Q\longrightarrow M_{l}(I_{d})$ is a unital homomorphism, then
there is an $l_{1}\geq l$ and a unital homomorphism $\widetilde{\phi}: M_{k}(C(X))\longrightarrow M_{l_{1}}(I_{d})$ with property that $QM_{k}(C(X))Q$ and
$M_{l}(I_{d})$ can be idenfied as corner subalgebras of $M_{k}(C(X))$ and $M_{l_{1}}(I_{d})$ respectively (i.e. $Q$ and $\textbf{1}_{M_{l}(I_{d})}$ can be considered
to be subprojections of $\textbf{1}_{k}$ and of $\textbf{1}_{M_{l_{1}}(I_{d})}$ respectively) and furthermore, in such a way that $\phi$ is the restriction of $\widetilde{\phi}$.

If $\phi_{t}: QM_{k}(C(X))Q\longrightarrow M_{l}(I_{d})$ $(0\leq t\leq 1)$ is a path of unital homomorphisms, then there are $M_{l_{1}}(I_{d})$ (as above) and a path of unital homomorphisms $\widetilde{\phi}_{t}: M_{k}(C(X))\longrightarrow M_{l_{1}}(I_{d})$ such that $QM_{k}(C(X))Q$ and $M_{l}(I_{d})$ are corner sub-algebras of $M_{k}(C(X))$ and $M_{l_{1}}(I_{d})$, respectively,  and $\phi_{t}$ is the restriction of $\widetilde{\phi}_{t}$.
\begin{proof}
The proof of this lemma is the same as the above proposition---Lemma 2.13 of [EG2]. See 2.12 of [EG2]. One only needs to notice that for any projection $P\in M_{n}(I_{d})\otimes(C[0,1])$, we can identify $P(M_{n}(I_{d})\otimes C[0,1])P$ with $M_{l_{1}}(I_{d})\otimes C[0,1]$.\\
\end{proof}

\noindent\textbf{1.12.}~~If $A\in \mathcal{HD}$ is a basic building block, then $K_{0}(A)=\mathbb{Z}$ (when $A$ is not of type $T_{\uppercase\expandafter{\romannumeral2},k}$), or $\mathbb{Z}\oplus\mathbb{Z}/k\mathbb{Z}$ (when $A=PM_{n}(T_{\uppercase\expandafter{\romannumeral2},k})P)$. So there is a natural map $rank$: $K_{0}(A)\longrightarrow\mathbb{Z}$. If $A=PM_{n}(C(X))P$, then this map is induced by any irreducible
representation of $A$ as $(\pi_{x})_*: K_{0}(A)\longrightarrow K_{0}(M_{rank(P)}(\mathbb{C}))=\mathbb{Z}$, where $\pi_{x}: PM_{n}(C(X))P\longrightarrow M_{rank(P)}(\mathbb{C})$ is a irreducible representation corresponding to $x\in X$. If $A=M_{l}(I_{k})$, then this  rank map is the map induced on $K_{0}$ by the evaluation map at 0 or 1. Let $\widetilde{K}_{0}(A)$ denote the kernel of rank map. Then $\widetilde{K}_{0}(A)=0$ if $A$ is $M_{r}(C(X))$ for $X$ being one of $\{pt\}$, $[0,1]$, $S^{1}$ or $A=M_{l}(I_{k})$ and $\widetilde{K}_{0}(A)=\mathbb{Z}/k\mathbb{Z}$ if $A=PM_{r}C(T_{\uppercase\expandafter{\romannumeral2},k})P$ .

Let $A,B\in\mathcal{HD}$ be basic building blocks and let $\alpha\in KK(A,B)$. Then $\alpha$ induces a map $\alpha_{*}$ on $K_{0}$-group. The $KK$ element $\alpha$ is called $m$-{\bf large} if $rank$ $\alpha_{*}(\textbf{1}_{A})\geq m$ $rank$$(\textbf{1}_{A})$, where $m>0$ is a positive integer.\\
(a) A unital homomorphism $PM_{\bullet}(C(X))P\longrightarrow QM_{\bullet}(C(Y))Q$ is $m$-large if and only if \\$rank(Q)/rank(P)\geq m.$\\
(b) A unital homomorphism $M_{l}(I_{k})\longrightarrow QM_{\bullet}(C(Y))Q$ is $m$-large if and only if ${rank(Q)}/{l}\geq m.$
\\
(c) A unital homomorphism $M_{l_{1}}(I_{k_{1}})\longrightarrow M_{l_{2}}(I_{k_{2}})$ is $m$-large if and only if $l_{2}/l_{1}\geq m.$\\
(d) A unital homomorphism $PM_{\bullet}(C(X))P\longrightarrow M_{l_{2}}(I_{k_{2}})$ is $m$-large if and only if ${l_{2}}/{rank(P)}\geq m.$

\noindent\textbf{1.13.}~~Let $A\in\mathcal{HD}$, and $p,q\in A$ be two projections. If $[p]=[q]\in K_{0}(A)$, then there is a unitary $u\in A$ such that  $p=uqu^{*}$. This fact is evidently true for the case  $A=M_{l}(I_{k})$,  and is also true for the case  $A=PM_{n}(C(X))P$, where $X$ is an arbitrary  3-dimensional $CW$ complex , by  3.26 of [EG2].

\noindent\textbf{1.14.}~~Let $X$ be a compact space and $\psi: C(X)\longrightarrow PM_{k_{1}}(C(Y))P$ ( $rank(P)=k$) be a unital homomorphism. For any point $y\in Y$, there are $k$ mutually orthogonal rank $1$  projections $p_{1},p_{2},\cdots,p_{k}$ with $\sum\limits_{i=1}\limits^{k}p_{i}=P(y)$ and $\{x_{1}(y),x_{2}(y),\cdots,x_{k}(y)\}\subset X$ (may be repeat) such that $\psi(f)(y)=\sum\limits_{i=1}\limits^{k}f(x_{i}(y))p_{i}, \forall f\in C(X).$ \\
We  denote the set $\{x_{1}(y),x_{2}(y),\cdots,x_{k}(y)\}$ (counting multiplicities), by $Sp\psi_{y}$.We shall call $Sp\psi_{y}$ {\bf the spectrum of } $\psi$
{\bf at the point} $y$.

\noindent\textbf{1.15.}~~For any $f\in I_{k}\subset M_{k}(C[0,1])=C([0,1],M_{k}(\mathbb{C}))$ as in [EGS, 3.2], let function $\underline{f}: [0,1]\longrightarrow \mathbb{C}\sqcup M_{k}(\mathbb{C})$ (disjoint union) be defined by
$$\underline{f}(t)=\begin{dcases}
 \lambda, & if\; t=0\; and\; f(0)=\lambda\textbf{1}_{k}\\
  \mu, & if\; t=1\; and \;f(1)=\mu\textbf{1}_{k}\\
  f(t), & if\; 0<t<1~~~~~~~~~~~~~~~~.
   \end{dcases}$$
That is, $\underline{f}(t)$ is the value of irreducible representation of $f$ corresponding to the point  $t$. Similarly,  for $f\in M_{l}(I_{k})$, we can
define $\underline{f}: [0,1]\longrightarrow M_{l}(\mathbb{C})\sqcup M_{lk}(\mathbb{C})$, by
$$\underline{f}(t)= \begin{cases}
  a, & if\;t=0\;and\;f(0)=a\otimes\textbf{1}_{k} \\
   b, & if\;t=1\;and\;f(1)=b\otimes\textbf{1}_{k} \\
    f(t),& if\;0<t<1~~~~~~~~~~~~~~~~~~.
   \end{cases}$$

\noindent\textbf{1.16.}~~Suppose that $\phi: I_{k}\longrightarrow PM_{n}(C(Y))P$ is  a unital homomorphism. Let $r=rank(P)$. For each $y\in Y$, there are $t_{1},t_{2},\cdots,t_{m}\in[0,1]$ and a unitary
$u\in M_{n}(\mathbb{C})$ such that
$$P(y)=u\left(
          \begin{array}{cc}
            \textbf{1}_{rank(P)} & 0 \\
            0 & 0 \\
          \end{array}
        \right)u^{*}$$
and$$ (*)~~~~~\;\;\;\;\phi(f)(y)=u\left(
                               \begin{array}{ccccc}
                                 \underline{f}(t_{1}) &  & & & \\
                                  & \underline{f}(t_{2}) &  &  &  \\
                                  &  & \ddots &  &  \\
                                  &  &  & \underline{f}(t_{m}) &  \\
                                  &  &  &  &  {\bf 0}_{n-r} \\
                               \end{array}
                             \right)
u^{*}\in P(y)M_{n}(\mathbb{C})P(y)\subset M_{n}(\mathbb{C})$$
for all $f\in I_{k}$.

\noindent\textbf{1.17.}~~Let $\phi$ be the homomorphism  defined by  $(*)$ in 1.16 with $t_{1},t_{2},\cdots,t_{m}$ as appeared in the diagonal of the matrix. We define the set $Sp\phi_{y}$ to be the points $t_{1},t_{2},\cdots,t_{m}$ with possible
fraction multiplicity. If $t_{i}=0$ or $1$, we will assume that the multiplicity of $t_{i}$ is  $\frac{1}{k}$; if $0<t_{i}<1$, we will assume that the multiplicity of  $t_{i}$ is $1$. For example if we assume $$t_{1}=t_{2}=t_{3}=0<t_{4}\leq t_{5}\leq \cdots \leq t_{m-2}<1=t_{m-1}=t_{m},$$ then
$Sp\phi_{y}=\{0^{\sim\frac{1}{k}},0^{\sim\frac{1}{k}},0^{\sim\frac{1}{k}},t_{4},t_{5},\cdots,t_{m-2},1^{\sim\frac{1}{k}},1^{\sim\frac{1}{k}}\}$, which  can also be written
as $$Sp\phi_{y}=\{0^{\sim\frac{3}{k}},t_{4},t_{5},\\ \cdots,t_{m-2},1^{\sim\frac{2}{k}}\}.$$ Here we emphasize that, for $t\in (0,1)$, we do not allow the multiplicity of $t$ to be non-integral. Also for 0 or 1, the multiplicity must be multiple of $\frac{1}{k}$ (other fraction numbers are not allowed).

Let $\psi: C[0,1]\longrightarrow PM_{n}(C(Y))P$ be defined by the following composition $$\psi:~C[0,1]\hookrightarrow I_{k}\xrightarrow {\phi}PM_{n}(C(X))P,$$ where
the first map is the canonical inclusion. Then we have $Sp\psi_{y}=\{Sp\phi_{y}\}^{\sim k}$---that is, for each element $t\in(0,1)$, its multiplicity in $Sp\psi_{y}$ is exactly $k$ times of the multiplicity in $\phi_{y}$.

\noindent\textbf{1.18.}~~Let $A=M_{l}(I_{k})$. Then every  point $t\in(0,1)$ corresponds to an irreducible representation $\pi_{t}$,  defined by $\pi_{t}(f)=f(t)$. The representations  $\pi_{0}$ and $\pi_{1}$ defined by $$\pi_{0}=f(0)~~~~~~\mbox{ and}~~~~~~~\pi_{1}=f(1)$$ are no longer irreducible. We use $\underline{0}$ and $\underline{1}$ to denote the corresponding points for the irreducible representations. That is, $$\pi_{\underline{0}}(f)=\underline{f}(0),~~~~~~~\mbox {and}~~~~~~~~\pi_{\underline{1}}(f)=\underline{f}(1).$$ Or we can also write $\underline{f}(0)\triangleq f(\underline{0})$ and $\underline{f}(1)\triangleq f(\underline{1})$. Then the equation $(*)$ could be written as
$$\phi(f)(y)=u\left(
                \begin{array}{ccccc}
                  f(t_{1}) &  &  &  &  \\
                   & f(t_{2}) &  &  &  \\
                  &  & \ddots &  &  \\
                   & &  & f(t_{m}) &  \\
                   &  &  &  & {\bf 0}_{n-r} \\
                \end{array}
              \right)u^{*},$$
where some of $t_{i}$ may be $\underline{0}$ or $\underline{1}$. In this notation,  $f(0)$ is is equal to  diag$(\underbrace{f(\underline{0}),f(\underline{0}),\ddots,f(\underline{0})}\limits_{k})$ up to unitary equivalence.

% Use notation in 1.15$a$, we can also write $0^{\sim\frac{1}{k}}$ as $\underline{0}$, then $$\{0^{\sim\frac{1}{k}},0^{\sim\frac{1}{k}},0^{\sim\frac{1}{k}},t_{4},t_{5},\cdots,t_{m-2},1^{\sim\frac{1}{k}},1^{\sim\frac{1}{k}}\}\\=\{\underline{0},\underline{0},\underline{0},
%t_{4},t_{5},\cdots,t_{m-2},\underline{1},\underline{1}\}.$$

Under this notation, we can also write $0^{\sim\frac{1}{k}}$ as $\underline{0}$. Then the example of $Sp\phi_{y}$ in 1.17 can be  written as
$$Sp\phi_{y}=\{0^{\sim\frac{1}{k}},0^{\sim\frac{1}{k}},0^{\sim\frac{1}{k}},t_{4},t_{5},\cdots,t_{m-2},1^{\sim\frac{1}{k}},1^{\sim\frac{1}{k}}\}\\=\{\underline{0},\underline{0},\underline{0},
t_{4},t_{5},\cdots,t_{m-2},\underline{1},\underline{1}\}.$$

\noindent\textbf{1.19.}~~For a homomorphism $\phi: A\longrightarrow M_{n}(I_{k})$, where $A=I_{k}$ or $C(X)$, and for any $t\in[0,1]$, define $Sp\phi_{t}=Sp\psi_{t}$, where
$\psi$ is defined by the composition $$\psi: A\xrightarrow {\phi} M_{n}(I_{l})\rightarrow M_{nl}(C[0,1]).$$
Also $Sp\phi_{\underline{0}}=Sp(\pi_{\underline{0}}\circ \phi)$. Hence, $Sp\phi_{0}=\{Sp\phi_{\underline{0}}\}^{\sim k}.$

\noindent\textbf{1.20.}~~Let $\phi: M_{n}(A)\longrightarrow B$ be a unital homomorphism. It is well known (see 1.34 and 2.6 of [EG2]) that there is an identification of $B$
with $(\phi(e_{11})B\phi(e_{11}))\otimes M_{n}(\mathbb{C})$ such that $$\phi=\phi_{1}\otimes id_{n}: M_{n}(A)=A\otimes M_{n}(\mathbb{C})\longrightarrow (\phi(e_{11})B\phi(e_{11}))\otimes M_{n}(\mathbb{C})=B,$$ where $e_{11}$ is the matrix unit of upper left corner of $M_{n}(A)$ and $\phi_{1}=\phi|_{e_{11}M_{n}(A)e_{11}}: A\longrightarrow \phi(e_{11})B\phi(e_{11})$.

If we further assume that $A=I_{k}$ or $C(X)$ (with $X$ being a connected $CW$ complex) and $B$ is either $QM_{n}(C(Y))Q$ or $M_{l}(I_{k_{1}}),$ then for any
$y\in SpB$, define $Sp\phi_{y}\triangleq Sp(\phi_{1})_{y}.$
Here,  we use the standard notation that if $B=PM_{m}(C(Y))P$ then $SpB=Y$; and if $B=M_{l}(I_{k})$, then $Sp(B)=[0,1].$

If $\phi: PM_{n}(C(X))P\longrightarrow B$ is a homomorphism with $B$ being either $QM_{n}(C(Y))Q$ or $M_{l}(I_{k})$, then by Proposition 1.10 and Lemma 1.11, there is a dilation $\widetilde{\phi}$ defined on $M_{n}(C(X))$. We define $Sp\phi_{y}\triangleq Sp\widetilde{\phi}_{y}.$

\noindent\textbf{1.21.}~~Suppose that both  $A$ and $B$ are  of form $PM_{n}(C(X))P$ (with $X$ path connected) and $\phi: A\longrightarrow B$ be a unital homomorphism, we define $SPV(\phi)$
({\bf spectral variation} of $\phi$) by $$SPV(\phi)\triangleq \sup\limits_{y_{1},y_{2}\in Sp(B)}dist(Sp\phi_{y_{1}},Sp\phi_{y_{2}}),$$ where for two sets
$\{x_{1},x_{2},\cdots,x_{n}\}$, $\{x'_{1},x'_{2},\cdots,x'_{n}\}\subset Sp(A)$, counting multiplicity, the distance $$dist(\{x_{1},x_{2},\cdots,x_{n}\},\{x'_{1},x'_{2},\cdots,x'_{n}\})$$ is defined by $\min\limits_{\sigma}\max(dist(x_{i},x'_{\sigma(i)})$ with
$\sigma$ running over all permutations of $\{1,2,\cdots,n\}$.

\noindent\textbf{1.22.}~~Let $A=M_{l}(I_{k})$, $B=PM_{n}(C(X))P$ with $X$ connected, and $\phi: A\longrightarrow B$ be a unital homomorphism. Then there are two nonnegative
integers $k_{0},k_{1}\leq k-1$ such that for any $y\in X$ $$Sp\phi_{y}=\{0^{\thicksim\frac{k_{0}}{k}},1^{\thicksim\frac{k_{1}}{k}},t_{1},t_{2},\cdots,t_{m}\},$$
where $m=(\frac{rank(P)}{l}-k_{0}-k_{1})\cdot \frac{1}{k}$, and $\{t_{1},t_{2},\cdots,t_{m}\}\subset[0,1]$ are $m$ points depending on $y$ (some of them
may be 0 or 1). If $y,y'\in X$, with $$Sp\phi_{y}=\{0^{\thicksim\frac{k_{0}}{k}},1^{\thicksim\frac{k_{1}}{k}},t_{1},t_{2},\cdots,t_{m}\},
~~~~\mbox{and}~~~Sp\phi_{y'}=\{0^{\thicksim\frac{k_{0}}{k}},1^{\thicksim\frac{k_{1}}{k}},t_{1}',t_{2}',\cdots,t_{m}'\},$$ then
$dist(Sp\phi_{y},Sp\phi_{y^{\prime}})=\min\limits_{\sigma}\max\limits_{i}(dist(t_{i},t_{\sigma(i)}^{\prime})),$ where $\sigma$ runs over all permutations of $\{1,2,\cdots,m\}$.
And we define $$SPV(\phi)=\sup\limits_{y,y'\in X}dist(Sp\phi_{y},Sp\phi_{y'}).$$

\noindent\textbf{1.23.}~~Let $A$ be of form $PM_{n}(C(X))P$ or $M_{l}(I_{k})$, $B=M_{l_{1}}(I_{k_{1}})$, and $\phi: A\longrightarrow B$ be a unital homomorphism. Then
$SPV(\phi)$ is defined to be $SPV(\phi')$, where  $\phi'=\imath\circ\phi :~A\xrightarrow {\phi}M_{l_{1}}(I_{k_{1}})\xrightarrow {\imath} M_{l_{1}k_{1}}(C[0,1]).$

\noindent\textbf{1.24.}~~Let $A$ and $B$ be either of form $PM_{n}(C(X))P$ (with $X$ path connected) or of form $M_{l}(I_{k})$. Let $\phi: A\longrightarrow B$ be a unital
homomorphism, we say that $\phi$ has property $sdp(\eta,\delta)$ ({\bf spectral distribution property with respect to $\eta$ and  $\delta$}) if for any $\eta$-ball $B_{\eta}(x)=\{x^{\prime}\in X~|~~dist(x^{\prime},x)<\eta)\} \subset  X(=Sp(A))$ and any point $y\in Sp(B)$, $$\sharp(Sp\phi_{y}\cap B_{\eta}(x))\geq\delta\sharp Sp\phi_{y},$$ counting multiplicity. If $\phi$ is not unital, we say that $\phi$ has $sdp(\eta,\delta)$ if the corresponding
unital homomorphism $\phi: A\longrightarrow \phi(\textbf{1}_{A})B\phi(\textbf{1}_{A})$ has property $sdp(\eta,\delta)$.

\noindent\textbf{1.25.}~~Let $A=M_{n}(C(X))$ or $M_{l}(I_{k})$ and let $F\subset A$ be a finite subset. Let $\varepsilon>0$ and $\eta>0$. We say that $F$ is {\bf $(\varepsilon,\eta)$
uniformly continous} if for any pair of points $x_{1},x_{2}\in X$ (or $x_{1},x_{2}\in [0,1]$ for $A=M_{l}(I_{k}))$, $dist(x_{1},x_{2})<\eta$, implies
$\|f(x_{1})-f(x_{2})\|<\varepsilon, \forall f\in F$.

For a finite set  $F\subset M_{n}(C(X))$ (or $F\subset M_{l}(I_{k})$, respectively), each $f\in F$ corresponds to a matrix $(f_{ij})_{n\times n}$ (or $(f_{ij})_{lk\times lk}$
if $F\subset M_{l}(I_{k})\subset M_{lk}(C[0,1])$, respectively). We define $\widetilde{F}\subset C(X)$ (or $\widetilde{F}\subset C[0,1]$, if $F\subset M_{l}(I_{k})$, respectively) to be the collection of
all $f_{ij}$ for all $f=(f_{ij})\in F$. Then,  that  $\widetilde{F}$ is $(\frac{\varepsilon}{n},\eta)$ uniformly continuous (or $(\frac{\varepsilon}{lk},\eta)$ uniformly continuous, respectively) implies that $F$ is $(\varepsilon,\eta)$ uniformly continuous.

For each $F\subset PM_{n}(C(X))P$, $F$ is called $(\varepsilon,\eta)$ uniformly continuous if it is $(\varepsilon,\eta)$ uniformly continuous regarded
as a subset $F\subset M_{n}(C(X))$.

\noindent\textbf{Lemma 1.26.}~~ Let $F\subset A(\in\mathcal{HD})$ be any finite set and $\varepsilon>0$. There is an $\eta>0$ such that if
$\phi: A\longrightarrow B(\in\mathcal{HD})$ satisfies $SPV(\phi)<\eta$, then $\phi(F)$ is weakly approximately constant to within $\varepsilon$---that
is,  $\omega(\phi(F))<\varepsilon$.
\begin{proof}
For $\varepsilon>0$, by uniform continuity of $F$, there is an $\eta>0$ such that if $x_{1},x_{2}\in Sp(A)$ (but $x_1$, or $ x_2\not={\underline{0}}$, or $ {\underline{1}} $) satisfy $dist(x_{1},x_{2})<\eta$,
%(and for the case of $A=M_{l}(I_{k})$, we need to assume $x_{1},x_{2}\notin \{0,1\}$),
then there is a unitary $u$
such that
 $$\|\pi_{x_{1}}(f)-u\pi_{x_{2}}(f)u^{*}\|<\varepsilon~~~\forall f\in F,$$
  where $\pi_{x_{1}},\pi_{x_{2}}$ are corresponding to irreducible
representations. For the case of $A=M_{l}(I_{k})$, we will allow  $x_1$ or $x_2$ to be $0$ or $1$, for which the representation $\pi_{x_1}$ or $\pi_{x_2}$ is no longer irreducible.  It is routine to verify that the lemma holds for such $\eta$.\\
\end{proof}
\noindent\textbf{1.27.}~~Let $X$ be a finite $CW$ complex. Set $F^{k}X=Hom(C(X),M_{k}(\mathbb{C}))_{1}$. The space $F^{k}X$ is compact and metrizable. We can endow the space $F^{k}X$
with a fixed metric $d$ as below. Choose a finite set $\{f_{i}\}^{n}_{i=1}\subset C(X)$ which generates $C(X)$ as a $C^{*}$-algebra. For any $\phi,\psi\in F^{k}X$ which, by definition, are unital homomorphisms from $C(X)$ to $M_{k}(\mathbb{C})$, define
$d(\phi,\psi)=\sum\limits_{i=1}\limits^{n}\|\phi(f_{i})-\psi(f_{i})\|.$

\noindent\textbf{1.28.}~~Set $P^{n}X=\underbrace{X\times X\times\cdots\times X}\limits_{n}/\thicksim$, where the equivalence relation $\thicksim$ is defined by $$(x_{1},x_{2},\cdots,x_{n})\thicksim(x^{\prime}_{1},x^{\prime}_{2},\cdots,x^{\prime}_{n})$$ if there is a permutation $\sigma$ of $\{1,2,\cdots,n\}$ such that
$x_{i}=x^{\prime}_{\sigma(i)}$ for each $1\leq i\leq n$. A metric $d$ on $X$ can be extended to a metric on $P^{n}X$ by $$d([x_{1},x_{2},\cdots,x_{n}],[x^{\prime}_{1},x^{\prime}_{2},\cdots,x^{\prime}_{n}])=\min\limits_{\sigma}\max\limits_{1\leq i\leq n}d(x_{i},x^{\prime}_{\sigma(i)}),$$
where $\sigma$ is taken from the set of all permutations, and $[x_{1},x_{2},\cdots,x_{n}]$ denote the equivalence class of $(x_{1},x_{2},\cdots,x_{k})$ in $P^{k}X$.

\noindent\textbf{1.29.}~~Let $X$ be a metric space with metric $d$. Two $k$-tuple of (possible repeating) points $\{x_{1},x_{2},\cdots,x_{n}\}\subset X$ and
$\{x'_{1},x'_{2},\cdots,x'_{n}\}\subset X$ are said to {\bf be paired within} $\eta$ if there is a permutation $\sigma$ such that
$$d(x_{i},x^{\prime}_{\sigma(i)})<\eta,~~~ i=1,2,\cdots,k. $$ This is equivalent to the following statement. If one regards $[x_{1},x_{2},\cdots,x_{n}]$ and
$[x^{\prime}_{1},x^{\prime}_{2},\cdots,x^{\prime}_{n}]$ as points in $P^{n}X$, then $d([x_{1},x_{2},\cdots,x_{n}],[x^{\prime}_{1},x^{\prime}_{2},\cdots,x^{\prime}_{n}])<\eta.$

\noindent\textbf{1.30}~~For $X=[0,1]$, let $P^{(n,k)}X$, where $n,k\in\mathbb{Z_{+}}\backslash \{0\}$, denote the set of $\frac{n}{k}$ elements from $X$, in which only 0 or 1 may appear
fractional times. That is,  each element in $X$ is of the form $$\{0^{\thicksim\frac{n_{0}}{k}},t_{1},t_{2},\cdots,t_{m},1^{\thicksim\frac{n_{1}}{k}}\}\;\;\;\;\;\;\;(*)$$
with $0<t_{1}\leq t_{2}\leq\cdots\leq t_{m}<1$ and $\frac{n_{0}}{k}+m+\frac{n_{1}}{k}=\frac{n}{k}.$\\
An element in $P^{(n,k)}X$ can always be written as $$\{0^{\thicksim\frac{k_{0}}{k}},t_{1},t_{2},\cdots,t_{i},1^{\thicksim\frac{k_{1}}{k}}\}, \;\;\;\;\;\;\;(**)$$
where $0\leq k_{0}<k$, $0\leq k_{1}<k$, $0\leq t_{1}\leq t_{2}\leq\cdots\leq t_{i}\leq1$ and $\frac{k_{0}}{k}+i+\frac{k_{1}}{k}=\frac{n}{k}.$ (Here $t_{i}$
could be 0 or 1.) In the above two representations $(*)$ and $(**)$, we know $k_{0}\equiv n_{0}(mod\;k)$, $k_{1}\equiv n_{1}(mod\;k)$. Let $$y=[0^{\thicksim\frac{k_{0}}{k}},t_{1},t_{2},\cdots,t_{i},1^{\thicksim\frac{k_{1}}{k}}]\in P^{(n,k)}X~~~\mbox{ and }~~~y^{\prime}=[0^{\thicksim\frac{k^{\prime}_{0}}{k}},t^{\prime}_{1},t^{\prime}_{2},\cdots,t^{\prime}_{i},1^{\thicksim\frac{k^{\prime}_{1}}{k}}]\in P^{(n,k)}X,$$ with $k_{0},k_{1},k^{\prime}_{0},k^{\prime}_{1}\in\{0,1,\cdots,k-1\}.$
We define $dist(y,y')$ as the following: if $k_{0}\neq k^{\prime}_{0}$ or $k_{1}\neq k^{\prime}_{1}$, then $dist(y,y^{\prime})=1$; if $k_{0}=k^{\prime}_{0}$ and $k_{1}= k^{\prime}_{1}$ (consequently $i=i^{\prime}$), then $$dist(y,y^{\prime})=\max\limits_{1\leq j\leq i}|t_{j}-t^{\prime}_{j}|,$$ as we order the $\{t_{j}\}$ and $\{t^{\prime}_{j}\}$ as
$t_{1}\leq t_{2}\leq\cdots\leq t_{i}$ and $t^{\prime}_{1}\leq t^{\prime}_{2}\leq\cdots\leq t^{\prime}_{i}$, respectively.

Note that $P^{(n,1)}X=P^{n}X$ with the same metric. Let $\phi,\psi: I_{k}\longrightarrow M_{n}(\mathbb{C})$ be two unital homomorphisms. Then $Sp\phi$ and
$Sp\psi$ define two elements in $P^{(n,k)}[0,1]$. We say that $Sp\phi$ and $Sp\psi$ can be paired within $\eta$, if $dist(Sp\phi,Sp\psi)<\eta$.

Note that
if $dist(Sp\phi,Sp\psi)<1$, then $KK(\phi)=KK(\psi)$.

\noindent\textbf{1.31.}~~Let $A=PM_{k}(C(X))P$ or $M_{l}(I_{k})$ and $X_{1}\subset Sp(A)$ be a closed subset---that is,  $X_{1}$ is a closed subset of $X$ or of $[0,1]$.We define $A|_{X_{1}}$ to be the quotient algebra $A/I$, where $I=\{f\in A, f|_{X_{1}}=0\}.$
Evidently $Sp(A|_{X_{1}})=X_{1}$. If $B=QM_{k}(C(Y))Q$ or $M_{l_{1}}(I_{k_{1}})$, $\phi: A\longrightarrow B$ is a homomorphism, and $Y_{1}\subset Sp(B)(=Y\;or\;[0,1])$ is a closed subset, then we use $\phi|_{Y_{1}}$ to denote the composition $\phi|_{Y_{1}}: A\xrightarrow{\phi} B\rightarrow B|_{Y_{1}}.$\\
If $Sp(\phi|_{Y_{1}})\subset X_{1}\cup X_{2}\cup\cdots\cup X_{k}$, where $ X_{1},X_{2},\cdots,X_{k}$ are mutually disjoint closed subsets of $X$, then the
homomorphism $\phi|_{Y_{1}}$ factors as
$A\longrightarrow A|_{X_{1}\cup X_{2}\cup\cdots\cup X_{n}}=\bigoplus\limits^{n}\limits_{i=1}A|_{X_{i}}\longrightarrow B|_{Y_{1}}.$\\
We will use $\phi|^{X_{i}}_{Y_{1}}$ to denote the part of $\phi|_{Y_{1}}$ corresponding to the map $A|_{X_{i}}\longrightarrow B|_{Y_{1}}.$ Hence
$\phi|_{Y_{1}}=\bigoplus\limits_{i}\phi|^{X_{i}}_{Y_{1}}.$

\noindent\textbf{1.32.}~~In 1.31, if $Y_{1}=\{y\}\subset Y$, then we also denote $\phi|_{\{y\}}$ by $\phi|_{y}$ (or some time $\phi_y$). That is,  we do not differentiate the point $y$ and set $\{y\}$ of
single point $y$. If $\phi: A\longrightarrow M_{l}(I_{k})$, then for any $t\in[0,1]$, $\phi|_{t}$ is the homomorphism from $A$ to $M_{lk}(\mathbb{C})$.
But $\phi|_{0}$ and $\phi|_{1}$ are homomorphisms to the subalgebra $M_{l}(\mathbb{C})\otimes \textbf{1}_{k}$. We will use $\phi|_{\underline{0}}$ and $\phi|_{\underline{1}}$ to denote the corresponding homomorphisms to $M_{l}(\mathbb{C})$ (rather than $M_{lk}(\mathbb{C})$). That is, if $\pi_{\underline{0}}$ and $\pi_{\underline{1}}$ are the irreducible representations corresponding to the spectra 0 and 1, then
$$\phi|_{\underline{0}}=\pi_{\underline{0}}\circ\phi: A\longrightarrow M_{l}(\mathbb{C}),~~~\mbox{and} ~~~\phi|_{\underline{1}}=\pi_{\underline{1}}\circ\phi: A\longrightarrow M_{l}(\mathbb{C}).$$
We still reserve the notation $\pi_{0}$ and $\pi_{1}$, as the evaluation $\pi_{0}(f)=f(0)$ and $\pi_{1}(f)=f(1)$. Note that as notation in 1.15, $\pi_{\underline{0}}(f)=\underline{f}(0)=f(\underline{0})$ and $\pi_{\underline{1}}(f)=\underline{f}(1)=f(\underline{1})$.

From now on, we will not use the notation $\underline{f}$ any more---since $\underline{f}(0)$ and $\underline{f}(1)$ are already written as $f(\underline{0})$ and $f(\underline{1})$; and for $t\in (0,1)$, $\underline{f}(t)$ is written as $f(t)$.

\noindent\textbf{1.33.}~~Two inductive limit systems ($A_{n},\phi_{nm}$) and ($B_{n},\psi_{nm}$) are said to be shape equivalent if there are subsequences $l_{1}<l_{2}<l_{3}<\cdots$ and $k_{1}<k_{2}<k_{3}<\cdots$ and homomorphisms $$\alpha_{n}: A_{l_{n}}\longrightarrow B_{k_{n}}~~~\mbox{and}~~~\beta_{n}:B_{k_{n}}\longrightarrow A_{l_{n+1}}$$ such that $$\beta_{n}\circ\alpha_{n}\thicksim_{h}\phi_{l_{n},l_{n+1}}, ~~~\mbox{ and }~~~\alpha_{n+1}\circ\beta_{n}\thicksim_{h} \psi_{k_{n},k_{n+1}}.$$

\noindent\textbf{2. The invariant}

The invariant used in  the classification in this paper,  called $Inv(A)$ is introduced in [GJL], we will give brief summary here. For detailed discussion of  this invariant, we refer the readers to [GJL].

\noindent\textbf{2.1.}~~Let A be a $C^{*}$-algebra. $K_{0}(A)^{+}\subset K_{0}(A)$ is defined to be the semigroup of $K_{0}(A)$ generated
by $[p]\in K_{0}(A)$, where $p\in M_{\infty}(A)$ are projections. For all $C^{*}$-algebras considered in this paper, for
example, $A\in \mathcal{HD}$, or $A$ is $A\mathcal{HD}$ algebra, or $A=B\otimes C(T_{\uppercase\expandafter{\romannumeral2},k}\times S^{1})$, where B is $\mathcal{HD}$ or $A\mathcal{HD}$
algebra, we will always have
$$K_{0}(A)^{+}\bigcap (-K_{0}(A)^{+})=\{0\}~~~\mbox{and}~~~K_{0}(A)^{+}-K_{0}(A)^{+}=K_{0}(A).$$
Therefore $(K_{0}(A),K_{0}(A)^{+})$ is an ordered group. Define $\Sigma A \subset K_{0}(A)^{+}$ to be
$$\Sigma A = \{[p]\in K_{0}(A)^{+}, p~ is ~a ~projection ~in ~A\}.$$
Then $(K_{0}(A),K_{0}(A)^{+},\Sigma A)$ is a scaled ordered group.

\noindent\textbf{2.2.}~~Let $\underline K (A)=K_{*}(A)\bigoplus \bigoplus_{k=2}^{+\infty}K_{*}(A,\mathbb{Z}/ k\mathbb{Z})$ be as in [DG].
Let $\wedge$ be the Bockstein operation on $\underline K (A)$(see 4.1 of [DG]). It is well known that
$K_{*}(A,Z\oplus \mathbb{Z}/ k\mathbb{Z})=K_{0}(A\otimes C(W_{k}\times S^{1})),$
where $W_{k}=T_{II,k}$. \\As in [DG], let
$K_{*}(A,Z\oplus \mathbb{Z}/ k\mathbb{Z})^{+}=K_{0}(A\otimes C(W_{k}\times S^{1})^{+})$
and let $\underline K (A)^{+}$ be the semigroup generated by $\{K_{*}(A,\mathbb{Z}\oplus\mathbb{Z}/ k\mathbb{Z})^{+},k=2,3,\cdots\}$.

\noindent\textbf{2.3.}~~Let $Hom_{\wedge} (\underline K (A),\underline K (B))$ be the set of homomorphisms between $\underline K (A)$ and $\underline K (B)$
compatible with Bockstein operation $\wedge$. There is   a surjective map (see [DG])
$$\Gamma: KK(A,B)\rightarrow Hom_{\wedge}(\underline K (A),\underline K (B)).$$
Following  [R] and [DL], we denote $KL(A,B)\triangleq KK(A,B)/ker\;\Gamma.$ The triple
 $(\underline{K}(A),\underline{K}(A)^{+},\Sigma A)$ is part of our invariant.
For two $C^{*}$-algebras $A$ and $B$, by a``homomorphism"
$$\alpha:~(\underline{K}(A),\underline{K}(A)^{+},\Sigma A) \to (\underline{K}(B),\underline{K}(B)^{+},\Sigma B),$$
 we mean a system of maps:
$$\alpha^{i}_{k}: ~K_{i}(A,\mathbb{Z}/k\mathbb{Z})\longrightarrow K_{i}(B,\mathbb{Z}/k\mathbb{Z}),~~i=0,1,~~k=0,2,3,\cdots$$ which are compatible with Bechstein
operations and $\alpha=\oplus_{k,i}\alpha^{i}_{k}$ satisfies $\alpha(\underline{K}(A)^{+})\subset\underline{K}(B)^{+}$. And finally,  $\alpha^{0}_{0}(\Sigma A)\subset\Sigma B$.

\noindent\textbf{2.4.}~~ As in 2.8 of [GJL], for a unital $C^{*}$-algebra $A$, let $TA$ denote the space of tracial states of $A$, $AffTA$ denote the Banach space of  all affine continuous maps from $TA$ to $\mathbb{R}$.
Let $(AffTA)_{+}$  be the subset of $AffTA$ consisting of all nonnegative affine functions.
An element
$\textbf{1}\in AffTA$, defined by $\textbf{1}(\tau)=1$ for all $\tau\in TA$,  is called the unit of $AffTA$. ($AffTA$, $(AffTA)_{+}$, \textbf{1}) forms a scaled ordered real Banach space. Any unital homomorphism $\phi: A\longrightarrow B$ induces a continuous affine map $T\phi: TB\longrightarrow TA$, which, in turn,  induces a unital  positive linear map $AffT\phi: AffTA\longrightarrow AffTB.$

If $\phi: A\longrightarrow B$ is not unital, we still use $AffT\phi$ to denote the unital  positive linear map $$AffT\phi: AffTA\longrightarrow AffT(\phi(\textbf{1}_A)B\phi(\textbf{1}_A))$$
by regarding $\phi$ as the unital homomorphism from $A$ to $\phi(\textbf{1}_A)B\phi(\textbf{1}_A)$. The homomorphism $\phi$ also induce a positive linear map, denoted by $\phi_T$ to avoid the confusion, from $AffTA$ to $AffTB$. But this map  will not preserve the
unit $\textbf{1}$. It has the property that $\phi_T(\textbf{1}_{AffTA})\leq\textbf{1}_{AffTB}$. In particular, we will often use the notation $\phi_T$ for the following situation: If $p_1<p_2$ are two projections in $A$, and $\phi=\imath:~p_1Ap_1 \longrightarrow p_2Ap_2 $ is the inclusion, then $\imath_T$ will denote the (not necessarily unital) map from $AffT(p_1Ap_1)$ to $AffT(p_2Ap_2)$ induced by $\imath$. (See 2.8 of [GJL] for details.)

For a unital $C^*$ algebra $A$, define $\rho_A: K_0(A) \to AffT(A)$ by $\rho_A([p])(\tau)=\sum_{i=1}^n\tau(p_{ii})$ for any element  $[p]\in K_0(A)^+\subset  K_0(A)$ represented by the projection $p=(p_{ij})\in M_n(A)$. Let us use $\widetilde{\rho K}_{0}(A)$ to denote the closure of the real vector
space spanned by $\overline{\rho K_{0}(A)}$. That is,  $$\widetilde{\rho K}_{0}(A)=\overline{\{\Sigma\lambda_{i}\phi_{i}~|~~\lambda_{i}\in\mathbb{R},\phi_{i}\in\rho_{0} K_{0}(A)}\}.$$
Follow 2.26 of [GJL], let $\widetilde{d'}$ denote the
quotient metric on $AffTA/\widetilde{\rho K}_{0}(A)$ of $AffTA$, that is,
$$\widetilde{d'}(f,g)=inf\{\|f-g-h\|~|~~h\in\widetilde{\rho K_{0}(A)}\}~~~~~\forall f,g\in AffTA/\widetilde{\rho K}_{0}(A).$$ Define $\widetilde{d}_{A}$ by
$$\widetilde{d}_{A}(f,g)=
 \left\{ \begin{lgathered}
  2,\;\;\; \;\;\;\;\;\;\;\;\;\;\;\;\;\;\;\;\;\;if\;\widetilde{d'}(f,g)\geq\frac{1}{2} \\
  |e^{2\pi i\widetilde{d'}(f,g)}-1|,\;\;\;if\;\widetilde{d'}(f,g)<\frac{1}{2}~~~.
   \end{lgathered} \right.$$

\noindent\textbf{2.5.}~~Recall that the commutator subgroup of a group $G$ is the subgroup generated by all elements of the form $aba^{-1}b^{-1}$, where $a,b\in G$. Let $A$ be a unital $C^{*}$-algebra. Let $U(A)$ denote the group of unitaries of $A$ and,  $U_{0}(A)$, the connected component of $\textbf{1}_{A}$ in $U(A)$. Let $DU(A)$ and $DU_{0}(A)$ denote the commutator subgroups of $U(A)$ and $U_{0}(A)$, respectively. By Proposition 2.23 of [GJL], for any  $A\in\mathcal{HD}$ or $A\in A\mathcal{HD}$, we have $\overline{DU_{0}(A)}=\overline{DU(A)}$.

\noindent\textbf{2.6.} Let us  introduce the extended commutator group $DU^{+}(A)$, which is generated by $DU(A)\subset U(A)$ and the set
$\{e^{2\pi itp}=e^{2\pi it}p+(\textbf{1}-p)\in U(A)~|~~t\in\mathbb{R},p\in A ~\mbox{ is a projection}\}.$
Let $\widetilde{DU}(A)$ denote the closure of $DU^{+}(A)$. That is,  $\widetilde{DU}(A)=\overline{DU^{+}(A)}$. Define
$$\widetilde{SU}(A)\triangleq\{\overline{x\in U(A)~|~~x^{n}\in\widetilde{DU}(A)~~~~~\mbox{for some} ~n\in\mathbb{Z}_{+}\backslash\{0\}}\}.$$

For any $u,v\in U(A)$ $$\overline{D}_{A}([u],[v])=inf\{\|uv^{*}-c\|~|~~c\in\widetilde{DU}(A)\},$$
where $[u], [v]$ are regarded as elements in $ U(A)/\widetilde{DU}(A)$; and define  $$\widetilde{D}_{A}([u],[v])=inf\{\|uv^{*}-c\|~|~~c\in\widetilde{SU}(A)\},$$
where $[u], [v]$ are regarded as elements in $ U(A)/\widetilde{SU}(A)$.

Let $U_{tor}(A)$ denote the set of unitaries $u\in A$ such that $[u]\in tor\;K_{1}(A)$. Then for $A\in\mathcal{HD}$ or $A\in A\mathcal{HD}$, we have
$$\widetilde{SU}(A)\subset U_{tor}(A), ~~\widetilde{DU}(A)=U_{0}(A)\cap\widetilde{SU}(A)~~\mbox{ and}~~U_{tor}(A)=U_{0}(A)\cdot\widetilde{SU(A)}.$$ Evidently, we have $U_{0}(A)/\widetilde{DU}(A)\cong U_{tor}(A)/\widetilde{SU}(A).$
%The metric $\overline{D_{A}}$ on $U(A)/\widetilde{DU(A)}$ induces a metric $\widetilde{D}_{A}$ on $U(A)/\widetilde{SU(A)}$.
The identification $U_{0}(A)/\widetilde{DU(A)}$ with $U_{tor}(A)/\widetilde{SU(A)}$ is an isometry with respect to $\overline{D}_A$ and $\widetilde{D}_{A}$. (See 2.26, 2.28 and 2.30 of [GJL].)\\

The following lemma is Lemma 2.31 in [GJL]

\noindent\textbf{Lemma 2.7.}~~ Let $A\in\mathcal{HD}$ or $A\in A\mathcal{HD}$.\\
(1)~~ There is a split exact sequence $$0\rightarrow AffTA/\widetilde{\rho K_{0}(A)}\xrightarrow{\widetilde{\lambda}_{A}}U(A)/\widetilde{SU(A)}\xrightarrow{\pi_{A}} K_{1}(A)/tor\;K_{1}(A)\rightarrow0.$$
(2)~~ $\widetilde{\lambda}_{A}$ is an isometry with respect to the metrics $\widetilde{d}_{A}$ and $\widetilde{D}_{A}$.

\noindent\textbf{2.8.} (see 2.32 of {GJL])~~ For each pair of projections $p_{1},p_{2}\in A$ with $p_{1}=up_{2}u^{*},$ $$U(p_{1}Ap_{1})/\widetilde{SU(p_{1}Ap_{1})}\cong U(p_{2}Ap_{2})/ \widetilde{SU(p_{2}Ap_{2})}.$$ Also,  since in any unital $C^{*}$-algebra $A$ and unitaries $u,v\in U(A)$, $v$ and $uvu^{*}$ represent
 a same element in $U(A)/\widetilde{SU(A)}$, and the above identification does not depend on the choice
of $u$ to implement $p_{1}=up_{2}u^{*}$. That is for any $[p]\in\Sigma A$, the group $U(pAp)/\widetilde{SU(pAp)}$ is well defined, which does not depend on choice of $p\in[p]$. We will include this group (with metric) as part of our invariant. If $[p]\leq[q]$, then we can choose $p,q$ such that $p\leq q$. In this case,  there is a natural inclusion map $\imath: pAp\longrightarrow qAq$  which induces $$\imath_{*}: U(pAp)/\widetilde{SU(pAp)}\longrightarrow U(qAq)/\widetilde{SU(qAq)},$$ where $\imath_{*}$ is defined by
$$\imath_{*}(u)=u\oplus(q-p)\in U(qAq),~~~~ \forall u\in U(pAp).$$
~~~A unital homomorphism $\phi: A\to B$ induces a contractive group homomorphism
$$\phi^{\natural}: U(A)/\widetilde{SU(A)} \longrightarrow U(B)/\widetilde{SU(B)}.$$
If $\phi$ is not  unital, then the map $\phi^{\natural}: U(A)/\widetilde{SU(A)} \longrightarrow U(\phi(\textbf{1}_A)B\phi(\textbf{1}_A))/\widetilde{SU}(\phi(\textbf{1}_A)B\phi(\textbf{1}_A))$ is induced by the corresponding unital homomorphism.  In this case,  $\phi$ also induces the map $\imath_*\circ \phi^{\natural}: U(A)/\widetilde{SU(A)} \longrightarrow U(B)/\widetilde{SU(B)}$, which is denoted by $\phi_*$ to avoid confusion.

Since $U(A)/\widetilde{SU(A)} $ is an Abelian group, we will call the unit $[{\bf 1}]\in U(A)/\widetilde{SU(A)} $  the zero element. If $\phi: A \to B$ satisfies $\phi(U(A))\subset \widetilde{SU}(\phi(\textbf{1}_A)B\phi(\textbf{1}_A))$, then $\phi^{\natural}=0$. In particular, if the image of $\phi$ is of finite dimensional, then $\phi^{\natural}=0$.

\noindent\textbf{Definition 2.9.}~~ (see 2.33 of [GJL]) In this paper, we will denote $$(\underline{K}(A),\underline{K}(A)^{+},\Sigma A,\{AffT(pAp)\}_{[p]\in\Sigma A},\{U(pAp)/\widetilde{SU(pAp)}\}_{[p]\in\Sigma A})$$ by $Inv(A)$. By a map from $Inv(A)$ to $Inv(B)$, we mean
$$\alpha: (\underline{K}(A),\underline{K}(A)^{+},\Sigma A)\longrightarrow(\underline{K}(B),\underline{K}(B)^{+},\Sigma B)$$ as in 2.3, and for each pair $([p], [\overline{p}]) \in\Sigma A\times \Sigma B$ with $\alpha([p])=[\overline{p}]$, there are an associate unital positive linear map
$$\xi^{p,\overline{p}}: AffT(pAp)\longrightarrow AffT(\overline{p}B\overline{p})$$ and an associate contractive group homomorphism
$$\chi^{p,\overline{p}}: U(pAp)/\widetilde{SU(pAp)}\longrightarrow U(\overline{p}B\overline{p})/\widetilde{SU(\overline{p}B\overline{p})}$$
satisfy the following compatibility conditions.\\
(a)~~ If $p<q$, then the diagrams
$$\CD
  AffT(pAp) @>\xi^{p,\overline{p}}>> AffT(\overline{p}B\overline{p}) \\
  @V\imath_T  VV @V \imath_T VV \\
  AffT(qAq) @>\xi^{q,\overline{q}}>> AffT(\overline{q}B\overline{q})
\endCD\;\;\;\;\;\;\;\;\;\;\;\;\;\;\;\;\;\;\;\;\;\;\;\;\;~~~~~~~~~(\uppercase\expandafter{\romannumeral1})$$ and
$$\CD
  U(pAp)/\widetilde{SU(pAp)} @>\chi^{p,\overline{p}}>> U(\overline{p}B\overline{p})/\widetilde{SU(\overline{p}B\overline{p})} \\
  @V\imath_*  VV @V \imath_* VV   \\
  U(qAq)/\widetilde{SU(qAq)} @>\chi^{q,\overline{q}}>> U(\overline{q}B\overline{q})/\widetilde{SU(\overline{q}B\overline{q})}
\endCD\;\;\;\;\;\;\;~~~~~~~~~~~(\uppercase\expandafter{\romannumeral2})$$
commutes, where the vertical maps are induced by inclusions.\\
(b)~~ The following diagram commutes
$$\CD
  K_{0}(pAp) @>\rho>> AffT(pAp) \\
  @V \alpha VV @V \xi^{p,\overline{p}} VV  \\
  K_{0}(\overline{p}B\overline{p}) @>\rho>> AffT(\overline{p}B\overline{p})
\endCD\;\;\;\;\;\;\;\;\;\;\;\;\;\;\;\;\;\;\;\;\;\;\;\;\;\;\;\;\;\;\;\;\;\;~~~~~~~(\uppercase\expandafter{\romannumeral3})$$
and therefore $\xi^{p,\overline{p}}$ induces a map (still denoted by $\xi^{p,\overline{p}}$):
$$\xi^{p,\overline{p}}: AffT(pAp)/\widetilde{\rho K_{0}(pAp)}\longrightarrow AffT(\overline{p}B\overline{p})/\widetilde{\rho K_{0}(\overline{p}B\overline{p})}.$$
(The commutativity  of $(\uppercase\expandafter{\romannumeral3})$ follows from the commutativity of $(\uppercase\expandafter{\romannumeral1})$, by 1.20 of [Ji-Jiang]. So this is not an extra requirement.)\\
(c)~~ The following diagrams
$$\CD
  AffT(pAp)/\widetilde{\rho K_{0}(pAp)} @>>> U(pAp)/\widetilde{SU(pAp)} \\
  @V \xi^{p,\overline{p}} VV @V \chi^{p,\overline{p}} VV  \\
  AffT(\overline{p}B\overline{p})/\widetilde{\rho K_{0}(\overline{p}B\overline{p})} @>>> U(\overline{p}B\overline{p})/\widetilde{SU(\overline{p}B\overline{p})}
\endCD\;\;\;\;\;~~~~~~~~~~~(\uppercase\expandafter{\romannumeral4})$$ and
$$\CD
  U(pAp)/\widetilde{SU(pAp)} @>>> K_{1}(pAp)/tor\;K_{1}(pAp) \\
  @V\chi^{p,\overline{p}} VV @V \alpha_{1} VV  \\
  U(\overline{p}B\overline{p})/\widetilde{SU(\overline{p}B\overline{p})} @>>> K_{1}(\overline{p}B\overline{p})/tor\;K_{1}(\overline{p}B\overline{p})
\endCD\;\;\;\;\;\;\;\;\;~~~~~~~~~~~(\uppercase\expandafter{\romannumeral5})$$
commute, where $\alpha_{1}$ is induced by $\alpha$.

We will denote the map from $Inv(A)$ to $Inv(B)$ by
$$(\alpha,\xi,\chi): (\underline{K}(A),\{AffT(pAp)\}_{[p]\in\Sigma A},\{U(pAp)/\widetilde{SU(pAp)}\}_{[p]\in\Sigma A})\longrightarrow~~~~~~~~~~~~~~~~~~~~$$
$$~~~~~~~~~~~~~(\underline{K}(B),\{AffT(\overline{p}B\overline{p})\}_{[\overline{p}]\in\Sigma B},\{U(\overline{p}B\overline{p})/\widetilde{SU(\overline{p}B\overline{p})}\}_{[\overline{p}]\in\Sigma B}).$$

The  part
$(\underline{K}(A),\underline{K}(A)^{+},\Sigma A,\{AffT(pAp)\}_{[p]\in\Sigma A},)$ of $Inv(A)$ is introduced in [Jiang1] and we will denote it by  $Inv^0(A)$. By a map from $Inv^0(A)$ to $Inv^0(B)$, we mean a map $\alpha: (\underline{K}(A),\underline{K}(A)^{+},\Sigma A)\longrightarrow(\underline{K}(B),\underline{K}(B)^{+},\Sigma B)$ as in the above, and for each pair $([p], [\overline{p}]) \in\Sigma A\times \Sigma B$ with $\alpha([p])=[\overline{p}]$, there is an associate unital positive linear map
$\xi^{p,\overline{p}}: AffT(pAp)\longrightarrow AffT(\overline{p}B\overline{p})$ as the above with the diagram $(I)$ commutes.

In the above definition of $Inv(A)$ one can replace  $U(pAp)/\widetilde{SU}(pAp)$ by  $U(pAp)/\overline{DU(pAp)}$ as the part of the invariant. This is called $Inv'(A)$ in 2.37 of [GJL]. We choose the formulation of $Inv(A)$, since it is much more convenient for the proof of the main theorem here and it is formally a weaker requirement than the one to require the isomorphism between $Inv'(A)$ and $Inv'(B)$, and the isomorphism theorem is formally stronger.

For the convenient of reader, we will quote the following results and notations from [GJL] (see 2.39--2.41 of [GJL]).

\noindent\textbf{2.10.}~~In general, for $A=\oplus A^{i}$, $\widetilde{SU}(A)=\oplus_{i}\widetilde{SU}(A^{i}).$ For $A=PM_{l}(C(X))P\in\mathcal{HD}$, $\widetilde{SU}(A)=\widetilde{DU}(A).$ For $A=M_{l}(I_{k})$, $\widetilde{SU}(A)=\widetilde{DU}(A)\oplus K_{1}(A).$
For both cases, $U(A)/\widetilde{SU}(A)$ can be identified with \\$C_{1}(X,S^{1}):=C(X,S^{1})/\{constant\; functions\}$ or with
$C_{1}([0,1],S^{1})=C([0,1],S^{1})/\{constant\; functions\},$ for $A=M_l(I_{k}).$

Furthermore, $C_{1}(X,S^{1})$ can be identified as the set of continuous functions from $X$ to
$S^{1}$ such that $f(x_{0})=1$ for certain fixed base point $x_{0}\in X$. For $X=[0,1]$, we choose 0 to be the base point. For $X=S^{1}$, we choose $1\in S^{1}$ to be the base point.

\noindent\textbf{2.11.}~~Let $A=\oplus^{n}_{i=1}A^{i}\in\mathcal{HD}$, $B=\oplus^{m}_{j=1}B^{j}\in\mathcal{HD}$. In this subsection we will discuss some consequences of the compatibility of the maps between $AffT$ spaces. Let $$p=\oplus p^{i}<q=\oplus q^{i}\in A~~~\mbox{ and }~~~\overline{p}=\oplus_{j=1}^{m} \overline{p}^{j}<\overline{q}=\oplus_{j=1}^{m} \overline{q}^{j}\in B$$ be projections satisfying $\alpha([p])=[\overline{p}]$ and $\alpha([q])=[\overline{q}]$. Suppose that two unital positive linear maps $\xi_{1}: AffTpAp\longrightarrow AffT\overline{p}B\overline{p}$ and
$\xi_{2}: AffTqAq\longrightarrow AffT\overline{q}B\overline{q}$ are compatible with $\alpha$ (see diagram $(III)$ in 2.9) and compatible with each
other (see diagram $(I)$ in 2.9). Since the (not necessarily unital) maps $AffTpAp\longrightarrow AffTqAq$ and $AffT\overline{p}B\overline{p}\longrightarrow AffT\overline{q}B\overline{q}$ induced by inclusions, are injective, we know that the map $\xi_{1}$
is completely decided by $\xi_{2}$. Let $$\xi_{2}^{i,j}: AffTq^{i}Aq^{i}\longrightarrow AffT\overline{q^{j}}B^{j}\overline{q^{j}}~~(\mbox{ or}~~ \xi_{1}^{i,j}: AffTp^{i}Ap^{i}\longrightarrow AffT\overline{p^{j}}B^{j}\overline{p^{j}})$$ be the corresponding component of the map $\xi_{2}$ (or $\xi_{1}$). If $p^{i}\neq0$
and $\overline{p}^{j}\neq0$, then $\xi_{1}^{i,j}$ is given by the following formula, for any $f\in AffTp^{i}A^{i}p^{i}=C_{\mathbb{R}}(SpA^{i})(\cong AffTq^{i}Aq^{i})$,
 $$\xi_{1}^{i,j}(f)=\frac{rank\;\overline{q_{j}}}{rank\;\overline{p_{j}}}\cdot\frac{rank\;\alpha^{i,j}(p^{i})}{rank\;\alpha^{i,j}(q^{i})}\cdot\xi_{2}^{i,j}(f).$$
In particular, if $q=\textbf{1}_{A}$ with $\overline{q}=\alpha_{0}[\textbf{1}_{A}]$,
and $\xi_{2}=\xi: AffTA\longrightarrow Aff\alpha_{0}[\textbf{1}_{A}]B\alpha_{0}[\textbf{1}_{A}]$ (note that since $AffTQBQ$ only depends on the unitary equivalence class of $Q$,  it is convenient to denote it as $AffT[Q]B[Q]$), then we will denote $\xi_{1}$ by  $\xi|_{([p],\alpha[p])}$. Even for the general case, we can also write $\xi_{1}=\xi_{2}|_{([p],\alpha[p])}$, when $p<q$ as above.

\noindent\textbf{2.12.}~~As in 2.40, let $A=\oplus_{i=1}^{n}A^{i}$, $B=\oplus_{j=1}^{m}B^{j}$ and $p<q\in A$, $\overline{p}<\overline{q}\in B$, with $\alpha_{0}[p]=[\overline{p}]$ and $\alpha_0[q]=[\overline{q}]$. If $$\gamma_{1}: U(pAp)/\widetilde{SU}(pAp)\longrightarrow U(\overline{p}B\overline{p})/\widetilde{SU}(\overline{p}B\overline{p})$$ is compatible with
$$\gamma_{2}: U(qAq)/\widetilde{SU}(qAq)\longrightarrow U(\overline{q}B\overline{q})/\widetilde{SU}(\overline{q}B\overline{q}),$$ then $\gamma_{1}$ is completely determined by $\gamma_{2}$ (since both maps $$U(pAp)/\widetilde{SU}(pAp)\longrightarrow U(qAq)/\widetilde{SU}(qAq),~~~~ U(\overline{p}B\overline{p})/\widetilde{SU}(\overline{p}B\overline{p})\longrightarrow U(\overline{q}B\overline{q})/\widetilde{SU}(\overline{q}B\overline{q})$$
are injective). Therefore we can denote $\gamma_{1}$ by $\gamma_{2}|_{([p],\alpha[p])}$.

\vspace{0.3in}

\noindent\textbf{\S3. Shape equivalence}

The following dichotomy result is due to Pasnicu.

\noindent\textbf{Proposition 3.1.}~~(see [Pa1, Lemma 2.9]) Let $A\!=\!lim(A_{n}\!=\!\oplus A_{n}^{i},~\phi_{n,m})$ be an $A\mathcal{HD}$ inductive limit. Suppose that $A$ has the ideal property. For any $\delta>0$ and $A_{n}$, there is $m_{0}>n$ such that the following is true:

For any closed subset $X\subset SpA^{i}_{n}$ and any $m>m_{0}$, the homomorphism $\phi_{n,m}^{i,j}$ satisfy the dichotomy condition: $$\mbox{either}~~~~
Sp(\phi_{n,m}^{i,j})_{y}\cap X=\emptyset,~\forall y\in SpA_{m}^{j}~~~~\mbox{ or}~~~~Sp(\phi_{n,m}^{i,j})_{y}\cap B_{\delta}(X)\neq\emptyset,~ \forall y\in SpA_{m}^{j}.$$
\begin{proof}
Lemma 2.9 of [Pa1] says the theorem is true for $AH$ inductive limit with the ideal property. One can verify that  the proof of Lemma 2.9 works for our case. Notice that the test functions used in both Lemma 2.8 and Lemma 2.9 of [Pa1] (the functions $\widetilde{f}$ and $f_{1}$ in the
proof of Lemma 2.8 and $f$ in the proof of Lemma 2.9) are in the center of $A_{n}^{i}=P_{n,i}M_{[n,i]}(C(X_{n,i}))P_{n,i}$. If $A_{n}^{i}=M_{l}(I_{k})$, we can choose the same functions in $M_{lk}(C[0,1])$---it is in $M_{l}(I_{k})$ automatically since the functions are in the center. Then the result follows the same
proof words by words.\\
\end{proof}
The following two results follow the same way as in the proofs of Lemma 2.10 and Theorem 2.6 of [Pa1],  provided that we replace Lemma 2.9 of [Pa1] by the above proposition.

\noindent\textbf{Proposition 3.2.}~~[Pa1, Lemma 2.10] Let $A=lim(A_{n}\!=\!\oplus A_{n}^{i},~\phi_{n,m})$ be an $A\mathcal{HD}$ inductive limit system. Suppose that $A$ has the ideal property. Fix $A_{n}, F=\oplus_i F_{n}^{i}\subset\oplus_i A_{n}^{i}=A_{n}$, and $ \varepsilon>0.$ Let $\eta>0$ such that $F$ is $(\varepsilon,\eta)$ continuous (see 1.25 for the terminology). Let $L>0$. There is $m_{0}>n$ such that each partial map $\phi_{n,m}^{i,j}(m>m_{0})$ satisfies either\\
(a)~~ $\phi_{n,m}^{i,j}$  is $L$-large (see 1.12 for $L$-large), or\\
(b)~~ There exist finitely many points $x_{1},x_{2},\cdots,x_{\bullet}\in Sp(A_{n}^{i})$ and contractible subspaces $X_{1},X_{2},\cdots,X_{\bullet}$ such that
$X_{k}\cap X_{l}=\emptyset$ for $k\neq l$ and $x_{k}\in X_{k}\subset B_{\eta}(x_{k})$, and such that $Sp\phi_{n,m}^{i,j}\subset \cup_{k=1}X_{k}.$\\
Consequently, there is a unital homomorphism $$\psi_{n,m}^{i,j}: A_{n}^{i}\longrightarrow\phi_{n,m}^{i,j}(\textbf{1}_{A_{n}^{i}})A_{m}^{j}\phi_{n,m}^{i,j}(\textbf{1}_{A_{n}^{i}})$$ factoring through a finite dimensional $C^{*}$-algebra as
$$A_{n}^{i}\longrightarrow A_{n}^{i}|_{\{x_{1},x_{2},\cdots,x_{\bullet}\}}\longrightarrow\phi_{n,m}^{i,j}(\textbf{1}_{A_{n}^{i}})A_{m}^{j}\phi_{n,m}^{i,j}(\textbf{1}_{A_{n}^{i}})$$ such that
$\phi_{n,m}^{i,j}$ is homotopic to $\psi_{n,m}^{i,j}$ (within the projection $\phi_{n,m}^{i,j}(\textbf{1}_{A_{n}^{i}})$) and $$
\|\phi_{n,m}^{i,j}(f)-\psi_{n,m}^{i,j}(f)\|<\varepsilon,~~~~ \forall f\in F_{n}^{i}.$$

\noindent\textbf{Proposition 3.3.}~~[Pa1, Lemma 2.6] Let $A\!=\!lim(A_{n}=\oplus A_{n}^{i},~\phi_{n,m})$ be an $A\mathcal{HD}$ inductive limit system.
Suppose that $A$ has the ideal property. Then there exist a sequence $l_{1}<l_{2}<\cdots<l_{n}<\cdots$, and a homomorphism
$\widetilde{\phi}_{n, n+1}: A_{l_{n}}\longrightarrow A_{l_{n+1}}$ such that\\
(a)~~ $\widetilde{\phi}_{n, n+1}$ is homotopic to $\phi_{l_{n},l_{n+1}}$ and\\
(b)~~ $SPV(\widetilde{\phi}_{n, m})\leq 2^{-m}$ for all $n<m$.\\
Consequently,  $\widetilde{A}\!=\!lim(A_{l_{l}},~\widetilde{\phi}_{n, n+1})$ is of real rank zero.

\noindent\textbf{3.4.}~~Let $A\!=\!lim(A_{n}\!=\!\oplus A_{n}^{i},~\phi_{n,m})$ and $B\!=\!lim(B_{n}\!=\!\oplus B_{n}^{i},~\psi_{n,m})$ be two $A\mathcal{HD}$ inductive limit systems whose limit
algebras $A$ and $B$ have the ideal property. If $$\alpha: (\underline{K}(A),\underline{K}(A)^{+},\Sigma A)\longrightarrow(\underline{K}(B),\underline{K}(B)^{+},\Sigma B)$$ is an isomorphism (compatible with Bockstein operations), then by 3.3, $A\!=\!lim(A_{n},~\phi_{n,m})$ is shape equivalent to $\widetilde{A}\!=\!lim(A_{l_{n}},~\widetilde{\phi}_{n,m})$ and $B\!=\!lim(B_{n},~\psi_{n,m})$ is shape equivalent to $\widetilde{B}\!=\!lim(B_{k},~\widetilde{\psi}_{n,m})$. Evidently the shape equivalences induce the isomorphisms of $$(\underline{K}(A),\underline{K}(A)^{+},\Sigma A)\cong(\underline{K}(\widetilde{A}),\underline{K}(\widetilde{A})^{+},\Sigma \widetilde{A})~~\mbox{and}~~(\underline{K}(B),\underline{K}(B)^{+},\Sigma B)\cong(\underline{K}(\widetilde{B}),\underline{K}(\widetilde{B})^{+},\Sigma \widetilde{B}).$$ It follows that  the
isomorphism $\alpha$ induces an isomorphism $$\widetilde{\alpha}: (\underline{K}(\widetilde{A}),\underline{K}(\widetilde{A})^{+},\Sigma \widetilde{A})\longrightarrow(\underline{K}(\widetilde{B}),\underline{K}(\widetilde{B})^{+},\Sigma \widetilde{B}).$$ By Theorem 9.1 of [DG] and its proof, $\widetilde{A}\cong\widetilde{B}$ and the isomorphism is induced by the shape equivalence between $\widetilde{A}\!=\!lim(A_{l_{n}},~\widetilde{\phi}_{n,m})$
and $\widetilde{B}\!=\!lim(B_{k_{n}},~\widetilde{\psi}_{n,m})$. Hence we have the following theorem (passing to subsequence if necessary).

\noindent\textbf{Theorem 3.5.}~~ Let $A\!=\!lim(A_{n},~\phi_{n,m})$ and $B\!=\!lim(B_{n},~\psi_{n,m})$ be two $A\mathcal{HD}$ inductive limit systems with the ideal property. Let $$\alpha: (\underline{K}(A),\underline{K}(A)^{+},\Sigma A)\longrightarrow(\underline{K}(B),\underline{K}(B)^{+},\Sigma B)$$ be
an isomorphism. Then there are subsequences $$l_{1}<l_{2}<l_{3}<\cdots~~\mbox{ and }~~k_{1}<k_{2}<k_{3}<\cdots$$ and homomorphisms
$\nu_{n}: A_{l_{n}}\longrightarrow B_{k_{n}}~~~\mbox{and}~~~
 \mu_{n}: B_{k_{n}}\longrightarrow A_{l_{n+1}}$ such that\\
(a)~~ $\mu_{n}\circ\nu_{n}\sim_{h}\phi_{l_{n}, l_{n+1}}$ and  $\nu_{n+1}\circ\mu_{n}\sim_{h}\psi_{k_{n},k_{n+1}}.$\\
(b)~~ $\alpha\circ\underline{K}(\phi_{l_{n},\infty})=\underline{K}(\psi_{k_{n},\infty})\circ\underline{K}(\nu_{n})$ and
$\alpha^{-1}\circ\underline{K}(\psi_{k_{n},\infty})=\underline{K}(\phi_{l_{n+1},\infty})\circ\underline{K}(\mu_{n})$, where $\underline{K}(\phi)$ is the map on $\underline{K}$ induced by homomorphism $\phi$.

\noindent\textbf{Remark 3.6.} ~~It is easy to see that for any family of mutually orthogonal projections $\{p_{i} \}\subset A_{l_{n}}$, there is a unitary
$u\in A_{l_{n+1}}$ such that for each $p_{i}$,
$\phi_{l_{n},l_{n+1}}(p_{i})=u^{*}(\mu_{n}\circ\nu_{n})(p_{i})u,$  Hence we can modify $\mu_{n}$ to $Adu\circ\mu_{n}$ to obtain
$\phi_{l_{n},l_{n+1}}(p_{i})=(\mu_{n}\circ\nu_{n})(p_{i}).$

\vspace{0.2in}

\noindent\textbf{\S4. The spectral distribution property and the decomposition theorems}

In the case of simple inductive limit $A\!=\!lim(A_{n},~\phi_{n,m})$ with $A_{n}\in \mathcal{HD}$ (and $\phi_{n,m}$ injective), one can prove that for each $A_{n}$ and each $\eta>0$, there is an $m>n$ and $\delta>0$, such that all partial maps $\phi_{n,m}^{i,j}: A_{n}^{i}\longrightarrow \phi_{n,m}^{i,j}(\textbf{1}_{A_{n}^{i}})A_{m}^{j}\phi_{n,m}^{i,j}(\textbf{1}_{A_{n}^{i}})$ have property $sdp(\eta,\delta)$ (see 1.24 for the definition of $sdp(\eta,\delta)$). For the case of $AH$ inductive limit, this is proved in [Li2] and for the case of inductive limit of dimension drop $C^{*}$-algebras, this is proved in [EGJS].) For the case of non simple inductive limits with the ideal property, $\phi_{n,m}^{i,j}$ itself will not have
property $sdp(\eta,\delta)$ in general no matter how large the $m$ is---we will prove that one can factor the map $\phi_{n,m}^{i,j}$ as
$A_{n}^{i}\rightarrow\bigoplus\limits_{s}(A_{n}^{i}|_{Y_{i}^{j,s}})\xrightarrow{\oplus\phi^{s}} A_{m}^{j},$ where $Y_{i}^{j,s}\subset X_{n,i}$ are connected simplicial sub-complexes and $\phi^{s}$ has property $sdp(\eta,\delta)$. For the case of $AH$ inductive limit, this is Lemma 2.8 of [GJLP1]. The general case
(which also involves dimension drop interval algebras) can be proved similarly (see below). We will use such theorem to prove several decomposition
theorems, by decomposing the connecting homomorphisms into several parts, with the major part factoring through interval algebras.

\noindent\textbf{4.1.}~~ Let $A=PM_{k}(C(x))P$ or $M_{l}(I_{k})$. Let $Y\subset Sp(A)$ be a closed subspace. As in 1.31, recall that $A|_{Y}\triangleq\{f|_{Y},f\in A\}$
which is a quotient algebra of $A$.

\noindent\textbf{Proposition 4.2.}~~ (c.f. Theorem 2.2 of [GJLP1])  Let $A\!=\!lim(A_{n}=\oplus_{i=1}^{t_{n}}A_{n}^{i},~\phi_{n,m})$ be an $A\mathcal{HD}$ inductive limit algebra with the ideal property. For any $A_{n},\eta>0$, let $\kappa$ be a simplicial decomposition of $\coprod_{i=1}^{t_n} X_{n,i}$ such that every simplex of $(\coprod_{i=1}^{t_n} X_{n,i},\kappa)$ has diameter less than $\frac{\eta}4$. Then there is an $m$,  and there are sub-complexes $X_{n,i}^{1},X_{n,i}^{2},\cdots,X_{n,i}^{t_{m}}\subset Sp(A_{n}^{i})$ of $(X_{n,i}, \kappa)$ such that the following hold.\\
(1)~~ $Sp(\phi_{n,m}^{i,j})\subset X_{n,i}^{j}$;\\
(2)~~ For any $x_{0}\in X_{n,i}^{j}$ and $y\in Sp(A_{m}^{j})$, $$Sp(\phi_{n,m}^{i,j})_{y}\cap B_{\eta/2}(x_{0},X_{n,i}^{j})\neq\phi,$$ where
$B_{\eta/2}(x_{0},X_{n,i}^{j})=\{x\in X_{n,i}^{j}~|~~dist(x,x_{0})<\eta/2\}$.
\begin{proof}
The proof is completely the same as the proof of Theorem 2.2 of [GJLP1] replacing Lemma 2.9 of [Pa1] by Proposition 3.1 above. Notice that, in the construction of sub-complexes $\{X_{n,i}^j\}_j$ of $X_{n,i}$ for certain simplicial structure $\kappa$ in the proof of  Theorem 2.2 of  [GJLP1], we only use the condition that every simplex of $\kappa$ has diameter less than $\frac{\eta}4$.\\
\end{proof}
\noindent\textbf{4.3.}~~ In the above proposition, $X_{n,i}^{j}(1\leq j\leq t_{m})$ may not be connected. We need to reduce the version which involves only connected sub-complexes. One can write $$X_{n,i}^{j}=Z_{1}\cup Z_{2}\cdots\cup Z_{\bullet},$$ where $Z_{k}$ are connected simplicial sub-complexes of $X_{n,i}^{j}$.
But the fact that
\begin{center}
$Sp(\phi_{n,m}^{i,j})_{y}\cap B_{\eta/2}(x_{0},X_{n,i}^{j})\neq\emptyset$~~~~~~$\forall x_{0}\in X_{n,i}^{j}$
\end{center}
does not imply
\begin{center}
$Sp(\phi_{n,m}^{i,j})_{y}\cap B_{\eta/2}(x_{0},Z_{k})\neq\emptyset$~~~~~~$\forall x_{0}\in Z_{k}$
\end{center}
 because even for $x_{0}\in Z_{k}$, the set $B_{\eta/2}(x_{0},Z_{k})$ could
be a proper subset of $B_{\eta/2}(x_{0},X_{n,i}^{j})$. In the paper [GJLP1], we use a technique to prove that $Sp(\phi_{n,m}^{i,j})_{y}\cap B_{\eta}(x_{0},Z_{k})\neq\emptyset$ (use $\eta$ instead of $\eta/2$). This technique works equally well for $M_{l}(I_{k})$ with $SpM_{l}(I_{k})=[0,1]$. In fact, the case that $A\!=\!M_{l}(I_{k})$ is similar to the interval algebras $A\!=\!M_{n}(C[0,1])$, which has been dealt with in [Ji-Jiang] by a different, but more
straight forward way. That is, we have the following theorem.

\noindent\textbf{Theorem 4.4.}~ (c.f Lemma 2.8 of [GJLP1])
  Let $A\!=\!lim(A_{n}=\oplus_{i=1}^{t_{n}}A_{n}^{i},~\phi_{n,m})$ be an $A\mathcal{HD}$ inductive limit with the ideal property. For $A_{n},\eta>0$,  let $\kappa$ be a simplicial decomposition of $\coprod_{i=1}^{t_n} X_{n,i}$ such that every simplex of $(\coprod_{i=1}^{t_n} X_{n,i},\kappa)$ has diameter less than $\frac{\eta}{128}$. Then there exist a $\delta>0$, a positive integer $m>n$, connected finite simplicial complexes $Z_{i}^{1},Z_{i}^{2},\cdots,Z_{i}^{\bullet}\subset (X_{n,i}, \kappa),  i=1,2,\cdots,t_{n}$ and a homomorphism $$\phi: B=\oplus_{i=1}^{t_{n}}(\oplus_{s}A_{n}^{i}|_{Z_{i}^{s}})\longrightarrow A_{m}$$ such that\\
(1)~~ $\phi_{n,m}$ factors through $B$ as $A_{n}\xrightarrow {\pi} B\xrightarrow{\phi} A_{m}$, where $\pi$ is defined by
$$\pi(f)=(f|_{Z_{i}^{1}},f|_{Z_{i}^{2}},\cdots,f|_{Z_{i}^{\bullet}})\in\oplus_{s}A_{n}^{i}|_{Z_{i}^{s}}\subset B$$ for any $f\in A_{n}^{i}$;\\
(2)~~ The homomorphism $\phi$ satisfies the dichotomy condition

$(*)$: for each $Z_{i}^{s}$, the partial map
$$\phi^{(i,s),j}: A_{n}^{i}|_{Z_{i}^{s}}\longrightarrow A_{m}^{j}$$ either has the property $sdp(\eta/32,\delta)$ or it is a zero map. Furthermore,  for any
$m'>m$, each partial map of $\phi_{m,m'}\circ \phi$ satisfies the dichotomy condition $(*)$: either it has property $sdp(\eta/32,\delta)$ or it is a
zero map.

\noindent\textbf{4.5.}~~ Recall $KK(I_{k_{1}},I_{k_{2}})=\mathbb{Z}\oplus\mathbb{Z}/k_{1}\oplus\mathbb{Z}/gcf(k_{1},k_{2})$, where $gcf(k_{1},k_{2})$ is the greatest common factor of $k_{1}$ and $k_{2}$. Let $\psi_{0},\psi_{1}: I_{k_{1}}\longrightarrow I_{k_{2}}$ be defined by
$$\psi_{0}(f)=f(\underline{0})\cdot1_{k_{2}},~~ ~~~\psi_{1}(f)=f(\underline{1})\cdot1_{k_{2}}
~~~~~\forall f\in I_{k_{1}}.$$
Let $\psi: I_{k_{1}}\longrightarrow M_{k_{1}/gcf(k_{1},k_{2})}(I_{k_{2}})$ be defined by
$$\psi(f)=diag(\underbrace{f,f,\cdots,f}\limits_{k_{2}/gcf(k_{1},k_{2})}) ~~~~~~\forall f\in I_{k_{1}}.$$
As in Theorem 3.1 of [EGJS], in $KK(I_{k_{1}},I_{k_{2}})(=\mathbb{Z}\oplus\mathbb{Z}/k_{1}\oplus\mathbb{Z}/gcf(k_{1},k_{2}))$, $[\psi_{0}]$ can be identified with $(1,0,0)$, $[\psi_{1}]-[\psi_{0}]$ can be identified with $(0,1,0)$, and $[\psi]-k_{1}/gcf(k_{1},k_{2})[\psi_{0}]$ can be identified with $(0,0,1)$.

\noindent\textbf{4.6.}~~ Let $e_{0},e_{1}: M_{l_{2}}(I_{k_{2}})\longrightarrow M_{l_{2}}(\mathbb{C})$ be the irriducuble  representations defined by evaluations at points 0 and 1 respectively.
That is, $e_{0}(f)=f(\underline{0})=\pi_{\underline{0}}(f)~~~~\mbox{and}~~~~e_{1}(f)=f(\underline{1})=\pi_{\underline{1}}(f).$ \\Let
$\phi: M_{l_{1}}(I_{k_{1}})\longrightarrow M_{l_{2}}(I_{k_{2}})$ be a unital homomorphism. As in 1.31 and 1.32, we use $\phi|_{\underline{0}}$ and
$\phi|_{\underline{1}}$ to denote the homomorphisms $e_{0}\circ\phi$ and $e_{1}\circ\phi$ respectively. Then $$KK(\phi|_{\underline{0}})\in KK(I_{k_{1}},\mathbb{C})=\mathbb{Z}\oplus\mathbb{Z}/k_{1}~~~~\mbox{ and}~~~ KK(\phi|_{\underline{1}})\in KK(I_{k_{1}},\mathbb{C})=\mathbb{Z}\oplus\mathbb{Z}/k_{1}.$$
Let $\frac{n_{0}}{k_{1}}$ ($\frac{n_{1}}{k_{1}}$ respectively) be the number of $Sp\phi|_{\{\underline{0}\}}\cap\{1\}$ ($Sp\phi|_{\{\underline{1}\}}\cap\{1\}$, respectively) counting multiplicity. Namely,
 $$\frac{n_{0}}{k_{1}}=\sharp(Sp\phi|_{\{\underline{0}\}}\cap\{1\})~~~~\mbox{and }~~~~\frac{n_{1}}{k_{1}}=\sharp(Sp\phi|_{\{\underline{1}\}}\cap\{1\})$$ counting multiplicity. Then $$KK(\phi|_{\{0\}})=(\frac{l_{2}}{l_{1}},n_{0})\in
\mathbb{Z}\oplus\mathbb{Z}/k_{1} ~~~~~\mbox{and}~~~~~KK(\phi|_{\{1\}})=(\frac{l_{2}}{l_{1}},n_{1})\in\mathbb{Z}\oplus\mathbb{Z}/k_{1}$$---of course $n_{0}$ and $n_{1}$
are considered to be integers modulo $k_{1}$. (Note that the multiplicity of 0 or 1 in $Sp\phi_{t}$ may appear fractional times.)

\noindent\textbf{Lemma 4.7.}~~ $n_{1}$-$n_{0}$ is a multiple of $k_{1}/gcf(k_{1},k_{2})$ (we denote $\widetilde{k}=k_{1}/gcf(k_{1},k_{2})$).
\begin{proof}
Recall that $\phi|_{0}$ (and $\phi|_{1}$ respectively): $M_{l_{1}}(I_{k_{1}})\longrightarrow M_{l_{2}k_{2}}(\mathbb{C})$ are defined by
$$\phi|_{0}=diag(\underbrace{\phi|_{\underline{0}},\phi|_{\underline{0}},\cdots,\phi|_{\underline{0}}}\limits_{k_{2}})~~~~\mbox{ and
}~~~~\phi|_{1}=diag(\underbrace{\phi|_{\underline{1}},\phi|_{\underline{1}},\cdots,\phi|_{\underline{1}}}\limits_{k_{2}}),$$ respectively. Then $\phi|_{0}$
is homotopic to $\phi|_{1}$ via $\phi|_{t}, 0<t<1$. Consequently,  $$k_{2}(n_{1}-n_{0})=0~~(mod\;k_{1})$$ as desired.\\
\end{proof}
It is easy to see
$$KK(\phi)=(\frac{l_{2}}{l_{1}},n_{0},\frac{n_{1}-n_{0}}{\widetilde{k}})\in\mathbb{Z}\oplus\mathbb{Z}/k_{1}\oplus\mathbb{Z}/gcf(k_{1},k_{2})$$ under the identification in 4.5.

\noindent\textbf{4.8.}~~ Suppose that a unital homomorphism $$\phi: M_{l_{1}}(I_{k_{1}})\longrightarrow M_{m}(\mathbb{C})$$ satisfies $$KK(\phi)=(m/l_{1},m_{1})\in\mathbb{Z}\oplus\mathbb{Z}/k_{1},$$ where $0\leq m_{1}<k_{1}$. Let $$m_{0}=\frac{m}{l_{1}}-m_{1}~~~~~(mod\;k_{1})$$ and $$0\leq m_{0}<k_{1}.$$
Then up to unitary equivalence, $\phi$ can be written as $\phi_{0}\oplus\phi'\oplus\phi_{1}$, where $$\phi_{0}(f)=diag(\underbrace{{f}(\underline{0}),{f}(\underline{0}),\cdots,{f}(\underline{0})}\limits_{m_{0}}),~~~~\phi_{1}(f)=diag(\underbrace{{f}(\underline{1}),{f}(\underline{1}),\cdots,{f}(\underline{1})}\limits_{m_{1}}),$$ and $\phi'$ factors through $M_{l_{l}k_{1}}C[0,1]$ as
$$\phi':M_{l_{1}}(I_{k_{1}})\xhookrightarrow{\imath} M_{l_{1}k_{1}}(C[0,1])\xrightarrow{\tilde{\phi}} M_{m-l_{1}(m_{0}+m_{1})}(\mathbb{C}),$$\\
where $\imath$ is the canonical inclusion.

The following proposition is a combination of [EGJS, 2.3-2.4] and the proof  of [EGJS,  Theorem 3.7]. In fact, it is essentially proved in [Ell1-2] and [Su].

\noindent\textbf{Proposition 4.9.} Let $\phi,\psi:M_{l_{1}}(I_{k_{1}})\rightarrow M_{l_{2}}(I_{k_{2}})$ be two unital homomorphisms. Suppose that $Sp\phi|_{\underline{i}}$ and $Sp\psi|_{\underline{i}}~ (i=0\ or\ 1)$ can be paired within $\eta$ and suppose that for any $t\in(0,1)$, $Sp\phi|_{t}$ and $Sp\psi|_{t}$ can be paired within $\eta<1$. Then for any finite subset $F\subset M_{l_{1}}(I_{k_1})$, there is a unitary $u\in M_{l_{2}}(I_{k_{2}})$, such that $\parallel\phi(f)-u\psi(f)u^{\ast}\parallel\leq sup\{\|f(x)-f(y)\|,|x-y|<\eta\}$.

That is,  if $F$ is $(\varepsilon,\eta)$ uniformly continuous, then
$\|\phi(f)-u\psi(f)u^{\ast}\|\leq\varepsilon~~$ for all $f\in F.$
 \begin{proof}
 Applying Theorem 2.3 of [EGJS], one can reduce the general case to the case that both $\phi$ and $\psi$ are of standard forms as described in Definition 2.2 of [EGJS]. Then the proposition follows from the proof of Theorem 3.7 of [EGJS] (also see Theorem 2.4 of [EGJS] and its proof).\\
\end{proof}

\noindent\textbf{Remark 4.10.}
(a) We do not need to assume $KK(\phi)=KK(\psi)$ since it follows from the condition that $Sp\phi|\underline{_{i}}$ and $Sp\psi|\underline{_{i}} (i=0,1)$ can be paired within $\eta<1$.

(b) Also in the proposition, one can replace the condition that $Sp\phi|\underline{_{i}}$ and $Sp\psi|\underline{_{i}} (i=0,1)$ can be paired within $\eta$ by $KK(\phi)=KK(\psi)$. In fact, the conditions $KK(\phi)=KK(\psi)$ and $Sp\phi|_{t}$ and $Sp\psi|_{t}$ (as $t\rightarrow0$) can be paired within $\eta$ imply that $Sp\phi|_{0}$ and $Sp\psi|_{0}$ can be paired within $\eta'$ for any $\eta'>\eta$. One can see this as follows. To simplify the notation, let $\phi,\psi: I_{k_{1}}\rightarrow M_{l}(I_{k_{2}})$(that is $l_{1}=1$). By $KK(\phi)=KK(\psi)$, there is $l_{0},l_{1}<k_{1}$ \\such that
$$
Sp\phi\underline{_{0}}=\{\underline{0}^{\sim l_{0}},\underline{1}^{\sim l_{1}},t_{1},t_{2},\cdots,t_{m}\}
=\{0^{\sim \frac{l_{0}}{k_{1}}},1^{\sim \frac{l_{1}}{k_{1}}},t_{1},t_{2},\cdots,t_{m}\}~~~\mbox{and}
$$
$$
Sp\psi\underline{_{0}}=\{\underline{0}^{\sim l_{0}},\underline{1}^{\sim l_{1}},t_{1}',t_{2}',\cdots,t_{m}'\}
=\{0^{\sim \frac{l_{0}}{k_{1}}},1^{\sim \frac{l_{1}}{k_{1}}},t_{1}',t_{2}',\cdots,t_{m}'\}
$$
for $t_{i}\in[0,1]$, and $t_{i}'\in[0,1]$. $Sp\phi|_{t}$ and $Sp\psi|_{t}$ can be paired within $\eta$ (for $t\rightarrow0$)
implies that
$$
Sp\phi_{0}=\{0^{\sim l_{0}\cdot\frac{k_{2}}{k_{1}}},1^{\sim l_{1}\cdot\frac{k_{2}}{k_{1}}},t_{1}^{\sim k_{2}},t_{2}^{\sim k_{2}},\cdots,t_{m}^{\sim k_{2}}\}
~~\mbox{and}~~
Sp\psi_{0}=\{0^{\sim l_{0}\cdot\frac{k_{2}}{k_{1}}},1^{\sim l_{1}\cdot\frac{k_{2}}{k_{1}}},t_{1}'^{\sim k_{2}},t_{2}'^{\sim k_{2}},\cdots,t_{m}'^{\sim k_{2}}\}
$$
can be paired within $\eta'(>\eta)$. This is equivalent to that $\{t_{1}^{\sim k_{2}},t_{2}^{\sim k_{2}},\cdots,t_{m}^{\sim k_{2}}\}$ and
$\{t_{1}'^{\sim k_{2}},t_{2}'^{\sim k_{2}},\cdots,t_{m}'^{\sim k_{2}}\}$ can be paired within $\eta'$ which is also equivalent to that $\{t_{1},t_{2},\cdots,t_{m}\}$ and $\{t_{1}',t_{2}',\cdots,t_{m}'\}$ can be paired within $\eta'$.

\noindent\textbf{Theorem 4.11.} Let $A=I_{k_{1}}$, for any positive numbers $\delta$, $\eta<1$, there is a positive number $L$ such that if $l\geq L$ and $\phi:I_{k_{1}}\rightarrow M_{l}(I_{k_{2}})$ is a unital homomorphism with property $sdp(\eta/16,\delta)$, then there are two homomorphisms $\varphi_{1}:I_{k_{1}}\rightarrow M_{l_{1}}(I_{k_{2}})$ and $\varphi_{2}:I_{k_{1}}\rightarrow M_{l_{2}}(I_{k_{2}})$ with $l_{1}+l_{2}=l$, such that for any finite set $F\subset I_{k_{1}}$,\\
(1) $\|\phi(f)-(\varphi_{1}\oplus\varphi_{2})(f)\|\leq sup\{\|f(x)-f(y)\|~|~~ |x-y|<\eta\}, \forall f\in F$;\\
(2) $SPV(\varphi_{1})<\eta$;\\
(3) $\varphi_{2}$ factors through $M_{k_{1}}(C[0,1])$ as~~
$\varphi_{2}:I_{k_{1}}\hookrightarrow M_{k_{1}}(C[0,1])\xrightarrow{\varphi_{2}'} M_{l_{2}}(I_{k_{2}}).$
\begin{proof}
Let $n$ be a positive integer with $\frac{1}{n}<\frac{\eta}{2}$. Let $L$ be an integer satisfying $\delta\cdot L>12nk_1$. Let us assume that $l\geq L$ and that $ \phi:I_{k_{1}}\rightarrow M_{l}(I_{k_{2}})$ is a unital homomorphism. Then $KK(\phi)=(l,d,s)\in\mathbb{Z}\oplus\mathbb{Z}/k_{1}\oplus\mathbb{Z}/(k_{1},k_{2})$(we use $(k_{1},k_{2})$ for the great common factor of $k_{1}$ and $k_{2}$). Let us assume that $d<k_{1}$ and $s<(k_{1},k_{2})$. Let $d_{0}=l-d-s\cdot\frac{k_{1}}{(k_{1},k_{2})}$ (mod $k_{1}$) and $d_{0}<k_{1}$. Set $l_{0}=d+d_{0}+s\cdot\frac{k_{1}}{(k_{1},k_{2})}$, then $l-l_{0}$ is a multiple of $k_{1}$, and $l_{0}<3k_{1}$.

Define $\varphi_{0}:I_{k_{1}}\rightarrow M_{l_{0}}(I_{k_{2}})$ by
$$\varphi_{0}=diag(\psi_{0}^{\sim d_{0}},\psi_{1}^{\sim d},\psi^{\sim s}),$$
where $\psi_{0}$, $\psi_{1}$ and $\psi$ are as in 4.5, which represent $(1, 0, 0), (1, 1, 0)$ and $(\frac{k_{1}}{(k_{1},k_{2})}, 0, 1)$, respectively,  in $KK(I_{k_{1}},I_{k_{2}})=\mathbb{Z}\oplus\mathbb{Z}/k_{1}\oplus\mathbb{Z}/(k_{1},k_{2}).$
Hence $KK(\varphi_{0})=(l_{0},d,s)$. Let $\varphi_{1}':I_{k_{1}}\rightarrow M_{3(n-1)k_{1}}(I_{k_{2}})$ be defined by $$\varphi_{1}'(f)=
diag \left (f(\frac{1}{n})^{\sim3},f(\frac{2}{n})^{\sim3},\cdots,f(\frac{n-1}{n})^{\sim3}\right )\otimes\textbf{1}_{k_{2}}\in M_{3(n-1)k_{1}}(C[0,1])\otimes\textbf{1}_{k_{2}}\subset M_{3(n-1)k_{1}}(I_{k_{2}} ).$$

Let $\varphi_{1}:I_{k_{1}}\rightarrow M_{l_{1}}(I_{k_{2}})$ be defined by $\varphi_{1}=diag(\varphi_{0},\varphi_{1}')$, where $l_{1}=l_{0}+3(n-1)k_{1}$. Since $\{\frac{i}{n}\}_{i=1}^{n-1}$ is $\frac{\eta}{2}$ dense in $[0,1]$, and $l_{0}<3k_{1}$, we have $SPV(\varphi_{1})<\eta$ (see Theorem 1.2.17 of [G5] and its proof).

By Theorem 2.3 of [EGJS], we can assume that $\phi: I_{k_{1}}\rightarrow M_{l}(I_{k_{2}})$ is given by
$$\phi(f)(t)=U diag \left({f}(\underline{0})^{\sim d_{0}k_{2}},f(\lambda_{1}(t)),\cdots,f(\lambda_{m}(t)),{f}(\underline{1})^{\sim dk_{2}} \right)U^{\ast},$$
where $$m=\frac{lk_{2}-(d_{0}+d)k_{2}}{k_{1}}=\frac{(l-l_{1}+3(n-1)k_{1}+s\frac{k_{1}}{(k_{1},k_{2})})\cdot k_{2}}{k_{1}},$$ and $$0\leq\lambda_{1}(t)\leq\lambda_{2}(t)\leq\cdots\leq\lambda_{m}(t).$$

Let $m_{1}=\frac{l- l_{1}}{k_{1}}\cdot k_{2}$. Then $m=m_{1}+3(n-1)k_{2}+s\cdot\frac{k_{2}}{(k_{1},k_{2})}$.

Let
\begin{eqnarray*}
&\mu_{1}(t)=\mu_{2}(t)=\cdots=\mu_{k_{2}}(t)=\lambda_{1}(t), \\
&\mu_{k_{2}+1}(t)=\mu_{k_{2}+2}(t)=\cdots=\mu_{2k_{2}}(t)=\lambda_{k_{2}+1}(t), \\
&\mu_{2k_{2}+1}(t)=\mu_{2k_{2}+2}(t)=\cdots=\mu_{3k_{2}}(t)=\lambda_{2k_{2}+1}(t), \\
&\vdots \\
&\mu_{(\widetilde{m}_{1}-1)k_{2}+1}(t)=\mu_{(\widetilde{m}_{1}-1)k_{2}+2}(t)=\cdots=\mu_{\widetilde{m}_{1}k_{2}}(t)=\lambda_{(\widetilde{m}_{1}-1)k_{2}+1}(t),\\ &\mu_{\widetilde{m}_{1}k_{2}+1}(t)=\mu_{\widetilde{m}_{1}k_{2}+2}(t)=\cdots=\mu_{m}(t)=\lambda_{\widetilde{m}_{1}k_{2}+1}(t),
\end{eqnarray*}
where $\widetilde{m}_{1}=\frac{l-l_{1}}{k_{1}}=\frac{m_{1}}{k_{2}}$. Note that $$m-\widetilde{m}_{1}k_{2}=m-m_{1}=3(n-1)+s\frac{k_{2}}{(k_{1},k_{2})}\leq3(n-1)+k_{2}.$$
Also note that $\phi$ has $sdp(\eta/16,\delta)$ property, so any interval of length $\eta/8$ containing at least $\delta\cdot l\cdot k_{2}(>3(n-1)k_{2})$ elements. It follows that $\{\lambda_{i}(t)\}_{i=1}^{m}$ and $\{\mu_{i}(t)\}_{i=1}^{m}$ can be paired to within $\eta/4$.

Define $\varphi_{2}':M_{k_{1}}(C[0,1])\rightarrow M_{\widetilde{m}_{1}k_{1}}(I_{k_{2}})$ by
$$\varphi_{2}'(f)(t)=diag \big (f(\mu_{k_{2}}(t)),f(\mu_{2k_{2}}(t)),\cdots,f(\mu_{\widetilde{m}_{1}k_{2}}(t))
\big )\otimes\textbf{1}_{k_{2}}\in M_{\widetilde{m}_{1}k_{1}}(C[0,1])\otimes\textbf{1}_{k_{2}}\hookrightarrow M_{\widetilde{m}_{1}k_{1}}(I_{k_{2}}).$$
And let $\varphi_{2}=\varphi_{2}'\circ\imath$, where $\imath:I_{k_{1}}\hookrightarrow M_{k_{1}}(C[0,1])$ is the canonical inclusion. Evidently for each $t\in(0,1)$, $Sp(\varphi_{1}\oplus\varphi_{2})|_{t}$ can be  obtained from the set $$\Theta(t)\triangleq\{\underline{0}^{\sim d_{0}k_{2}},\mu_{1}(t),\mu_{2}(t),\cdots,\mu_{m}(t),\underline{1}^{\sim dk_{2}}\}$$ by replacing $m-\widetilde{m}_{1}k_{2}=3(n-1)k_{2}+s\cdot\frac{k_{2}}{(k_{1},k_{2})}$ elements with the set  $$\{t^{\sim s\cdot\frac{k_{2}}{(k_{1},k_{2})}},(\frac{1}{n})^{\sim3k_{2}},(\frac{2}{n})^{\sim3k_{2}},\cdots,(\frac{n-1}{n})^{\sim3k_{2}}\}.$$

Let $p=m-\widetilde{m}_{1}k_{2}$. Evidently, each interval of length $\eta/4$ contains at least $\delta lk_{2}\geq12nk_{1}k_{2}>4p$ elements in the set $\Theta(t)$.

It is routine to verify the following fact: if $X\subset[0,1]$ is a finite set with multiplicity, such that, any interval $(a,b)\subset[0,1]$ of length $\eta/4$ contains at least $4p$ elements in $X$, and if  $Y$ is obtained from $X$ by replacing some $p$ elements of $X$ by any other $p$ elements from $[0,1]$, then $X$ and $Y$ can be paired within $\eta/2$. Applying this fact to $X=\Theta(t)$, we have that $Sp(\varphi_{1}\oplus\varphi_{2})|_{t}$ can be paired within $\Theta(t)$ to within $\eta/2$. On the other hand, we already have that $\Theta(t)$ can be paired within $Sp\phi|_{t}$ to within $\eta/4$. Hence conclusion (1) holds by  Proposition 4.9 and Remark 4.10(b) up to a unitary equivalence.\\
\end{proof}

\noindent\textbf{Remark 4.12.} Theorem 4.11 also holds for $\phi: I_{k}|_{[0,a]}\rightarrow M_{l}(I_{k})$ or $\phi: I_{k}|_{[b,1]}\rightarrow M_{l}(I_{k})$ with an even easier proof.

Also note that if $F$ is $(\varepsilon,\eta)$ continuous, then the conclusion (2) in Theorem 4.11 implies that $\omega(\varphi_{1}(F))<\varepsilon$.

The following Theorem is an analogy of Theorem 2.12 of [GJLP1],
 but involving dimension drop algebras.

\noindent\textbf{Theorem 4.13.} Let $A\!=\!lim(A_{n}\!=\!\oplus_{i=1}^{tn}A_{n}^{i},~\phi_{n,m})$ be $A\mathcal{HD}$ inductive limit with the ideal property. For any $A_{n}$, finite sets $F=\oplus_{i=1}^{tn}F_{n}^{i}\subset A_{n}$ and $G\subset AffTA_{n}$, and $\varepsilon>0$, there exists $m>n$, and there exist projections $Q_{0},Q_{1},Q_{2}\in A_{m}$ with $Q_{0}+Q_{1}+Q_{2}=\phi_{n,m}(\textbf{1}_{A_{n}})$ and unital map $\psi_{0}\in Map(A_{n},Q_{0}A_{m}Q_{0})_1$ and homomorphisms $\psi_{1}\in Hom(A_{n},Q_{1}A_{m}Q_{1})_1$, $\psi_{2}\in Hom(A_{n},Q_{2}A_{m}Q_{2})_1$ such that the following statements are true.\\
(1) $\parallel\phi_{n, m}(f)-\psi_{0}(f)\oplus\psi_{1}(f)\oplus\psi_{2}(f)\parallel<\varepsilon,~~ \forall f\in F$;\\
(2) $\psi_{1}$ is defined by point evaluation
%(or factoring through) finite dimensional ($C^{\ast}$-algebra),
and $$(\psi_{0}\oplus\psi_{1})(F)\subset (Q_{0}+Q_{1})A_{m}(Q_{0}+Q_{1})$$ is weakly approximately constant to within $\varepsilon$, i.e.,  $\omega((\psi_{0}\oplus\psi_{1})(F))<\varepsilon$;\\
(3) The homomorphism $\psi_{2}$ factors through $C$ --- a direct sum of matrix algebras $M_{\bullet}(C[a,b])$ over intervals $[a,b]$ (where $[a,b]\subset[0,1]$) or $\mathbb{C}$ as
$$\psi_{2}: A_{n}\xrightarrow{\xi_{1}}~ C\xrightarrow{\xi_{2}}~Q_{2}A_{m}Q_{2},$$
where $\xi_{1}$ and $\xi_{2}$ are unital homomorphisms;\\
(4) There exists a homomorphism $\widetilde{\psi}_{0}\in Hom(A_{n},Q_{0}A_{m}Q_{0})_1$ such that the following is true: Let $\psi: A_{n}\rightarrow\phi_{nm}(\textbf{1}_{A_{n}})A_{m}\phi_{nm}(\textbf{1}_{A_{n}})$ be defined by $\psi=\widetilde{\psi}_{0}\oplus\psi_{1}\oplus\psi_{2}$. Then $$KK(\phi_{n,m})=KK(\psi)\in KK(A_{n},A_{m})~~~\mbox{ and}~~~\|AffT\phi_{n,m}(f)-AffT\psi(f)\|<\varepsilon$$ for all
 $f\in G$ (see 2.4 for $AffT\phi$ for a
(possibly non unital) homomorphism $\phi$).
%where the maps at $AffT$ level are  the unital maps regarded as  maps from $AffTA_{n}$ to $AffT(\phi_{nm}(\textbf{1}_{A_{n}})A_{m}\phi_{nm}(\textbf{1}_{A_{n}}))$.

Furthermore,  by enlarging the set $F$, one can assume that the condition (1) implies that
% \begin{center}
 $$\|AffT\phi_{n,m}(g)-AffT\psi'(g)\|<{\varepsilon}/{2}, ~~~~~~\forall g\in G,$$
%\end{center}
where $\psi'=\psi_{0}\oplus\psi_{1}\oplus\psi_{2}$.

(If we further assume that $\frac{rank(Q_{0})}{rank(Q_{1})}<{\varepsilon}/{2}$, ~then
%\begin{center}
$$\|AffT\psi'(g)-AffT\psi(g)\|<{\varepsilon}/{2}~~~~~~\mbox{ for all}~ g ~\mbox{with}~~ \|g\|\leq1.)$$
%\end{center}

\begin{proof}
Using the dilation Lemma --- Lemma 1.9 (see also Lemma 1.11), one can assume that each block $A_{n}^{i}$ is of the form $M_{[n,i]}(C(X_{n,i}))$, $X_{n,i}=S^{1}$, $[0,1]$, $T_{\uppercase\expandafter{\romannumeral2}
,k}$ or of the  form $M_{[n,i]}(I_{k_{n,i}})$. Let $$\varepsilon'=\frac{\varepsilon}{\max_{1\leq i\leq t_{n}}[n,i]}.$$ For $F_{n}^{i}\subset A_{n}^{i}$, let $F'^{i}\subset C(X_{n,i})$ (or $I_{k_{n,i}}$) be the finite set consisting of all the entries of the elements in $F_{n}^{i}(\subset M_{[n,i]}C(X_{n,i}))$ or $M_{[n,i]}(I_{k_{n,i}})$. Let $\eta>0$ be such that if $x,x'\in X_{n,i}$ (or $x,x'\in[0,1]$) with dist$(x,x')<2\eta$, then
%\begin{center}
$\parallel f(x)-f(x^{\prime})\parallel<{\varepsilon'}/{3}~~\forall f\in F^{'i}.$
%\end{center}
Applying Theorem 4.4 to $A_{n}$ and $\eta$, there is a $\delta>0$ as in Theorem 4.4. That is, there are $m_{1}>n$, connected finite simplicial complexes
$Z_{i}^{1},Z_{i}^{2},\cdots,Z_{i}^{\bullet}\subset SpA_{n}^{i}$,  and a homomorphism
$\phi'=\oplus_{i}\oplus_{s}A_{n}^{i}|_{Z_{i}^{s}}\rightarrow A_{m_{1}}$ satisfying  the following conditions:\\
(1) $\phi_{n,m'}=\phi'\circ\pi$, where $\pi:A_{n}\rightarrow\oplus_{i}\oplus_{s}A_{n}^{i}|_{Z_{i}^{s}}$ is defined by the direct  sum of the restriction maps;\\
(2) For each $Z_{i}^{s}\subset Sp(A_{n}^{i})$ (either $X_{n,i}$ or $[0,1]$ for the case of $M_{[n,i]}(I_{k_{ni}})$), the partial map $$\phi'^{(i,s),j}:A_{n}^{i}|_{Z_{i}^{s}}\rightarrow A_{m_{1}}^{j}$$
either has $sdp(\eta/32,\delta)$ property or is zero map.

Let $L$ be the maximum of the number $L$  in Theorem 4.35 of [G5]~(corresponding to $Z_{i}^{s},\eta,\delta$) and the number  $L$  in Theorem 4.11 (see Remark 4.12).

Let $T= \max\limits_{i}[n,i]$, where $[n,i]$ are the sizes of $A_{n}^{i}\!=\!M_{[n,i]}(C(X_{n,i}))$ (or $A_{n}^{i}\!=\!M_{[n,i]}(I_{k_{n,i}})$). Let $$K=125T\cdot L^{2}\cdot2^{L+1}$$ (here $125=(2+2+1)^{3}$ is the largest possible number of $(\dim X+\dim Y+1)^{3}$ when we apply  Theorem 4.35 of [G5] to the spaces $X_{n,i}$).

Let $B=\oplus_{i}\oplus_{s}A_{n}^{i}\mid _{Z_{i}^{s}}$.

By Proposition 3.2, there is $m_{2}>m_{1}$ such that each partial map $\phi^{(i,s),j}: A_{n}^{i}|_{Z_{i}^{s}}\rightarrow A_{m_{2}}^{j}$ of $\phi=\phi_{m_{1},m_{2}}\circ\phi':B\rightarrow A_{m_{2}}$ satisfies one of the following two conditions:\\
(a) There is a unital homomorphism $\psi^{(i,s),j}:A_{n}^{i}|_{Z_{i}^{s}}\rightarrow\phi^{(i,s),j}(1)A_{m}^{j}\phi^{(i,s),j}(1)$ with finite dimensional image, such that $$\parallel\phi^{(i,s),j}(f)-\psi^{(i,s),j}(f)\parallel<\varepsilon$$
for all $f\in F^{(i,s)}=\pi_{s}(F^{i})$, where $\pi_{s}:A_{n}^{i}\rightarrow A_{n}^{i}|_{Z_{i}^{s}}$ is the restriction map; or\\
(b) $\phi^{(i,s),j}$ is $K$-large.

 We need the following observation: In the decomposition of $\phi:C(X)\rightarrow PM_{k}C(Y)P$ in Theorem 4.35 of [G5], if $\dim(X)\leq2,\dim(Y)\leq2$ with $H^{2}(X)$ and $H^{2}(Y)$ torsion groups as in our case (note that for a proper  sub-complex $X$ of $T_{\uppercase\expandafter{\romannumeral2}
,k}$, one has that  $H^{2}(X)=0)$, and $rank(Q_{0})\geqslant12$, then there is a homomorphism $\widetilde{\phi}_{0}\in Hom(C(X),Q_{0}M_{k}(C(Y))Q_{0})_1$ such that $KK(\phi)=KK(\widetilde{\phi}_{0}\oplus\phi_{1}\oplus\phi_{2})$. This is a consequence of Theorem 3.22 of [EG2]. (If rank$(Q_{0})<12$, we can move a part of point evaluation from  $\phi_{1}$ into $\phi_{0}$ to make rank$(Q_{0})\geq12)$.

Applying Theorem 2.12 of [GJLP1] to a homogenous copy of $A_{n}^{i}$ or applying Theorem 4.11 to a dimension drop copy of $A_{n}^{i}$, and using the above observation, we have the following claim.

\noindent\textbf{Claim:} Let $A=\lim(A_{n}\!=\!\oplus A_{n}^{i},~\phi_{n, m})$ be an $A\mathcal{HD}$ inductive limit with the ideal property. For any finite set $F=\oplus F_{n}^{i}\subset A_{n}^{i}$ and $\varepsilon>0$, there is an $m_{2}>0$, such that the following conclusions (1) and (2)  hold.

(1) For each block $A_{n}^{i}=M_{[n,i]}C(X_{n,i})$, and each $A_{m_{2}}^{j}$ (no matter $A_{m_{2}}^{j}$ is of the homogenous type or of the dimension drop type) if $\psi: A_{m_{2}}^{j}\rightarrow M_{\bullet}(C(Z))$ is a homomorphism, then
\begin{center}
$\psi\circ\phi_{n,m}^{i,j}:A_{n}^{i}\rightarrow (\psi\circ\phi_{n,m}^{i,j}(1))M_{\bullet}(C(Z)) (\psi\circ\phi_{n,m}^{i,j}(1))$
\end{center}
has the desired decomposition as in our theorem. That is, there are projections
% \begin{center}
 $Q_{0},Q_{1},Q_{2}\in M_{\bullet}(C(Z))$
%\end{center}
with
% \begin{center}
 $Q_{0}+Q_{1}+Q_{2}=\psi\circ\phi_{n,m}^{i,j}(\textbf{1})$
%\end{center}
and
% \begin{center}
 $\psi_{0}\in Map(A_{n}^{i},Q_{0}M_{\bullet}C(Z)Q_{0})_1$,
%\end{center}
% \begin{center}
 $\psi_{1}\in Hom(A_{n}^{i},Q_{1}M_{\bullet}C(Z)Q_{1})_1$,
%\end{center}
and
%\begin{center}
$\psi_{2}\in Hom(A_{n}^{i},Q_{2}M_{\bullet}C(Z)Q_{2})_1$
%\end{center}
 such that $$\parallel\psi\circ\phi_{n,m_{2}}^{i,j}(f)-(\psi_{0}(f)\oplus\psi_{1}(f)\oplus\psi_{2}(f))\parallel<\varepsilon, ~~\forall f\in F_{n}^{i},$$
where $(\psi_{0}\oplus\psi_{1})(F_{n}^{i})$ is weakly approximately constant to within $\varepsilon$, and $\psi_{2}$ factors through a direct sum of matrix algebras over intervals or points. (In particular, if $A_{m_{2}}^{j}$ is of homogeneous type, then this part of claim implies that we can get desired decomposition for $\phi_{n,m_{2}}^{i,j}$ itself by choosing $\psi$ to be identity). Note that we use the above observation to get required property (4) in our Theorem.

(2) For each block $A_{n}^{i}$ of dimension drop type, and  each $A_{m_{2}}^{j}$ (no matter $A_{m_{2}}^{j}$ is of dimension drop type or of homogeneous type) if $\psi:A_{m_{2}}^{j}\rightarrow M_{\bullet}(I_{k})$ is a homomorphism, then $\psi\circ\phi_{n,m_{2}}^{i,j}$ can be decomposed as desired in the theorem. Note for this case the property (4) already holds, since all parts $\psi_{0},\psi_{1},\psi_{2}$ of the decomposition are homomorphisms.

We apply the claim again to $A_{m_{2}}$ (in place of $A_{n}$) and $F'=\oplus_{j}\oplus_{i}\phi_{n,m}^{i,j}(F_{n}^{i})\subset\oplus A_{m_{2}}^{j}$ (in place of $F$) and $\varepsilon>0$, to obtain $A_{m}$ (in place of $A_{m_{2}}$). Then it is routine to verify that $A_{m}$ and $\phi_{n, m}$ satisfy the condition by decomposing $\phi_{n,m}^{i,j}=\bigoplus\limits_{k}\phi_{m_{2},m}^{k,j}\circ\phi_{n,m_{2}}^{i,k}$ --- that is, if $A_{n}^{i}$ and $A_{m_{2}}^{k}$ are of the same type (both homogeneous or both dimension drop), then $\phi_{n,m_2}^{i,k}$ can be decomposed with desired approximation on $F$, if $A_{m_{2}}^{k}$ and $A_{m}^{j}$ are of the same type, then we decompose $\phi_{m_{2},k}^{m,j}$, and finally if $A_{n}^{i}$ and $A_{m}^{j}$ are of the same type, then we decompose the composition $\phi_{m_{2},m}^{k,j}\circ\phi_{n,m_{2}}^{i,k}$.\\
\end{proof}
The following is more or less Theorem 3.7 of [EGJS] which is the uniqueness theorem involving dimension drop algebras (see also [Ell2]).

\noindent\textbf{Proposition 4.14.} Let $\eta<1$ and $\delta$ be the positive numbers. There exists a finite set $H\subseteq C_{\mathbb{R}}[0,1]$ such that the following is true. Let $\phi_{0},\phi_{1}:A\!=\!M_{l_{1}}(I_{k_{1}})\rightarrow B=M_{l_{2}}(I_{k_{2}})$ be two unital homomorphisms satisfying\\
(i) $[\phi_{0}]=[\phi_{1}]\in KK(A,B)$;\\
(ii) $\phi_{0}$ has $sdp(\frac{\eta}{8},\delta)$ property;\\
(iii) $\parallel AffT\phi_{0}(h)-AffT\phi_{1}(h)\parallel<\frac{\delta}{4}, ~~\forall h\in H$.\\Then there is a unitary $W\in B$ such that $$\parallel\phi_{1}(f)-W\phi_{0}(f)W^{\ast}\parallel\leq3\sup\{\|f(x)-f(y)\|~|~~|x-y|<\eta\}.$$ That is, for a fixed $F\subset A$, $\varepsilon>0$, if $F$ is $(\frac{\varepsilon}{3},\eta)$ uniformly continuous, then one can choose $W$ such that
%\begin{center}
$$\parallel\phi_{1}(f)-W\phi_{0}(f)W^{\ast}\parallel\leq\varepsilon ~~~~~~\forall ~f\in F.$$
%\end{center}

\begin{proof}
 One separates the fractional parts of $0$ and $1$ out from the $Sp\phi_{0}\mid_{t}$ and $Sp\phi_{1}\mid_{t}$ (which are the same since $KK(\phi_{0})=KK(\phi_{1})$. Then it follows from the same proof of [Li 2, 2.15], that $Sp\phi_{0}\mid_{t}$ and $Sp\phi_{1}\mid_{t}$ can be paired within $\eta$ for each $t\in(0,1)$). The proposition follows from Proposition 4.9 and Remark 4.10.\\
\end{proof}
\noindent\textbf{Remark 4.15.} In the above proposition, using the same set $H\subset C_{\mathbb{R}}[0,1]$ which only depends on $\eta$ and $\delta$, one has the following conclusion. For any closed interval $X\subset[0,1]$ ($X$ may be of forms $[0, b],[a,1]$ or $[a,b]$), if $\phi_{0},\phi_{1}:A\mid_{X}\rightarrow B=M_{l_{2}}(I_{k_{2}})$ such that\\
(1) $KK(\phi_{0})=KK(\phi_{1})$,\\
(2) $\phi_{0}$ has $sdp(\eta/8,\delta)$ property (with respect to $X$), and\\
(3) $\parallel AffT\phi_{0}(h\mid_{X})-AffT\phi_{1}(h\mid_{X})\parallel<\frac{\delta}{4}$,\\
then the conclusion of Proposition 4.14 still holds.

\vspace{0.2in}

%\noindent\textbf{\S5. Additional decomposition theorems}

The following theorems are proved in [JLW]

\noindent\textbf{Theorem 4.16.}~(Theorem 3.1 in [JLW])~ Let $F\subset I_{k}$ be a finite set, $\varepsilon > 0$.
There is an $\eta>0$, satisfying that if\\
 $\phi: I_{k}\rightarrow PM_{\bullet}(C(X))P~~(dim(X)\leq 2)$ is
a unital homomorphism such that for any $x\in X$,
\begin{center}
$\sharp(Sp\phi_{x}'\cap[0,\frac{\eta}4])\geq k$ and $\sharp(Sp\phi_{x}'\cap[1-\frac{\eta}4],1])\geq k,$
\end{center}
where
$\phi':C[0,1]\xrightarrow{\imath} I_{k}\xrightarrow{\phi}PM(C(X))P,$\\
then there are three mutually orthogonal projections ~$Q_{0},~Q_{1},~P_{1}~\in PM_{\bullet}(C(X))P$ with $Q_{0}+Q_{1}+P_{1}=P$
and a unital homomorphism $\psi_{1}:~M_{k}(C[0,1])\rightarrow P_{1}M_{\bullet}(C(X))P_{1}$ such that\\
(1) write $\psi(f)={f}(\underline{0})Q_{0}+{f}(\underline{1})Q_{1}+(\psi_{1}\circ\imath)(f),$
then
$\parallel\phi(f)-\psi(f)\parallel<\varepsilon
~~\forall f\in F\subset I_{k}\subset M_{k}(C[0,1]),$ and
 \\
(2) $rank(Q_{0})\leq k$ and $rank(Q_{1})\leq k$.\\

\noindent\textbf{Theorem 4.17.}~(Theorem 4.8 in [JLW])~ Let $X$ be a connected finite simplicial complex of dimension at most 2, $\varepsilon>0$ and $F\subset C(X)$, a finite set of generators. Suppose that $\eta\in(0,\varepsilon)$ satisfies that if $dist(x,x')\leq2\eta$, then $\|f(x)-f(x')\|<\frac{\varepsilon}{4}$ for all $f\in F$.
% (such $\eta$ is chosen in the proof of Theorem 5.30 for $F\subset C(X)$).

For any $\delta>0$ and  positive integer $J>0$, there exist an integer $L>0$ and a finite set $H\subseteq AffTC(X)(=C_{\mathbb{R}}(X))$ such that the following holds.

If $\phi,\psi:C(X)\rightarrow B=M_{K}(I_{k})$ (or $B=PM_{\bullet}(C(Y))P$) are unital homomorphisms with the properties: \\
(a) $\phi$ has $sdp(\eta/32,\delta)$;\\
(b) $K\geq L$ (or $rank(P)\geq L$);\\
(c) $\|AffT\phi(h)-AffT\psi(h)\|<\frac{\delta}{4}$, for all $h\in H$,\\
then there are three orthogonal projections $Q_{0},Q_{1}, Q_2 \in B$,
%two maps $\phi_{0},\psi_{0}\in Map(C(X),$\\$Q_{0}M_{k}(I_{k})Q_{0})_{1}$,
 two homomorphisms $\phi_{1}\in Hom(C(X),Q_{1}BQ_{1})_{1}$ and $\phi_{2}\in Hom(C(X),Q_{2}BQ_{2})_{1}$, and a unitary $u\in B$ such that\\
(1) $\textbf{1}_B=Q_{0}+Q_{1}+Q_2$;\\
(2) $\|\phi(f)-\big(Q_0\phi(f)Q_0+\phi_{1}(f)+\phi_{2}(f)\big)\|<\varepsilon$ ~~and \\ $\|(Adu\circ\psi)(f)-\big(Q_0(Adu\circ\psi)(f)Q_0+\phi_{1}(f)+\phi_{2}(f)\big)\|<\varepsilon$, for all $f\in F$;\\
(3) $\phi_{2}$ factors through $C[0,1]$;\\
(4) $Q_{1}=p_{1}+\cdots+p_{n}$ with $(rank(Q_{0})+2)J<rank(p_{i})$ $(i=1,2,\cdots,n)$, where rank:$K_{0}(B)\rightarrow \mathbb{Z}$ is defined as in 1.12  (which is $rank~ p_{i}(\underline{0})$ for $B=M_K(I_k)$, where $rank~p_{i}(\underline{0})$ is regarded as projections in $M_{K}(\mathbb{C})$ not $M_{K}(M_{k}(\mathbb{C}))$), and $\phi_1$ is defined by
%and $\phi_{0},\psi_{0}$ are defined by
%$$\phi_{0}(f)=p_{0}\phi(f)p_{0}+\phi_{1}(f),\ \psi_{0}(f)=p_{0}Adu\circ\psi(f)p_{0}+\phi_{1}(f),$$
 $$\phi_{1}(f)=\sum\limits_{i=1}\limits^{n}f(x_{i})p_{i},\ \forall f\in C(X),$$ where $p_{1},p_2, \cdots,p_{n}$ are mutually orthogonal projections and
$\{x_{1},x_{2},\cdots,x_{n}\}$ is an $\varepsilon$-dense subset of $X$.\\

\noindent\textbf{Theorem 4.18.}~(Theorem 4.9 in [JLW])~ Let $X$ be a connected finite simplicial complex of dimension at most 2, $\varepsilon>0$ and $F\subset C(X)$, a finite set of generators. Suppose that $\eta\in(0,\varepsilon)$ satisfies that if $dist(x,x')\leq2\eta$, then $\|f(x)-f(x')\|<\frac{\varepsilon}{4}$ for all $f\in F$. Let $\kappa$ be a fixed simplicial structure of $X$.
% (such $\eta$ is chosen in the proof of Theorem 5.30 for $F\subset C(X)$).

For any $\delta>0$ and  positive integer $J>0$, there exist an integer $L>0$ and a finite set $H\subseteq AffTC(X)(=C_{\mathbb{R}}(X))$ such that the following holds.

If $X_1$ is a connected sub-complex of $(X, \kappa)$, and  if $\phi,\psi:C(X_1)\rightarrow B=M_{K}(I_{k})$ (or $B=PM_{\bullet}(C(Y))P$) are unital homomorphisms with the following properties: \\
(a) $\phi$ has $sdp(\eta/32,\delta)$;\\
(b) $K\geq L$ (or $rank(P)\geq L$);\\
(c) $\|AffT\phi(h|_{X_1})-AffT\psi(h|_{X_1})\|<\frac{\delta}{4}$, for all $h\in H$,\\
then there are three orthogonal projections $Q_{0},Q_{1}, Q_2 \in B$, two homomorphisms $\phi_{1}\in Hom(C(X_1),Q_{1}BQ_{1})_{1}$ and $\phi_{2}\in Hom(C(X_1),Q_{2}BQ_{2})_{1}$, and a unitary $u\in B$ such that\\
(1) $\textbf{1}_B=Q_{0}+Q_{1}+Q_2$;\\
(2) $\|\phi(f|_{X_1})-\big(Q_0\phi(f|_{X_1})Q_0+\phi_{1}(f|_{X_1})+\phi_{2}(f|_{X_1})\big)\|<\varepsilon$ ~~and \\ $\|(Adu\circ\psi)(f|_{X_1})-\big(Q_0(Adu\circ\psi)(f|_{X_1})Q_0+\phi_{1}(f|_{X_1})+\phi_{2}(f|_{X_1})\big)\|<\varepsilon$ for all $f\in F$;\\
(3) $\phi_{2}$ factors through $C[0,1]$;\\
(4) $Q_{1}=p_{1}+\cdots+p_{n}$ with $(rank(Q_{0})+2)J<rank(p_{i})$ $(i=1,2,\cdots,n)$,
% where rank:$K_{0}(B)\rightarrow \mathbb{Z}$ is defined as in 1.12,
%which is $rank~ p_{j}(\underline{0})$ for $B=M_K(I_k)$, where $rank~p_{j}(\underline{0})$ is regarded as projections in $M_{K}(\mathbb{C})$ not $M_{K}(M_{k}(\mathbb{C}))$,
and $\phi_1$ is defined by
$$\phi_{1}(f)=\sum\limits_{i=1}\limits^{n}f(x_{i})p_{i}~~~ \forall f\in C(X_1),$$
where $p_{1},p_2, \cdots,p_{n}$ are mutually orthogonal projections and
$\{x_{1},x_{2},\cdots,x_{n}\}$ is an $\varepsilon$-dense subset of $X_1$.

\vspace{3mm}

\noindent\textbf{\S5. The existence theorem }

%\vspace{-2.4mm}

%\vspace{3mm}

In this section, we will present the main existence theorem.

\noindent\textbf{5.1.}~~Recall that a homomorphism $\phi: P M_{l} (C(X)) P\rightarrow B$ (or $\phi: M_{l} (I_{k}) \rightarrow B)$ is called to be defined by point evaluations if $\phi$ can factor through a finite dimensional algebra. For each $x\in X$, we can identify $P(x)M_{l}(\mathbb{C})P(x)\cong M_{rank(P)}(\mathbb{C})$ and regard $f(x)\in M_{rank(P)}(\mathbb{C})$ (the identification is unique up to inner automorphisms). Then $\phi: P M_{l}( C(X))P\rightarrow B$ is defined by point evaluations if there exist $x_{1},x_{2},\cdots,x_{k}\in X$ and homomorphisms $\phi_{i} :M_{rank(P)}(\mathbb{C})\rightarrow B ~(i=1,2,\cdots,k)$
such that
$$\phi(f)=\phi_{1}(f(x_{1}))\oplus \phi_{2}(f(x_{2}))\oplus\cdots \oplus\phi_{k}(f(x_{k})).$$
\noindent\textbf{5.2.}~~Recall from 1.20, if $\phi : M_{l}(\mathbb{C})\rightarrow B$ is a homomorphism, then $\phi=\phi_{1}\otimes id_{l}$, where
$\phi_{1}=\phi \mid_{e_{11}M_{l}(\mathbb{C})e_{11}} : \mathbb{C}\rightarrow \phi(e_{11}) B \phi(e_{11})$
under an identification
$
\phi(\textbf{1}_{M_{l}(\mathbb{C})})B \phi(\textbf{1}_{M_{l}(\mathbb{C})})=(\phi(e_{11}) B \phi(e_{11}))\otimes M_{l}(\mathbb{C}).
$

Let $Q=\phi(\textbf{1}_{M_{l}(\mathbb{C})})$ and $q=\phi(e_{11})$. Then $Q=diag(\underbrace{q,q,\cdots,q}_{l})$ (or denoted by $\underbrace{q\oplus q\oplus \cdots \oplus q}_{l}$) under the identifications $QBQ=M_{l}(\mathbb{C})\otimes qBq$ and $\phi(a)=a\otimes q $ (in the notation in 1.3.4 and 2.6 of [EG2]). We will use $a\cdot Q$ to denote $a\otimes q$ when $Q=diag(\underbrace{q,q,\cdots,q}_{l})$ and $a\in M_{l}(\mathbb{C})$.

Note that if $l=1$, $a\cdot Q$ is the ordinary multiplication.

\noindent\textbf{5.3.}~ For a homomorphism $\phi: PM_{l_1}(C(X))P\rightarrow B$ (with rank(P)=$l$) defined by point evaluations, there are a projection $q\in \phi(P)B\phi(P)$,
an identification $QBQ=M_{l}(\mathbb{C})\otimes qBq $ $(Q\triangleq\phi(P) )$, finitely many points $x_{1},x_{2},\cdots,x_{n}\in X$ (may be repeat), $n$ orthogonal projections $q_{1},q_{2},\cdots,q_{n}\in qBq $ with $\Sigma q_{i}=q $, and a unitary $u\in M_{l}(\mathbb{C})\otimes qBq $,
such that $\phi(f)=u(\sum_{i=1}^{n}f(x_{i})\otimes q_{i})u^{\ast}~~\forall f\in PM_{l}(C(X))P.$\\
In the notation of
5.2, if we let $P_{i}=\textbf{1}_{l}\otimes q_{i}\in M_{l}(\mathbb{C})\otimes qBq$,
then the above map can be written as $\phi(f)=u(\sum_{i=1}^{n}f(x_{i})P_{i})u^{\ast}$.

\noindent\textbf{Proposition 5.4.} Let $A=PM_{l}(C(X))P\in \mathcal{HD}$ or $M_{l}(I_{k})$. For any finite set $H\subseteq AffTA$ and $\varepsilon>0$, there exists an $L > 0$ satisfying following conditions.\\
Suppose that $B=\bigoplus^{s}_{i=1}B^{i}\in \mathcal{HD} $ and $\alpha \in KK(A,B)$, and suppose that $\xi : AffT A\rightarrow AffT B$ is a unital positive map.
Suppose that the following conditions (a) and (b) holds:\\
(a) $\alpha$ is $L-$full, that is, $rank~\alpha^{i}(1)\geq L\cdot rank(P)$
or $rank(\alpha^{i}(1))\geq L\cdot l$ (for $A=M_{l}(I_{k})$), where the function $rank$: $K_{0}(B_{i})\rightarrow \mathbb{Z}$ is defined in 1.13, and $\alpha$ can be realized by unital homomorphisms from $A$ to $B$.\\
(b) $\alpha$ and $\xi$ are compatible in the sense that $\xi\circ \rho_{A}(x)=\rho_{B}(\alpha_{\ast}(x))$, $\forall x\in K_{0}(A)$.

It follows that there exists a unital homomorphism $\psi: A\rightarrow B$ such that\\
1. $KK(\psi)=\alpha$.\\
2. $\|AffT\psi(f)-\xi(f)\|<\varepsilon$~ $\forall f \in H$.

\noindent\textbf{Proof.} One can reduce the proof to the case that B is a single block. For $A=PM_{l}C(X)P$ and $B$ is homogeneous, this is  [EGL, Theorem 3.15].

For the case that at least one of $A$ and $B$ is a dimension drop interval algebra, the proof is similar to that of [EGL, 3.15]. Note that, instead of using [EG2, \S 3], we  use [DG, Theorem 5.7] to realize the KK element $\alpha$ by a homomorphism.  Also  note that $AffT(M_{l}(I_{k}))=AffTC[0,1]=C_{\mathbb{R}}[0,1]$ (see the existence theorem in [EGJS, \S3]).
~~~~~~~~~~~~~~~~~~~~~~~~~~~~~~~~~~~~~~~~~~~~~~~~~~~~~~~~~~~~~~~~~~~~~~~~~~~~~~~~~~~~~~~~~~~~~~~~~~~~~~~~~~~$\Box$

\noindent\textbf{Remark 5.5.} Similar to Remark 3.16 of [EGL], for the case $A=M_{l}(C(S^{1}))$ and $B^{i}=M_{m}(I_{k})$, the homomorphism $\psi^{i}:~A\rightarrow B^{i}$ can be chosen to be of  the form $\psi^{i}=\psi^{i}_{1}\oplus \sigma^{i}$, where $\psi^{i}_{1}:~ A\rightarrow M_{m-l}(I_{k})$ and $\sigma^{i}:A\rightarrow M_{l}(I_{k})$ with $\sigma^{i}$ as below.

There is a continuous map $\gamma: [0,1]\rightarrow S^{1}$ such that
\begin{center}
$\sigma^{i}(f)(t)=f(\gamma(t))\otimes \textbf{1}_{k}\in M_{l}(\mathbb{C})\otimes \textbf{1}_{k}\subseteq M_{l_{k}}(\mathbb{C})$~~~ $\forall f\in M_{l}C((S^{1}))$.
\end{center}
 Furthermore, the same argument as in [EGL, 3.16] prove that $\gamma$ can be  replaced by any continuous map from $[0,1]$ to $S^{1}$.

\noindent\textbf{Proposition 5.6.} Suppose that $A=PM_{l}(\mathbb{C}(X))P$ (or $M_{l}(I_{k})$), $H(\subset AffTA)$  and $E(\subset U(A)/ \widetilde{SU}(A))$  are any finite sets (see \S2 for the notation). For any $\varepsilon>0,$ there is a positive integer $L$ satisfying the following conditions. Let $B=\bigoplus^{s}_{i=1}B^{i}\in \mathcal{HD}$ and $\alpha\in KK(A,B)$, and let $\xi: AffTA \rightarrow AffTB$ be a unital positive linear map, $\gamma: U(A)/ \widetilde{SU}(A) \rightarrow U(B) / \widetilde{SU}(B)$ be a group homomorphism which is a contraction with respect to $\widetilde{D}_A$ and $\widetilde{D}_B$. Suppose that the following statement hold: \\
(a) $\alpha\in KK(A,B)$ can be realized by a unital homomorphism;\\
(b) $\alpha$ is $L$-full;\\
(c) the map $\alpha_{0}: K_{0}(A) \rightarrow  K_{0}(B) $ induced by $\alpha$ is compatible with $\xi$, i.e., the diagram
$$\CD
  K_{0}(A) @>\rho_{A}>> AffTA \\
  @V \alpha_{0} VV @V \xi VV  \\
  K_{0}(B) @>\rho_{B}>> AffTB
\endCD $$
commutes. In other words, the map $\xi': AffTA / \widetilde{\rho K_{0}(A)} \rightarrow AffTB / \widetilde{\rho K_{0}(B)} $ is well defined. Recall that $\widetilde{\rho K_{0}(A)}$ is the closure of real vector spaces spanned by $\rho K_{0}(A)$;\\
(d) the map $\alpha_{1}: K_{1}(A)/ tor  K_{1}(A)\rightarrow  K_{1}(B)/ tor  K_{1}(B)$ induced by $\alpha$ is compatible with $\gamma$, i.e.,  the diagram
$$\CD
  U(A)/ \widetilde{SU(A)} @>\pi_{A}>> K_{1}(A) / tor  K_{1}(A)\\
  @V \gamma VV @V \alpha_{1} VV  \\
   U(B)/ \widetilde{SU(B)} @>\pi_{B}>> K_{1}(B) / tor  K_{1}(B)
\endCD$$
commutes;\\
(e) the maps $\xi$ and $\gamma$ are compatible,  i.e.,  the diagram
$$\CD
  AffTA/ \widetilde{\rho K_{0}(A)} @>\widetilde{\lambda}_{A}>> U(A)/ \widetilde{SU(A)} \\
  @V \xi' VV @V \gamma VV  \\
   AffTB/ \widetilde{\rho K_{0}(B)} @>\widetilde{\lambda}_{B}>> U(B)/ \widetilde{SU(B)}
\endCD$$
commutes, where $\xi'$ is as in $(c)$.

It follows that there exists  a unital homomorphism $\phi:A\rightarrow B $ such that\\
(1) $KK(\phi)=\alpha $,\\
(2) $\parallel AffT \phi(f)-\xi(f) \parallel < \varepsilon $, $\forall f\in H \subseteq AffTA$, and\\
(3) $dist(\phi^{\natural}(f),\gamma(f))\triangleq \widetilde{D}_{B}(\phi^{\natural}(f),\gamma(f))< \varepsilon $
for all $f\in E \subset S U(A)/ \widetilde{SU(A)}$, where $\phi^{\natural}$ denotes the map
$U(A)/ \widetilde{SU(A)} \rightarrow U(B)/ \widetilde{SU(B)}$ induced by $\phi$.

\noindent\textbf{Proof.} For the case $Sp(A)\neq S^{1}$,
$K_{1}(A)$ is either 0 or a torsion group.
 By Lemma 2.7, \\
$U(A)/ \widetilde{SU(A)}=tor U(A)/ \widetilde{SU(A)}$ is isomorphic to
$AffTA/ \widetilde{\rho K_{0}(A)}$ and therefore requirement (3) follows from (2). (In this case,  the diagram in (d) commutes  since $ K_{1}(A)/ tor K_{1}(A) $=0. Also the diagram in (e)
identifies $\xi'$ and $\gamma$.) If $A=M_{l}(C(S^{1}))$ and $B^{i}$ is homogeneous, then the proposition is just
[EGL, Theorem 3.17]. For the case that  $A=M_{l}(C(S^{1}))$ and $B^{i}$ is a dimension drop algebra, one can still adopt the proof of [EGL, Theorem 3.17],  with  [EGL, Remark 3.16] replaced by Remark 5.5 above. That is, one can change the map $\gamma$
in 5.5 suitably to make the estimation in (3) hold.\\
 $~~~~~~~~~~~~~~~~~~~~~~~~~~~~~~~~~~~~~~~~~~~~~~~~~~~~~~~~~~~~~~~~~~~~~~~~~~~~~~~~~~~~~~~~~~~~~~~~~~~~~~~~~~~~~~~~~~~~~~~~~~~~\Box$

The following lemma is a generalization of Lemma 2.5 of [Ji-Jiang]. We will give a little simpler proof. Recall again for a unital
homomorphism $\phi: A \rightarrow B $, we denote the induced map $AffTA \rightarrow AffTB $ by $AffT\phi$ and the induced map $U(A)/ \widetilde{SU(A)} \rightarrow U(B)/ \widetilde{SU(B)}$ by $\phi^{\natural}$.

\noindent\textbf{Lemma 5.7.} Let $A_{1}=PM_{l_{2}
}(C(X))P \in \mathcal{HD}$ or $A_{1}=M_{l_{1}}(I_{k})$,
and $H\subseteq AffTA_{1}$ and $E\subseteq U (A_{l})/ \widetilde{SU (A_{l})}$ be finite subsets. Let $\phi, \psi: A_{1}\rightarrow A_{2}$ ($\in \mathcal{HD}$) be two unital homomorphisms with $\psi$ being defined by point evaluation such that
$$KK(\phi)=KK(\psi),$$ $$\parallel AffT\phi(h)- AffT\psi(h)  \parallel < \varepsilon_{1}~~~~ \forall h\in H,~~~~\mbox{and} $$
$$dist(\phi^{\natural}(g),\psi^{\natural}(g))=dist((\phi^{\natural}(g),0)< \varepsilon_{2}~~~~\forall g\in F.$$
(Note that $\psi^{\natural}(g)=0$; see 2.8.)
%since $\psi(U (A_{1}))\subset \widetilde{SU (A_{2})}$ for any homomorphism $\psi$ defined by point evaluation.

Let $A_{3}\in \mathcal{HD}$, $\widetilde{\wedge}: A_{2}\rightarrow A_{3}$ be a unital homomorphism, $\xi: AffTA_{2}\rightarrow AffTA_{3} $ a unital positive linear map which is compatible with $K_{0}(\widetilde{\wedge})$ (see (c) of (5.6)), and let $\gamma:U(A_{2})/ \widetilde{SU(A_{2})} \rightarrow U(A_{3})/ \widetilde{SU(A_{3})}$ be any contractive (with respect to  $\widetilde{D}_{A_{2}}$ and $\widetilde{D}_{A_{3}}$) group homomorphism which is compatible with $$\xi': AffTA_{2}/ \widetilde{\rho K_{0}(A_{2})}\rightarrow AffTA_{3}/ \widetilde{\rho K_{0}(A_{3})}$$
(induced by $\xi$, see (e) of 5.6) and compatible with $$K_{1}(\widetilde{\wedge}):K_{1}(A_{2})/ tor K_{1}(A_{2}) \rightarrow K_{1}(A_{3})/ tor K_{1}(A_{3})$$ ( see ($d$) of 5.6 ). Let $\wedge: A_{1}\rightarrow A_{3}$ be defined by $\wedge=\widetilde{\wedge}\circ \psi $.
It follows that \\
(1) $KK(\wedge)=KK(\widetilde{\wedge}\circ \psi)$;\\
(2) $\parallel AffT\wedge(f)-(\xi\circ AffT\phi)(f) \parallel < \varepsilon_{1}$~~ $\forall f\in H$;\\
(3) $dist(\wedge^{\natural}(g),\gamma\circ\phi^{\natural}(g))< \varepsilon_{2}$~~ $\forall g\in E $.

\noindent\textbf{Proof.} Suppose that $B$ is a finite dimensional $C^{\ast}$~algebra and $C$ is a unital $C^{\ast}$~algebra.
If $\phi: B\rightarrow C$ is a unital homomorphism and $\xi: AffTB\rightarrow AffTC $ is an unital positive linear map compatible with $K_{0}(\phi)$, then $\xi=AffT\phi$. Namely, because $(K_{0}(B),K_{0}(B)_{+})=(\mathbb{Z}^{t},{\mathbb{Z}^{t}}_{+})$ and $AffTB={\mathbb{R}}^{t}$ for some $t\in \N$, there is only one map from $AffTB$ to $AffTC$ which is compatible with $K_{0}(\phi)$, and consequently it will be $AffT\phi$.

Let $\psi: A_{1}\rightarrow A_{2}$ be written as $\psi_{2}\circ \psi_{1}$, where $\psi_{1}:A_{1}\rightarrow B $ and $\psi_{2}:B\rightarrow A_{2}$ are
unital homomorphisms with $B$ being finite dimensional. Since $\xi:AffTA_{2}\rightarrow AffTA_{3}$ is compatible with $K_{0}(\widetilde{\wedge})$, we know that $\xi\circ AffTA(\psi_{2}):AffTB\rightarrow AffTA_{3}$ is compatible with $K_{0}(\widetilde{\wedge}\circ \psi_{2}):K_{0}(B)\rightarrow K_{0}(A_{3})$. Consequently, $\xi\circ AffT(\psi_{2})=AffT(\widetilde{\wedge}\circ \psi_{2})$. Hence $\xi\circ AffT\psi=AffT(\widetilde{\wedge}\circ \psi)=AffT\wedge$, which implies (2) above. The above (1) and (3) are obvious
( note that $\wedge^{\natural}=0$ since $\psi^{\natural}=0$ ).\\
$~~~~~~~~~~~~~~~~~~~~~~~~~~~~~~~~~~~~~~~~~~~~~~~~~~~~~~~~~~~~~~~~~~~~~~~~~~~~~~~~~~~~~~~~~~~~~~~~~~~~~~~~~~~~~~~~~~~~~~~~\Box$

\noindent\textbf{5.8.}~~Let $A=lim(A_{n},\phi_{n, m})$, $B=lim(B_{n},\psi_{n, m})$ be two unital $A\mathcal{HD}$ inductive limits. Suppose that there is an isomorphism
$$\alpha:(\underline{K}(A),\underline{K}(A)^{+},\textbf{1}_{A})\rightarrow (\underline{K}(B),\underline{K}(B)^{+},\textbf{1}_{B}))$$
which is induced by a shape equivalence as in 3.5 (we also regard $\alpha \in KL(A,B)$ see 2.3),
suppose that there is a positive linear isomorphism
$$\xi:AffTA\rightarrow AffTB$$
which is compatible with $\alpha_{0}$(induced by $\alpha$), and suppose that there is a isomorphism\\
$$\gamma:U(A) / \widetilde{SU(A)}\rightarrow U (B) / \widetilde{SU (B)}$$
which is compatible with both $\xi':AffTA/ \widetilde{\rho K_{0}(A)}\rightarrow AffTB/ \widetilde{\rho K_{0}(B)}$
(induced by $\xi$) and\\
 $\alpha_{1}:K_{1}(A)/ tor K_{1}(A)\rightarrow K_{1}(B)/ tor K_{1}(B)$
(induced by $\alpha$) (see (d) and (e) of 5.6).

\noindent\textbf{Proposition 5.9.} Let $A=lim(A_{n},\phi_{n, m}),B=lim(B_{n},\psi_{n, m})$, $\alpha$, $\xi$, $\gamma$ be as in 5.8. Then for any $A_{n}$, any finite sets $H\subseteq AffTA_{n}$, $E\subset U(A_{n})/ \widetilde{SU(A_{n})}$ and any $\varepsilon > 0$, there are an $m>n$, a $KK$-element $$\alpha_{n}\in KK(A_{n},B_{m})$$ which can be realized by a unital homomorphism, a unital positive linear map $$\xi_{n}:AffTA_{n}\rightarrow AffTB_{m}$$ and a contractive group
homomorphism $$\gamma_{n}:U (A_{n})/ \widetilde{SU (A_{n})}\rightarrow U (B_{m}) / \widetilde{SU (B_{m}) }$$
with the following properties\\
(1) $\alpha_{n}\times KL(\psi_{m \infty})=KL(\phi_{n \infty})\times \alpha \in KL(A_{n},B)$;\\
(2) $\parallel( AffT\psi_{m \infty}\circ \xi_{n})(f)-(\xi\circ AffT\phi_{n \infty})(f) \parallel < \varepsilon$ ~~$\forall f\in H$;\\
(3) $dist((\psi_{m \infty}^{\natural}\circ\gamma_{n})(f),(\gamma \circ \phi_{n \infty}^{\natural})(f))< \varepsilon $~~~ $\forall f\in E $;\\
(4) $\alpha_{n}^{0}$ and $\xi_{n}$ are compatible (see (c) of 5.6);\\
(5) $\gamma_{n}$ and $\xi_{n}$ are compatible (see (e) of 5.6);\\
(6) $\gamma_{n}$ and $\alpha_{n}^{1}$ are compatible (see (d) of 5.6).

\noindent\textbf{Proof.} Notice that Lemma 3.21 of [EGL] is true in our seting --- that is, with building blocks like $M_{l}(I_{k})$, it will not make any difference in the proof, since both $AffT$ space, $K_{0}$ group and pairing between them for $M_{l}(I_{k})$ are the same as the case $M_{l}([0,1])$ (the simplicity of the inductive limit is not used in the proof of [Li1, Lemma 6.6] and also [EGL, Lemma 3.21], as  pointed out in [Ji-Jiang Lemma 2.1 ]). Then the proof of this proposition is completely the same as [EGL, Lemma 3.23 ].\\
$~~~~~~~~~~~~~~~~~~~~~~~~~~~~~~~~~~~~~~~~~~~~~~~~~~~~~~~~~~~~~~~~~~~~~~~~~~~~~~~~~~~~~~~~~~~~~~~~~~~~~~~~~~~~~~~~~~~~~\Box$

The following is the main existence theorem of this article:

\noindent\textbf{Theorem 5.10.}  Let $A=Lim(A_{n},\phi_{n, m})$, $B=Lim(B_{n},\psi_{n, m})$
be (not necessary unital) $A\mathcal{HD}$ inductive limit algebras, with the ideal property. For any
 \begin{center}
 $H=\oplus H^{i}\subseteq \oplus_{i}AffTA_{n}^{i}$, $E=\oplus E^{i}\subseteq \oplus_{i} U(A_{n}^{i})/ \widetilde{SU(A_{n}^{i})}=U(A_{n}) / \widetilde{SU(A_{n})} $, $ \varepsilon_{1}> 0 $, $\varepsilon_{2}> 0$,
 \end{center}
 there is an $m>n$ with the following property.

Suppose that $\alpha \in KK(A_{m},B_{l})$ can be realized by a homomorphism, suppose that
$$\xi:AffTA_{m}\rightarrow AffT(\alpha [\textbf{1}_{A_{m}}]B_{l}\alpha [\textbf{1}_{A_{m}}] )$$
(here, $\alpha [\textbf{1}_{A_{m}}]$  denotes any projection representing $\alpha_{\ast} [\textbf{1}_{A_{m}}])$ is a unital positive linear map compatible with\\
~~~$\alpha_{0}:K_{0}(A_{m})\rightarrow K_{0}(\alpha [\textbf{1}_{A_{m}}])B_{l}\alpha [\textbf{1}_{A_{m}]})$ (see (c) of 5.6), and suppose that
~~~$$\gamma:U (A_{m})/ \widetilde{SU(A_{m})}\rightarrow U (\alpha [\textbf{1}_{A_{m}}]B_{l}\alpha [\textbf{1}_{A_{m}}])/ \widetilde{SU (\alpha [\textbf{1}_{A_{m}}]B_{l}\alpha [\textbf{1}_{A_{m}}])}$$ is a contractive group
homomorphism compatible with\\
$\bar{\xi}: AffTA_{m} / \widetilde{\rho K_{0}(A_{m})}\rightarrow AffT (\alpha [\textbf{1}_{A_{m}}])B_{l}\alpha [\textbf{1}_{A_{m}}])/ \widetilde{\rho_{0} K_{0}(\alpha [\textbf{1}_{A_{m}}])B_{l}\alpha [\textbf{1}_{A_{m}}])}$
(induced by $\xi$)\\ and comptible with $\alpha_{1}: K_{1}(A_{m})/ tor K_{1}(A_{m})\rightarrow K_{1}(B_{l})/ tor K_{1}(B_{l}).$

Then there is a homomorphism $\wedge: A_{n}\rightarrow B_{l}$ with the following properties:\\
(1) $KK(\wedge)=KK(\phi_{n,m})\times\alpha \in KK(A_{n},B_{l})$;\\
(2) Let $P=\phi_{n,m}(\textbf{1}_{A_{n}})\leq \textbf{1}_{A_{m}}$ and $\xi'=\xi\mid_{([P],\alpha[P])}: AffTPA_{m}P\rightarrow AffT\alpha[P]B_{l}\alpha[P]$ be the map  determined by $\xi$
(as in 2.11),  and let $$\widetilde{\xi'}:AffTA_{n}\rightarrow AffT\alpha[P]B_{l}\alpha[P]$$ be geven by $$\widetilde{\xi'}=\xi'\circ AffT\phi_{nm},$$ where $AffT\phi_{nm}:AffTA_{n}\rightarrow AffTPA_{m}P$ is induced by $\phi_{nm}$ (see 2.8). Finally let $$\xi''^{i}=\widetilde{\xi'}\mid_{([\textbf{1}_{A_{n}^{i}}],\wedge [\textbf{1}_{A_{n}^{i}}])}:AffTA_{n}^{i}\rightarrow AffT\wedge(\textbf{1}_{A_{n}^{i}}) B \wedge(\textbf{1}_{A_{n}^{i}}).$$
Then we have
$$
 \parallel (AffT\wedge)\mid_{(\textbf{1}_{A_{n}^{i}},\wedge [\textbf{1}_{A_{n}^{i}}])}(h)-\xi''^{i}(h)\parallel < \varepsilon_{1},\;\forall h\in H^{i};
$$
(3) Let
$
\gamma^{i}:U(A_{n}^i)/ \widetilde{SU(A_{n}^i)}\rightarrow U(\wedge(\textbf{1}_{A_{n}^{i}}) )B (\wedge(\textbf{1}_{A_{n}^{i}}))/ \widetilde{SU}(\wedge(\textbf{1}_{A_{n}^{i}}) B (\wedge(\textbf{1}_{A_{n}^{i}})))
$
be defined as \\
$(\gamma \mid_{([P],\alpha[P])}\circ \phi_{n,m}^{\natural})_{([\textbf{1}_{A_{n}^{i}}],\wedge[\textbf{1}_{A_{n}^{i}}])}$,
then
%\begin{center}
$$\parallel\wedge^{\natural}|_{(\textbf{1}_{A_{n}^{i}},\wedge(\textbf{1}_{A_{n}^i}))}(g)-\gamma^{i}(g) \parallel < \varepsilon_{2}~~~~~~\forall g\in E^{i}$$
%\end{center}
Here,  for a (not necessary unital) homomorphism,
$\phi: A\rightarrow B$, the maps $AffT\phi$ and $\phi^{\natural}$ are defined for the corresponding unital homomorphisms as in 2.4 and 2.8.
%will denote the unital maps:$AffTA\rightarrow AffT\phi(\textbf{1}_{A})B\phi(\textbf{1}_{A})$ and
%$U(A)\widetilde{SU(A)}\rightarrow U(\phi(\textbf{1}_{A})B\phi(\textbf{1}_{A}))/ \widetilde{SU} (\phi(\textbf{1}_{A})B\phi(\textbf{1}_{A}))$, induced by
%$\phi: A\rightarrow \phi(\textbf{1}_{A})B\phi(\textbf{1}_{A}).$

\noindent\textbf{Proof.} The theoerm follows from the dichotomy theorem Proposition 3.2,  Proposition 5.6 (for $L$-large $KK$ component) and Lemma 5.7 (for $KK$ component which is not $L$-large). Note that the estimation for $\parallel \phi_{n, m}^{ij}-\psi_{n, m}^{ij} \parallel$ in Proposition 3.2 (b) also give the estimation for $\parallel AffT\phi_{n, m}^{ij} - AffT\psi_{n, m}^{ij} \parallel$.\\
$~~~~~~~~~~~~~~~~~~~~~~~~~~~~~~~~~~~~~~~~~~~~~~~~~~~~~~~~~~~~~~~~~~~~~~~~~~~~~~~~~~~~~~~~~~~~~~~~~~~~~~~~~~~~~~~~~~~~~~~~~~~\Box$

\vspace{3mm}

\noindent\textbf{\S6. The uniqueness theorem }

%\vspace{-2.4mm}

%\vspace{3mm}

In this section, we will prove several uniqueness theorems --- one of the main ingredients in the classification theory.

\noindent\textbf{6.1.}~~Let $X$ be a 2-dimensional connected simplicial complex. It is well known that the isomorphism class of
vector bundle $E$ is completely determined by rank(E) and Chern class $c_{1}(E)\in H^{2}(X)$.
That is,  $K^{0}(X)\cong \mathbb{Z}\oplus H^{2}(X)$. The following facts are also well known. If $dim(E)<dim(F)$, then $E$ is isomorphic to a subbundle of
$F$. Furthermore,  if $E$ is a subbundle of $F$ (denoted by $E\subset F $) with $dim(E)<k<dim(F)$, then there is a trivial bundle $G$ such that $E\subset G\subset F$ and $dim(G)=k$.

\noindent\textbf{6.2.}~~Let $[a,b]\subset [0,1] (a<b)$ be a closed interval. Let $X$ be any connected finite simplicial complex and $\phi:C[a,b]\rightarrow PM_{\bullet}(C(X))P$ be a unital homomorphism. Then for any $x\in X$, one can order the set $Sp\phi_{x}$ as
$$a\leq t_{1}(x)\leq t_{2}(x)\leq \cdots \leq t_{n}(x)\leq b~~~~~~ (n=rank(p)).$$
In particular, all $t_{i}(x)$ are continuous functions from $X$ to $[a,b]$. If $dim(X)\leq 2 $, then
by [Choi-Ell], for any finite set $F\subset C[a,b]$ and $\varepsilon>0$, there is a $\psi\!:C[a,b]\rightarrow PM_{\bullet}(C(X))P$
such that $Sp\psi_{x}$ are distinct for each $x\in X$ and $\parallel\phi(f)-\psi(f)\parallel < \varepsilon$ for all $f\in F$.
In this case, $Sp\psi_{x}$ can be written as $$a\leq s_{1}(x)<s_{2}(x)<\cdots < s_{n}(x)\leq b.$$
One can define $p_{i}(x)$ to be the spectrum projection corresponding to the spectrum $s_{i}(x)$ (see [G5  1.2.8 and 1.2.9])
Then $p_{i}\in PM_{\bullet}(C(X))P$ and $\Sigma p_{i} = P$. Furthermore,  $\psi(f)(x)=\sum_{i=1}^{n}f(s_{i}(x))p_{i}(x)$.

\noindent\textbf{Lemma 6.3.} Let $[a,b]\subset[0,1]$, and  $F\subset C[a,b]$ be $(\frac{\varepsilon}{3},2\eta)$ uniformly continuous (see 1.25).
Let $\phi: C[a,b]\rightarrow PM_{\bullet}(C(X))P$ ($dim X=2$) be a homomorphism such that for any
$x\in X$ and any interval $(c,d)\subset[a,b]$ with length at least $\eta$,  $Sp\phi_{x}\cap(c,d)\neq\emptyset$.
Then there are mutually othogonal projections $p_{1},p_{2},\cdots,p_{n}$ ($n=rank(P)$) and continuous
functions $t_{1},t_{2},\cdots,t_{n}:X\rightarrow[a,b]$ such that $p_{1},p_{2},\cdots,p_{n-1}$ are trivial projections
and $\sum_{i=1}^{n}p_{i}=P$ ($p_{n}$ trivial if and only if $P$ trivial)
and $\parallel \phi(f)-\xi(f)\parallel < \varepsilon $, where $\xi(f)(x)=\sum_{i=1}^{n}f(t_{i}(x))p_{i}(x)$ for all $f\in C[a,b]$.

\noindent\textbf{Proof.} As mentioned in 6.2, one can find a homomorphism $\psi$ such that $\parallel\psi(f)-\phi(f)\parallel< \frac{\varepsilon}{6}$,
such that $Sp\psi_{x}\cap(c,d)\neq\emptyset$ if $length(c,d)\geq \eta$ and,  such that $Sp\psi_{x}$ are distinct.

Write $\psi(f)(x)=\sum f (s_{i}(x))q_{i}(x)$ as in 6.2, where $s_{i}\!:X\rightarrow[a,b]$ are continuous and $q_{i}\in PM_{\bullet}(C(X))P$ are rank 1 projections. Then $s_{1}(x)<s_{2}(x)<\cdots <s_{n}(x)$ and $\mid s_{i+1}(x)-s_{i}(x)\mid<\eta $ for each $i$.

By 6.1  there are trivial projections $P_{1},P_{3},\cdots,P_{2k-1}$ (with $2k-1$ being largest odd integer with $2k-1\leq n-1$), such that \\
$P_{-1}=0<P_{1}< q_{1}+q_{2}<P_{3}< q_{1}+q_{2}+q_{3}+q_{4}<P_{5}< q_{1}+q_{2}+q_{3}+q_{4}+q_{5}+q_{6}<\cdots <q_{1}+q_{2}+\cdots+q_{2k-2}  < P_{2k-1}<q_{1}+q_{2}+\cdots+q_{2k}\leq P$.

Let $\xi:C[a,b]\rightarrow PM_{\bullet}(C(X))P$ be defined by\\
$\xi(f)(x)=\sum_{i=0}^{k-1}f(s_{2i+1}(x))(P_{2i+1}(x)-P_{2i-1}(x))+f(s_{n}(x))(P(x)-P_{2k-1}(x))$.

Let~~ $\widetilde{p}_{1}=P_{1}$,
~~~$\widetilde{p}_{2}=q_{1}+q_{2}-P_{1}$,
~~~$\widetilde{p}_{3}=P_{3}-(q_{1}+q_{2})$,
~~~$\widetilde{p}_{4}=(q_{1}+q_{2}+q_{3}+q_{4})-P_{3}$,~~ $\cdots\cdots$,\\
$\widetilde{p}_{2k-1}=P_{2k-1}-(q_{1}+q_{2}+\cdots+q_{2k-2})$,~~~ and ~~~
$\widetilde{p}_{2k}=(q_{1}+q_{2}+\cdots+q_{2k})-P_{2k-1}$~~~ and\\
$\widetilde{p}_{n}=P-(q_{1}+q_{2}+\cdots+q_{2k})$ if $n=2k+1$.

Then according to the decomposition $\widetilde{p}_{1},\widetilde{p}_{2},\cdots,\widetilde{p}_{n}$,
we have
$$
\xi(f)=f(s_{1}(x))\widetilde{p}_{1}(x)+f(s_{3}(x))(\widetilde{p}_{2}(x)+\widetilde{p}_{3}(x))+\cdots
+f(s_{2k-1}(x))(\widetilde{p}_{2k-2}(x)+\widetilde{p}_{2k-1}(x))+f(s_{k}(x))(\widetilde{p}_{2k}(x)+\ast),
$$
where $\ast=\widetilde{p}_{n}(x)$ or 0 depending on $n=2k+1$, or $n=2k$.

On the other hand, $\widetilde{p}_{1}+\widetilde{p}_{2}=q_{1}+q_{2}$,
$\widetilde{p}_{3}+\widetilde{p}_{4}=q_{3}+q_{4}$, $\cdots$, which are commutative with $\psi(f)$ and
$(\widetilde{p}_{2i-1}(x)+\widetilde{p}_{2i}(x))\psi(f)(x)(\widetilde{p}_{2i-1}(x)+\widetilde{p}_{2i}(x))=f(s_{2i-1}(x))q_{2i-1}(x)+f(s_{2i}(x))q_{2i}(x)$
which is closed to $f(s_{2i}(x))(\widetilde{p}_{2i-1}(x)+\widetilde{p}_{2i}(x))$ to within $\frac{\varepsilon}{3}$ (note that $|s_{2i}(x)-s_{2i-1}(x)|\leq \eta )$. Also, we have
\begin{center}
$(\widetilde{p}_{2i-1}(x)+\widetilde{p}_{2i}(x))\xi(f)(\widetilde{p}_{2i-1}(x)+\widetilde{p}_{2i}(x))=f(s_{2i-1}(x))\widetilde{p}_{2i-1}(x)+f(s_{2i}(x))\widetilde{p}_{2i}(x)$
\end{center}
which is also close to $f(s_{2i}(x))(\widetilde{p}_{2i-1}(x)+\widetilde{p}_{2i}(x))$ to within $\frac{\varepsilon}{3}$.
Hence $\parallel\psi(f)(x)-\xi(f)(x)\parallel<\frac{2\varepsilon}{3}$.

Evidently, $\xi$ is as desired.

%Finally note that $\xi$ is as desired, since each trivial projection of rank 2 can be written as sum of two trivial projection %of rank 1, and if $n=2k+1$, then the last projection $P-p_{2k-1}$ can be written as sum of two projection of rank 1, with %the first of them being trivial.\\
$~~~~~~~~~~~~~~~~~~~~~~~~~~~~~~~~~~~~~~~~~~~~~~~~~~~~~~~~~~~~~~~~~~~~~~~~~~~~~~~~~~~~~~~~~~~~~~~~~~~~~~~~~~~~~~~~~~~~~~~~~~\Box$

\noindent\textbf{Remark 6.4.} One can see the following from the proof of Lemma 6.3. The homomorphism $\psi:C[a,b]\rightarrow PM_{\bullet}(C(X))P$ can be chosen such that $Sp\psi_{x}\cap(c,d)\neq\emptyset$ for any interval $(c,d)\subset [a,b]$ of length $\eta$, and $Sp\psi_x$ can be paired
with $Sp\phi_{x}$ (for any $x\in X $) to within arbitrarily given small number. Also $\xi: C[a,b]\rightarrow PM_{\bullet}(C(X))P $ can be chosen so
that $Sp\xi_{x}\cap(c,d)\neq\emptyset$ for any interval $(c,d)\subset [a,b]$ of length $2\eta$, and $Sp\xi_{x}$ can be paired with $Sp\psi_{x}$
to within $\eta$. Consequently, we can choose $\xi$ such that $Sp\xi_{x}$ and $Sp\phi_{x}$ can be paired to within $\eta$.

\noindent\textbf{Lemma 6.5.} Let $F\subset C[a,b]$ be $(\frac{\varepsilon}{9},2\eta)$ uniformly continuous. Let $X$ be a connected simplicial complex of dimension at most 2. Let $\phi_{1},\phi_{2}:C[a,b]\rightarrow PM_{\bullet}(C(X))P$ be unital homomorphisms
such that for any $x\in X$, $Sp(\phi_{1})_{x}$ and $Sp(\phi_{2})_{x}$ can be paired to within $\frac{\eta}{5}$ and $Sp(\phi_{1})_{x}\cap(c,d)\neq\emptyset$ for any interval $(c,d)\subset [a,b]$ of length $\frac{\eta}{5}$. Then there is a unitary $u\in  PM_{\bullet}(C(X))P $ such that $$\parallel \phi_{1}(f)-u\phi_{2}(f)u^{\ast}\parallel<\varepsilon~~~~ \forall f\in F.$$

\noindent\textbf{Proof.} Note that $Sp(\phi_{2})_{x}\cap(c,d)\neq\emptyset$ for any $(c,d)$ of length $\frac{3\eta}{5}$. By Lemma 6.3, there
are $\xi_{1},\xi_{2}:C[0,1]\rightarrow PM_{\bullet}(C(X))P$ of forms $\xi_{1}(f)=\sum_{i=1}^{n}f(s_{i}(x))p_{i}$, $\xi_{2}(f)=\sum_{i=1}^{n}f(t_{i}(x))q_{i}$ with
$\{ p_{i}\}_{i=1}^{n-1}$ and $\{ q_{i}\}_{i=1}^{n-1}$ being trivial projections, $s_{i},t_{i}:~X\rightarrow [a,b]$ being continuous functions
with $a\leq s_{1}(x)\leq  s_{2}(x) \cdots \leq  s_{n}(x)\leq b$, $a\leq t_{1}(x)\leq t_{2}(x) \cdots \leq t_{n}(x)\leq b,$
such that  $\parallel\phi_{1}(f)-\xi_{1}(f)\parallel < \frac{\varepsilon}{3}$, $\parallel\phi_{2}(f)-\xi_{2}(f)\parallel < \frac{\varepsilon}{3}$.

Furthermore,  as in Remark 6.4, $Sp(\phi_{1})_{x}$ and $Sp(\xi_{1})_{x}$ can be paired with $\frac{\eta}{5}$, and $Sp(\phi_{2})_{x}$ and $Sp(\xi_{2})_{x}$ can be paired with $\frac{3\eta}{5}$. Consequently $Sp(\xi_{1})_{x}$ and $Sp(\xi_{2})_{x}$ can be paired to within $\frac{\eta}{5}+\frac{3\eta}{5}+\frac{\eta}{5}=\eta$ for each $x\in X $. That is,  $|t_{i}(x)-s_{i}(x)|<\eta$.
Since $\{ p_{i}\}_{i=1}^{n-1}$ and $\{ q_{i}\}_{i=1}^{n-1}$ are trivial, we have that $p_{n}=P-\sum_{i=1}^{n-1}p_{i}$ is unitarily equivalent
to $q_{n}=P-\sum_{i=1}^{n-1}q_{i}$. There is a unitary $u\in PM_{\bullet}(C(X))P$ such that $uq_{i}u^{*}=p_{i}$ for~~$i=1,2,\cdots n$.
Consequently,  $\|\xi_{1}(f)-u\xi_{2}(f)u^{\ast}\|<\frac{\varepsilon}{3}$~for all $f\in F$.
Hence
%\begin{center}
$$\|\phi_{1}(f)-u\phi_{2}(f)u^{\ast}\|<\frac{\varepsilon}{3}+\frac{\varepsilon}{3}+\frac{\varepsilon}{3}=\varepsilon~~~~ \forall f\in F.$$
%\end{center}
$~~~~~~~~~~~~~~~~~~~~~~~~~~~~~~~~~~~~~~~~~~~~~~~~~~~~~~~~~~~~~~~~~~~~~~~~~~~~~~~~~~~~~~~~~~~~~~~~~~~~~~~~~~~~~~~~~~~~~~~~~~\Box$

\noindent\textbf{Remark 6.6.} In Lemma 6.5 (Lemma 6.3 respectively) if length$[a,b]<\frac{\eta}{5}$ (length$[a,b]<\eta$, respectively),
then the requirement $Sp\phi_{x}\cap(c,d)\not=\emptyset$ for any $(c,d)\subset[a,b]$ of length $\frac{\eta}{5}$ (of length $ \eta$, respectively) is empty requirement. In this case, the theorem is evidently true.

The following lemma is an easy consequence of proof of Proposition 4.9 (also see [Ell1, 4.4] and [EGJS]).

\noindent\textbf{ Lemma 6.7.} Let $\phi_{1},\phi_{2}:C[a,b]\rightarrow M_{l}(I_{k})$ be two unital  homomorphisms such that $Sp(\phi_{1})_{t}$ and $Sp(\phi_{2})_{t}$ can be paired within $\eta$. Suppose that $F\subset C[a,b]$ is $(\varepsilon,2\eta)$ uniformly continuous. Then there is a unitary $u\in M_{l}(I_{k}) $ such that
$\|\phi_{1}(f)-u\phi_{2}(f)u^{\ast}\|<\varepsilon$~ for all $f\in F$.

\noindent\textbf{Lemma 6.8.} Let $F\subset I_{k}$ be a finite set and $\varepsilon>0$. There is an $\eta>0$, satisfying the following condition.

 Let $a>0$ and $A=I_{k}|_{[a,1]}$. If $\phi:A\rightarrow PM_{\bullet}(C(X))P$ is a homomorphism which satisfies that for any $x\in X$, $\sharp(Sp\phi_{x}\cap[1-\frac{\eta}{4},1])\geq k$ and $Sp\phi_{x}\cap (c,d)\neq\emptyset$ if $(c,d)\subset[a,1]$ is of length at least  $\eta$, then there is a homomorphism
$\psi:~A\rightarrow PM_{\bullet}(C(X))P$ such that $$\|\phi(f|_{[a,1]})-\psi(f|_{[a,1]})\|<\varepsilon ~~~~~\forall f\in F,$$ and $\psi$ is of the following form: There are mutually orthogonal projections $P_{1},P_{2}$ with $P_{1}+P_{2}=P$ and with rank $(P_{1})=kl$, $k+1\leq rank(P_{2})\leq2k$, and there are mutually orthogonal and mutuality unitarily equivalent projections $p_{1}\!\!=\!\!p,p_{2},\cdots,p_{k}$ with $P_{1}=\sum_{i=1}^{k}p_{i}$ and an identification $P_{1}M_{\bullet}(C(X))P_{1}\cong pM_{\bullet}(C(X))p\otimes M_{k}$, and a homomorphism $\psi_{1}:C[a,1]\rightarrow pM_{\bullet}(C(X))p$ of form
$\psi_{1}(f)(x)=\sum_{i=1}^{l}f(t_{i}(x))q_{i}$, where $q_{1},q_{2},\cdots,q_{l}$ are trivial rank 1 projections
with $\sum_{i=1}^{l}q_{i}=p$ and $a\leq t_{1}(x)\leq t_{2}(x)\leq\cdots\leq t_{l}(x)$,with
each $t_{i}(x)$ being a continuous function such that
%\begin{center}
 $$\psi(f)=(\psi_{1}\otimes id_{k}\circ \imath)(f)+{f}(\underline1)P_{2}~~~\forall f\in A,$$
%\end{center}
where $\imath:I_{k}|_{[a,1]}\rightarrow C[a,1]\otimes M_{k}$ is the canonical inclusion.

\noindent\textbf{Proof.}  This is a direct consequence of
Theorem 4.16 and Lemma 6.3. Namely, by Theorem 4.16, $\phi$ is approximated by $\phi'=\phi'_{1}\oplus\phi'_{2}$, where $\phi'_{1}$ factors through as
%\begin{center}
$$A|_{[a,1]}\rightarrow M_{k}(C[a,1]) \xrightarrow{\phi''_{1}} P_{1}^{\prime}M_{\bullet}C(X)P'_{1}$$
%\end{center}
 and $\phi'_{2}$ is defined
by $\phi'_{2}(f)=f(\underline1)P'_{2}$ with $1\leq rank(P'_{2})\leq k$. Hence $rank(P'_{1})=k(l+1)$---some multiple of $k$. Write $\phi''_{1}=\xi\otimes id_{k}$ and then apply Lemma 6.3 to $\xi$. That is, up to a small perturbation, $\xi$ can be written as the form of $\psi_{1}$ above with one more item $f(t_{l+1})(x)q_{l+1}$, and this last projection $q_{l+1}$ may not be trivial. Finally let $P_{2}=q_{l+1}\otimes \textbf{1}_{k}+P'_{2}$ and $P_{1}=P'_{1}-q_{l+1}\otimes \textbf{1}_{k}$ to get our conclusion.\\
$~~~~~~~~~~~~~~~~~~~~~~~~~~~~~~~~~~~~~~~~~~~~~~~~~~~~~~~~~~~~~~~~~~~~~~~~~~~~~~~~~~~~~~~~~~~~~~~~~~~~~~~~~~~~~~~~~~~~~\Box$

\noindent\textbf{Remark 6.9.} (a) The condition $\sharp(Sp\phi_{x}\cap[1-\frac{\eta}{4},1])\geq k$ follows from $Sp\phi_{x}\cap (c,d)\neq\emptyset$ for any interval of $(c,d)\subset[a,1]$ of length $\eta'=\frac{\eta}{4k}$. We can change $\eta$ to a smaller $\eta'$ and to remove the requirement $\sharp(Sp\phi_{x}\cap[1-\frac{\eta}{4},1])\geq k$.\\
(b) By symmetry, the lemma is also true for $A|_{[0,b]}$ (instead of $A|_{[a,1]}$) for some $b<1$.

\noindent\textbf{Lemma 6.10.} Let $F\subset I_{k}$ be a finite set and $\varepsilon>0$. There is an $\eta>0$ satisfies the following condition.

Let $A=I_{k}|_{[a,1]}$ for certain $a>0$ (or $A=I_{k}|_{[0,b]}$ for certain $b<1$). If $\phi,\psi:A\rightarrow PM_{\bullet}(C(X))P$ are two homomorphisms such that for any $x\in X$,  $Sp\phi_{x}\cap(c,d)\neq\emptyset$ for any
interval $(c,d)\subset[a,1]$ (or $(c,d)\subset [0,b]$) of length $\eta$, and $Sp\phi_{x}$ and $Sp\psi_{x}$ can be paired to within $\eta$ (for any $x\in X$),
then there is a unitary $u\in PM_{\bullet}(C(X))P$ such that
$\|\phi(f)-u\psi(f)u^{\ast}\|<\varepsilon$ ~for all $f \in \pi(F)$,  where $\pi: I_{k}\rightarrow A=I_{k}|_{[a,1]}$ (or $I_{k}|_{[0,b]}$) is the canonical quotient map.

\noindent\textbf{Proof.} This is a consequence of Lemma 6.8 (also see the proof of Lemma 6.5)\\
$~~~~~~~~~~~~~~~~~~~~~~~~~~~~~~~~~~~~~~~~~~~~~~~~~~~~~~~~~~~~~~~~~~~~~~~~~~~~~~~~~~~~~~~~~~~~~~~~~~~~~~~~~~~~~~~~~~~~~~~~~~~~\Box$

For the special case that both $A$ and $B$ are interval algebras, the following theorem is  [Ji-Jiang,  3.5].

\noindent\textbf{Theorem 6.11.}  Let $A$ be one of $M_{l}(C[0,1])$, $M_{l}(C(S^{1}))$, or $M_{l}(I_{k})$,  $F\subset A $ be any finite set, and  $\varepsilon>0$.
There is an $\eta>0$, such that for any $\delta>0$, there is a finite set $H\subset AffTA$  with following property. Let $X\subset Sp(A)$ be a  connected closed set
such that $X\neq Sp(A)$, ~ provided  $A=M_{l}(C(S^{1}))$ or $M_{l}(I_{k})$.  If  $B$ is any basic building block in $\mathcal{HD}$, and if $\phi,\psi:A|_{X}\rightarrow B$ are two unital homomorphisms such that\\
(1)$\phi$ has property $sdp(\frac{\eta}{32},\delta)$,\\
(2)$\|AffT\phi(h|_X)-AffT\psi(h_X)\|<\frac{\delta}{4}$ for all $h\in H$, and \\
(3)$K_{0}(\phi)=K_{0}(\psi)$,\\ then there is a unitary $u\in B$ such that $\|\phi(f)-u\psi(f)u^{\ast}\|<\varepsilon$ for all $f\in \pi(F)\subset A|_{X}$,
where $\pi:A\rightarrow A|_{X}$ is the canonical quotient map.

\noindent\textbf{ Proof.} Note that using the above condition (3), one can reduce the case to $A=C[0,1]$, $C(S^{1})$ or $I_{k}$,  that is, reduce to the case $l=1$ with $\phi(1)=\psi(1)$.
Then by   [Li2, Theorem 2.15] and [Ji-Jiang, Remark 3.4], conditions (1) and (2) (for suitable choice of $H$)  imply that, for any $x\in Sp(B)$,~ $Sp\phi_x$ and $Sp\psi_x$ can be paired to within $\frac{\eta}{4}$.
Then the theorem follows from Lemmas 6.5, 6.7 and 6.8 except the case that $A|_X=I_{k}|_{[a,1]}$ (or $A|_X=I_k|_{[0,b]}$) and $B$ is a dimension drop interval algebra.
For this case, we apply Proposition 4.14 and Remark 4.15.\\
$~~~~~~~~~~~~~~~~~~~~~~~~~~~~~~~~~~~~~~~~~~~~~~~~~~~~~~~~~~~~~~~~~~~~~~~~~~~~~~~~~~~~~~~~~~~~~~~~~~~~~~~~~~~~~~~~~~~~~~~~~~~\Box$

The following theorem is a special case of [Thm2, Lemma 3.1].
%(note that the dimension drop interval algebra is building block of type 4 in the paper [Thm2] which is simplier than the building block of type 2)

\noindent\textbf{Proposition 6.12.} ~For each pair $k,l\in \N$ such that $l>12$, there is a finite set $H\subseteq AffT(C(S^{1}))=C_{\mathbb{R}}(S^{1})$ with the following property.

Suppose that  $\phi,\psi: C(S^{1})\rightarrow M_{K}(I_{k_{1}})$ are two unital homomorphisms satisfying  the following conditions:\\
(1)  $\phi$ has $sdp ~(\frac{1}{4k},\frac{1}{l})$ and $sdp ~(\frac{1}{12l},2\delta)$;\\
(2)  $\|AffT\phi(h)-AffT\psi(h)\|<\delta$~~~for all $h\in H$;\\
(3)  For the canonical element $u\in C(S^{1})$ given by $u(z)=z$, there exist a continuous
function $\alpha:[0,1]\rightarrow \mathbb{R}$ with $|\alpha|<\frac{Kk_{1}}{l}$ and a constant $\lambda\in \mathbb{T}$
such that $Det(\psi(u)^{\mathcal{\ast}}\phi(u))(t)=\lambda e^{2\pi i\alpha(t)}$, for all $t\in [0,1]$ (this condition
is a consequence of the estimation $\widetilde{D}_{B}(\psi(u),\phi(u))<\frac{1}{l}$);\\
(4)  $K_{1}(\phi)=K_{1}(\psi)$.

It follows that for any $\varepsilon > 0 $, there is a unitary $w\in M_{K}(I_{k_{1}})$ such that $$\|\phi(u)-w\psi(u)w^{\ast}\|\leq (\frac{28}{k}+\frac{6}{l})\pi+\varepsilon.$$

Combining above proposition and [EGL, 2.13], we get

\noindent\textbf{Corollary 6.13}. Let $F\subset A=M_{\bullet}(C(S^{1}))$ be a finite set and $\varepsilon>0$. there is a number $\eta>0$ with property described below.

For any $\delta> 0$, there are $\delta'>0$ and $\widetilde{\eta}>0$ such that for any $\widetilde{\delta}>0$, there are
finite set $H\subseteq AffTA=C_{\mathbb{R}}(S^{1})$ and positive integer $L$ satisfying the following condition.

Suppose that two unital homomorphisms $\phi,\psi:A\rightarrow B\in \mathcal{HD}$ satisfy the following conditions:\\
(1)  $\phi$ has property $sdp(\frac{\eta}{32},\delta)$ and $sdp(\frac{\widetilde{\eta}}{32},\widetilde{\delta})$;\\
(2)  $\|AffT\phi(h)-AffT\psi(h)\|<\frac{\widetilde{\delta}}{4}$~~for all $h \in H$;\\
(3)  $\phi$ is $L$-large;\\
(4)  $\widetilde{D}_{B}(\phi(z),\psi(z))\leq\delta'$, where $z\!\in\! M_{k}(C(S^{1}))$ is defined by $z(e^{2\pi i t})\!=\!diag (e^{2\pi i t},1,\cdots,1)$; and \\
(5)  $K_{\ast}(\phi)=K_{\ast}(\psi)$.

%where $z\in M_{k}(C(S^{1}))$ is defined by
%$$z(e^{2\pi i t})=\left(
%                   \begin{array}{ccccc}
%                    e^{2\pi i t} &   &   &  &   \\
%                       & 1 &   &   &   \\
%                       &   & 1 &  &  \\
%                      &  &  & \cdots &  \\
%                      &  &  &  & 1 \\
%                   \end{array}
%                 \right)
%.$$
It follows that there exists a
unitary $u\in B$ such that $\|\phi(f)-u\psi(f)u^{\ast}\|<\varepsilon$~for all $f\in F$.

\begin{proof}
Note that $\widetilde{D}_{B}(\phi(z),\psi(z))<\delta'$ implies $D(\phi(z)\psi^{*}(z))<\delta'+\frac{2\pi }{L}$, where $D$ is defined in [EGL, 2.4]. (This is true because for any unitary $u\in B$ with constant determinant, we have $D(u)\leq\frac{2\pi }{L}$, for rank$(\textbf{1}_{B})\geq L$). For the block $B^i=QM_k(C(Y))Q$, this is [EGL, Proposition 2.13].  For the block $B^{i}=M_{k}(I_{k_{1}})$, the condition
$D(\phi(z))\psi^{*}(z)<\delta^{\prime}+\frac{2\pi }{L}$ implies the condition (3) in Proposition 6.12, provided $\delta^{\prime}+\frac{2\pi }{L}<\frac{1}{l}$; and therefore the proposition can be applied here. \\
\end{proof}

\noindent\textbf{6.14.}~~Fix $A=PM_{\bullet}(C(X))P$ for $X=T_{\uppercase\expandafter{\romannumeral2}
,k}$ or $A=M_{l}(I_{k})$. As in 5.16 of [G5],  $\alpha \in KK(A,B)$ is completely determined by
$\alpha_{0}^{0}:K_{0}(A)\rightarrow K_{0}(B)$,
$\alpha_{k}^{0}:K_{0}(A,\mathbb{Z}/ k)\rightarrow K_{0}(B,\mathbb{Z}/ k)$ and
$\alpha_{k}^{1}:K_{1}(A,\mathbb{Z}/ k)\rightarrow K_{1}(B,\mathbb{Z}/ k)$.

Note that, for any $C^{\ast}$-algebra $A$
$$K_{0}(A\otimes C(W_k\times S^{1}))=K_{0}(A)\oplus K_{1}(A)\oplus K_{0}(A,\mathbb{Z}/ k)
\oplus K_{1}(A,\mathbb{Z}/ k),$$
where $W_k=T_{\uppercase\expandafter{\romannumeral2}
,k}$.

Each projection $p\in M_{\infty}(A\otimes C(W_k\times S^{1}))$ defines an element $[p]\in K_{0}(A)\oplus K_{1}(A)\oplus K_{0}(A,\mathbb{Z}/ k)\oplus K_{1}(A,\mathbb{Z}/ k)$. For any finite set $\mathcal{P}\subset \cup_{k=1}^{\infty}M_{\infty}(A\otimes C(W_k\times S^{1}))$ of projections, denoted by $\mathcal{P}\underline{K}(A)$, the finite subset of $\underline{K}(A)$ consisting of elements coming from projection $p\in \mathcal{P} $.
That is,  $\mathcal{P}\underline{K}(A)=\{[p]\in \underline{K}(A), p\in \mathcal{P}\}$.

In particular,  if $A=PM_{\bullet}(C(T_{\uppercase\expandafter{\romannumeral2}
,k}))P$ or $M_{l}(I_{k})$, then we can choose a finite set of projections
$\mathcal{P}\subset M_{\infty}(A\otimes C(W_k\times S^{1})),$
such that the set
$$
 \{[p]\in K_{0}(A)\oplus K_{1}(A)\oplus K_{0}(A,\mathbb{Z}/ k)\oplus K_{1}(A,\mathbb{Z}/ k),\;p\in \mathcal{P}\}\subset \mathcal{P}\underline{K}(A)
$$
generates~~
$
K_{0}(A)\oplus K_{1}(A)\oplus K_{0}(A,\mathbb{Z}/ k)\oplus K_{1}(A,\mathbb{Z}/ k)\subseteq \underline{K}(A).
$

As in 5.17 of [G5], there are a finite set $G(\cal{P})\subset A$ and positive number $\delta({\cal{P}})>0$ such that if $\phi: A\rightarrow B$
is $G(\cal{P})$$-$$\delta(\cal{P})$
multiplicative (see Definition 1.2), then the map $\phi_{\ast}: \mathcal{P}\underline{K}(A)\rightarrow\underline{K}(B)$ is well defined.

The following results are quoted from [CJL].

\noindent\textbf{ Proposition 6.15.}~(Proposition 3.1 of [CJL])~ Let $A=PM_{\bullet}C(T_{\uppercase\expandafter{\romannumeral2}
,k})P$ or $M_{l}(I_{k})$ and $\mathcal{P}$ be as in 6.14. For any finite set $F\subset A$, $\varepsilon> 0$, there exist a  finite set $G\subset A$ ($G\supset G(\mathcal{P})$ large enough), a positive number $\delta > 0$ ($\delta<\delta(\mathcal{P})$ small enough)
such that the following statement is true.

If $B \in \mathcal{HD}$, $\phi,\psi\in Map(A,B)$ are $G$$-$$\delta$ multiplicative
and $\phi_{\ast}=\psi_{\ast}:\mathcal{P}\underline{K}(A)\rightarrow\underline{K}(B)$, then there is a homomorphism $\nu\in Hom (A,M_{L}(B))$ defined by point evaluations and there is a unitary $u\in M_{L+1}(B)$ such that
$\|(\phi\oplus \nu)(a)-u(\psi\oplus \nu)(a)u^{\ast}\|<\varepsilon$~for all $a\in F$.

\iffalse

\begin{proof}
For the special case that $A=PM_{\bullet}(C(T_{\uppercase\expandafter{\romannumeral2}
,k}))P$ and $B$ is homogeneous, this is [G5, Theorem 5.18].
 The proof of the general  case is completely same (see 6.14 above).
Note that, the calculation of
$K_{0}(\prod_{n=1}^{\infty}B_{n}/ \bigoplus_{n=1}^{\infty}B_{n})$,
$K_{0}(\prod_{n=1}^{\infty}B_{n}/ \bigoplus_{n=1}^{\infty}B_{n},\mathbb{Z}/ k)$
and
$K_{1}(\prod_{n=1}^{\infty}B_{n}/ \bigoplus_{n=1}^{\infty}B_{n},\mathbb{Z}/ k)$ in 5.12 of [G5] works well for the case that some of $B_{n}$ being of the form $M_{l}(I_{k})$. \\
\end{proof}

\fi

\noindent\textbf{Lemma 6.16.}~(Lemma 3.2 in [CJL])~Let $A=PM_{\bullet}(C(X))P$ or $M_{l}(I_{k})$. Let $F\subset A$ be approximately constant to within $\varepsilon$ (i.e.,  $\omega(F)<\varepsilon$). Then for any two homomorphisms $\phi,\psi:A\rightarrow B$ defined by point evaluations with $K_{0}\phi = K_{0}\psi $,   there exists a unitary $u\in B$ such that
$\|\phi(f)-u\psi(f)u^{\ast}\|< 2\varepsilon~~\forall~f\in F.$

\noindent\textbf{Lemma 6.17.}~(Lemma 3.3 in [CJL])~ Let $A=PM_{\bullet}(C(X))P$ or $M_{l}(I_{k})$, $\varepsilon> 0$, finite set $F\subset A $ with $\omega(A) < \varepsilon $. There exist a  finite set $G\subset A$ ($G\supset G(\cal{P})$), a number  $\delta>0$ ($\delta<\delta(\cal{P})$) and a positive integer $L$ such that the following statement is true. If $B\in \mathcal{HD}$, and $\phi,\psi \in Map (A,B)_{1}$ are $G\!-\!\delta$ multiplicative with $\phi_{\ast}=\psi_{\ast}: \mathcal{P}\underline{K}(A)\rightarrow \underline{K}(B)$
and $\nu: A\rightarrow M_{\infty}(B)$ is a homomorphism defined by point evaluations with $\nu([\textbf{1}_{A}])\geq L\cdot[\textbf{1}_{B}]\in K_{0}(B)$, then there is a unitary $u\in (\textbf{1}_{B}\oplus \nu(\textbf{1}_{A}))M_{\infty}(B)(\textbf{1}_{B}\oplus \nu(\textbf{1}_{A}))$ such that $$\|(\phi\oplus \nu)(f)-u(\psi\oplus \nu)(f)u^{\ast}\|<5\varepsilon~~~ \forall f\in F.$$

\iffalse

\begin{proof}
Let $L_{1}$ be as in Proposition 6.15 and $L=2L_{1}$. Since $[\nu(\textbf{1}_{A})]\geq L\cdot [\textbf{1}_{B}]=2L_{1}[\textbf{1}_{B}]\in K_{0}(B)$, there is a
projection $Q< \nu(\textbf{1}_{A})$ such that $Q$ is equivalent to $\textbf{1}_{M_{L_{1}}(B)}$. By Proposition 6.15, there exist  a homomorphism $\nu_{1}:A\rightarrow QM_{\infty}(B)Q$ defined by point evaluation and a unitary $w\in (\textbf{1}_{B}\oplus \nu_{1}(\textbf{1}_{A}))M_{\infty}(B)(\textbf{1}_{B}\oplus \nu_{1}(\textbf{1}_{A}))$
such that $\|(\phi\oplus \nu_{1})(f)-w(\psi\oplus \nu_{1})(f)w^{\ast}\|<\varepsilon ~~~\forall f\in F.$

Since $\nu_{1}$ (and $\nu$ respectively) is homotopic to a homomorphism factor through  $M_{rank(P)}(\mathbb{C})$ (for $A=PM_{\bullet}(C(X))P$) or  factor through  $M_{l}(\mathbb{C})$ (for $A=M_{l}(I_{k})$), there is a unital homomorphism \\
$\nu_{2}: A\rightarrow(\nu(\textbf{1}_{A})-\nu_{1}(\textbf{1}_{A}))M_{\infty}(B)(\nu(\textbf{1}_{A})-\nu_{1}(A))$ such that
$
K_{0}(\nu_{1}\oplus \nu_{2})=K_{0}(\nu).
$

By Lemma 6.16, $\nu$ is approximately unitarily equivalent to $\nu_{1}\oplus \nu_{2}$ to within $2\varepsilon$ on $F$. Hence $\phi \oplus \nu$ is approximately unitarily equivalent to $\psi\oplus \nu$ on $F$ to within $2\varepsilon+\varepsilon+2\varepsilon=5\varepsilon$.\\
\end{proof}

\fi

\noindent\textbf{Theorem 6.18.} Let $A=PM_{\bullet}(C(T_{\uppercase\expandafter{\romannumeral2},k}))P$ or $M_{l}(I_{k})$, $\varepsilon> 0$, and finite set $F\subset A$ with $\omega(F)<\varepsilon$. There is an $\eta> 0$ posessing the following property.

Let $\kappa$ be any fixed simplicial structure of $Sp(A)=X$, where $X=T_{\uppercase\expandafter{\romannumeral2},k}$ or $X=[0,1]$. For any $\delta>0$, there exist an integer $L>0$ and a finite set $H\subseteq AffTA(=C_{\mathbb{R}}(T_{\uppercase\expandafter{\romannumeral2}
,k}))$ or $C_{\mathbb{R}}[0,1])$ such that the following statement holds.

 Suppose that $X_{1},X_{2},\cdots,X_{s}\subset Sp(A)$ are connected  sub-complexes of $(X, \kappa)$,  $\phi_{i},\psi_{i}:A|_{X_{i}}\rightarrow B(\in \mathcal{HD})$ are homomorphisms satisfying\\
(a)  $\phi_{i}$ has $sdp(\frac{\eta}{32},\delta)$;\\
(b)  For each $i$, either both $\phi_{i}$ and $\psi_{i}$ are at least $L$-large or both $\phi_{i}$ and $\psi_{i}$ are defined by point evaluations;\\
(c)  $KK(\phi_{i})=KK(\psi_{i})$ and in particular,  $\phi_{i}(\textbf{1}_{A|_{X_{i}}})$ are unitarily equivalent to $\psi_{i}(\textbf{1}_{A|_{X_{i}}})$;\\
(d)  For $P_{i}=\phi_{i}(\textbf{1}_{A|_{X_{i}}})$, $Q_{i}=\psi_{i}(\textbf{1}_{A|_{X_{i}}}),$ the maps
$$
AffT\phi_{i}: AffT(A|_{X_{i}})\rightarrow AffTP_{i}BP_{i}~~~\mbox{ and}~~  AffT\psi_{i}: AffT(A|_{X_{i}})\rightarrow AffT Q_{i}BQ_{i}
$$ satisfy
$$
 \|AffT\phi_{i}(h|_{X_i})-AffT\psi_{i}(h|_{X_i})\|<\frac{\delta}{4},
$$
where  $AffTP_{i} BP_{i}$ is canonically identified with $AffTQ_{i}BQ_{i}$ by a unitary conjugacy;\\
%Regarding
%$$
%AffT\phi_{i}, AffT\psi_{i}: AffT(A|_{X_{i}})\rightarrow AffTP_{i}BP_{i}\cong  AffT Q_{i}BQ_{i},
%$$
%we have
%$$
 %\|AffT\phi_{i}(h|_{X})-AffT\psi_{i}(h|_{X})\|<\frac{\delta}{4}.
%$$
(e)  $\{\phi_{i}(\textbf{1}_{A|_{X_{i}}})\}_{i=1}^{s}$ and $\{\psi_{i}(\textbf{1}_{A|_{X_{i}}})\}_{i=1}^{s}$ are two sets of mutually orthogonal projections in $B$ and $\phi,\psi:A\rightarrow B$ are defined by
$$
 \phi(f)=\sum_{i=1}^{n}\phi_{i}(f|_{X_{i}}), \psi(f)=\sum_{i=1}^{n}\psi_{i}(f|_{X_{i}})~~~~\forall~ f\in A.
$$
%for all $f\in A$.\\
Then there is a unitary $u\in B$ such that
$$
\|\phi(f)-u\psi(f)u^{\ast}\|<5\varepsilon ~~~\forall f\in F.
$$

Furthermore,  if the condition (b) above is replaced by a weaker condition:\\
~~(b$'$)  for each $i$, either both $\phi_{i}$ and $\psi_{i}$ are $L$-large or there exist $ \phi'_{i}$ and $\psi'_{i}$ both defined by point evaluations such that $\|\phi_{i}(f|_{X_{i}})-\phi'_{i}(f|_{X_{i}})\|<\varepsilon$ and $\|\psi_{i}(f|_{X_{i}})-\psi'_{i}(f|_{X_{i}})\|<\varepsilon$ for all $f\in F$,\\
 then we have a weaker estimation
 $$\|\phi(f)-u\psi(f)u^{\ast}\|< 7\varepsilon~~~~\forall f\in F.$$
\begin{proof}
Applying the dilation lemma (Proposition 1.10), one can reduce the case of \linebreak
$PM_{\bullet}(C(T_{\uppercase\expandafter{\romannumeral2}
,k}))P$ to the case of $M_{\bullet}(C(T_{\uppercase\expandafter{\romannumeral2}
,k}))$ (see the proof of [EGL 2.14] also).
By (c), up to conjugacy  unitary, we can assume $\phi(\textbf{1}_{A_{|_{X_{i}}}})=\psi(\textbf{1}_{A_{|_{X_{i}}}})$. This reduces the case to the case of only one space $X_1$ (or $s=1$). Hence we will denote the homomorphisms $\phi_1$ and $\psi_1$ by $\phi$ and $\psi$, respectively, and reserve  $\phi_1$ and $\psi_1$ for other purpose.

If both $\phi$ and $\psi$  are defined by point evaluations, then  the conclusion follows  from 6.16. So we only need to prove the case
%can assume
that $\phi$ and $\psi$ are $L$-large. If both $A$ and $B$ are dimension drop algebras, one can apply Theorem 6.11 and Proposition 4.14 to get the conclusion. Hence, we assume one of them is not a dimension drop algebra.

For $A=M_{\bullet}(C(T_{\uppercase\expandafter{\romannumeral2},k}))$,  let $\eta$ be as in Theorem 4.18
and let $H\subseteq AffTA$ be also as in Theorem 4.18 for the given simplicial structure $\kappa$ of $T_{\uppercase\expandafter{\romannumeral2},k}$. Using the conditions (a), (c), (d) and $L$-large, by Theorem 4.18, we can decompose the homomorphisms $\phi$ and $Adu\circ \psi$ (for certain $u$ with $\phi(\textbf{1}_{A|_{X_{1}}}) =u^*(\psi(\textbf{1}_{A|_{X_{1}}}))u$)  simultaneously into three parts as below.
%(Here,  we apply Theorem 5.31 for the maps from $M_{\bullet}(C(X))$ to $M_{l}(I_{k})$,  apply [G5, Theorem 4.42] for the  maps from $M_{\bullet}(C(X))$ to $QM_{\bullet}(C(Y))Q$, or apply  Theorem 5.1 the  maps from  from $M_{l}(I_{k})$ to $QM_{\bullet}(C(Y))Q$, to reduce the maps  to  maps with major part factors through $M_{lk}([0,1])$.)

%It follows from (a), $L$-large, (c) and (d), applying Theorem 5.31 (for the part from $M_{\bullet}(C(X_{i}))$ to $M_{l}(I_{k})$) or [G5, Theorem 4.42] (for the part from $M_{\bullet}(C(X_{i}))$ to homongeneous algebra) or Theorem 5.1
%(which reduce the map from $M_{l}(I_{k})$ to $PM_{\bullet}(C(X_{i}))P$ to a map with major part factors through $M_{l_{k}}([0,1])$), we know that $\phi$ and $Adu\circ \psi$ (for certain $u$) can be decompose as three parts as below.

Part(1): $p_{0}\phi(f)p_{0}$, $p_{0}Adu\circ \psi(f)p_{0}$, where $p_{0}< \phi(\textbf{1}_{A|_{X_1}})= (Adu\circ \psi)(\textbf{1}_{A|_{X_1}})$ is a subprojection,  as $Q_{0}$ in Theorem 4.18.

Part(2):  $\phi_{1}$---a same homomorphism for both $\phi$ and $Adu\circ\psi$,  defined by point evaluations with $\phi_{1}(1)\geq Jp_{0}$ for any given $J$.

Part(3): $\phi_{2}$---a same homomorphism for both $\phi$ and $Adu\circ \psi$,  which factors through a direct sum of matrix algebras over intervals or points.

Then by 6.17,  applying  to $p_{0}\phi(f)p_{0}\oplus\phi_1(f)$ and $p_{0}(Adu\circ\psi(f))p_{0}\oplus\phi_1(f)$, we get the desired estimation.

For $A=M_l(I_k)$ and $B$  being homogeneous:~if $X_1\varsubsetneqq [0,1]$, then this is the conclusion of Theorem 6.11.

Let us assume that $X_1=[0,1]$,  $A|_{X_1}= A$, and $\phi, \psi: A\to B$ are  unital homomorphisms. By [Li2, 2.15],  the conditions (a) and (d) (for suitable choice of $H$) imply that $\psi$ has property $sdp(\eta/8, \delta)$. By Theorem 4.16, there are two sets of three projections $(P, P_0, P_1)$ and $(Q, Q_0, Q_1)$, with $P+P_0+P_1=Q+Q_0+Q_1={\bf 1}_B$ and $0<rank (P_i)\leq kl$ and $0<rank (Q_i)\leq kl$ for $i=0,1$,  and there are unital homomorphisms ${\tilde \phi}: M_{lk}(C[0,1])\to PBP$ and ${\tilde \psi}: M_{lk}(C[0,1])\to QBQ$ such that $\phi(f)$ and $\psi(f)$ are close to $\phi'(f)$ and $\psi'(f)$, respectively,  on $f \in F$, where
$$\phi'(f)=f(\underline{0})\cdot P_0+f(\underline{1})\cdot P_1+({\tilde \phi}\circ\imath)(f)~\mbox{and}~\psi'(f)=f(\underline{0})\cdot Q_0+f(\underline{1})\cdot Q_1+({\tilde \psi}\circ\imath)(f)~~~\forall f\in A$$
(see 5.2 for the convention). Without lose of generality, we  assume $\phi=\phi'$ and $\psi=\psi'$. It is easy to see from  the proof of Theorem 3.1 (Theorem 4.16 above) of [JLW] , we can assume $Sp({\tilde \phi}), Sp({\tilde \psi})\subset [a,1]$ for certain $a>0$. By Theorem 6.8, by increasing the rank of $P_1$ and $Q_1$ to
$$kl< rank (P_1)\leq 2kl ~~~\mbox{and}~~~~kl< rank (Q_1)\leq 2kl,$$
we can assume that both ${\tilde \phi}(e_{11})$ and ${\tilde \psi}(e_{11})$ (and $P$ and $Q$) are trivial projections. (Here $e_{11}\in M_{lk}(\C)$~($\subseteq M_{lk}(C[0,1])$) is the matrix unit of upper left corner.) Consequently, $KK({\tilde \phi})=KK({\tilde \psi})$. Up to unitary equivalence, we can assume $P=Q$ and $P_0+P_1= Q_0+Q_1=R$. Then $\phi={\tilde \phi}_0\oplus {\tilde \phi}$  and $\psi={\tilde \psi}_0\oplus {\tilde \psi}$, where  ${\tilde \phi}_0, {\tilde \psi}_0: A \to RBR$ are defined by
$${\tilde \phi}_0(f)=f(\underline{0})\cdot P_0+f(\underline{1})\cdot P_1~~\mbox{and}~~{\tilde \psi}_0(f)=f(\underline{0})\cdot Q_0+f(\underline{1})\cdot Q_1~~~~\forall f\in A.$$
Since $rank (P)\gg rank (R)$ if the number $L$ is large enough, by  increasing the number $L$ to $2L$, we can assume ${\tilde \phi}$ and ${\tilde \psi}$ still satisfy the conditions (a),  (c),  (d) and $L$ large in (b). Apply Theorem 4.17 to ${\tilde \phi}$ and ${\tilde \psi}$
to decompose ${\tilde \phi}$ and $Ad u\circ{\tilde \psi}$ (for certain unitary $u$), simultaneously, into three parts: $\phi_0=Q'_0{\tilde \phi}Q'_0$ ($\psi_0=Q'_0 (Ad u\circ{\tilde \psi})Q'_0$),
$\phi_1$ and $\phi_2$ (the last two parts are same for ${\tilde \phi}$ and $Ad u\circ{\tilde \psi}$). Then by 6.17 applying to ${\tilde \phi}_0(f)\oplus \phi_0(f)\oplus\phi_1(f)$ and ${\tilde \psi}_0(f)\oplus \psi_0(f)\oplus\phi_1(f)$, we get the desired estimation.

%Since $KK(\phi_1)=KK(\psi_1)$ from (c), we can assume that  $\phi_1( e_{ij}\otimes {\bf 1}_k)=\psi_1( e_{ij}\otimes {\bf 1}_k)$ for all matrix units $e_{ij} \in M_l$. Write $\phi_1=\phi'\otimes id_l$ and $\psi_1=\psi'\otimes id_l$, with $\phi', \psi': I_k \to qBq$, where $q=\phi_1( e_{ij}\otimes {\bf 1}_k)$.

%Finally for the partial map from dimension drop to dimension drop algebras, one can apply Theorem 7.11 and proposition 4.14.\\
\end{proof}
The following theorem is the main theorem of this section.

\noindent\textbf{Theorem 6.19.} Let $A=lim(A_{n},\phi_{nm})$ be (not necessarily  unital) $\mathcal{HD}$ inductive limit with the ideal property.
Suppose that $B=\bigoplus_{i=1}^{s} B^{i}\in \mathcal{HD}$ and $\xi: B\rightarrow A_{n}$ is a homomorphism. Let $\varepsilon>0$, $F=\bigoplus_{i=1}^{s}F^{i}\subset \bigoplus_{i=1}^{s}B^{i}$ be a finite set such that if $B^{i}=PM_{\bullet}(C(T_{\uppercase\expandafter{\romannumeral2}
,k}))P$ or $M_{l}(I_{k})$, then $F^{i}$ is weakly approximately constant to within $\varepsilon$ (i.e.,  $\omega(F^{i})<\varepsilon$).

There  exist a positive integer $m>n$,  two finite sets
$$H=\oplus H^{i}\subseteq AffTA_{m}=\oplus AffTA_{m}^{i}~~~\mbox{and}~~~E=\oplus E^{i}\subset U (A_{m})/ \widetilde{SU(A_{m})}=\bigoplus_{i}U(A_{m}^{i})/ \widetilde{SU(A_{m}^{i})},$$
and a number $\sigma> 0$ satisfying the following condition.

Suppose that two homomorphisms
$\psi_{1},\psi_{2}:A_{m}\rightarrow D=\bigoplus\limits_{i=1}\limits^{s}D^{i}(\in \mathcal{HD})$
satisfy  the following conditions:\\
(1)  $KK(\psi_{1})\!=\!KK(\psi_{2})$, consequently, there is a unitary $W\!\in\! D$ such that $\psi_{1}(\textbf{1}_{A_{m}^{i}})\!=\!W^{\ast} \psi_{2}(\textbf{1}_{A_{m}^{i}})W$ for all $i$;\\
(2)  For each pair $A_{m}^{i}$, $D^{j}$, denote $Q=\psi_{1}^{i,j}(\textbf{1}_{A_{m}^{i}})=W^{\ast}\psi_{2}^{i,j}(\textbf{1}_{A_{m}^{i}})W$
%and consider $\psi_{1}^{i,j}$, $AdW\circ \psi_{2}^{i,j}:A_{m}^{i}\rightarrow QD^{j}Q$, we have
$$\|AffT\psi_{1}^{i,j}(h)- AffT(Ad W\circ\psi_{2}^{i,j})(h)\|<\sigma ~~~~\forall h\in H^i,$$
where $AffT\psi_{1}^{i,j}$ and $AffT(AdW\circ\psi_{2}^{i,j})$ are regarded as unital maps from
$AffTA_{m}^{i}$ to $AffTQD^{j}Q$\\ (see 2.4);\\
(3) ~~~~~~~~~~~~~~~$dist((\psi_{1}^{i,j})^{\natural}(g),(AdW\circ \psi_{2}^{i,j})^{\natural}(g))<\sigma ~~~~\forall g\in E^{i},$\\
where
both $((\psi_{1}^{i,j})^{\natural})$ and $(AdW\circ \psi_{2}^{i,j})^{\natural}$ are regarded as maps from $U (A_{m}^{j})/\widetilde{SU (A_{m}^{j}})$
to $U (QD^{j}Q)/$ $ \widetilde{SU}(QD^{j}Q)$ (see 2.8). \\
It follows that  there is a unitary $u\in D$
such that $$\|\psi_{1}\circ\phi_{nm}\circ\xi(f)-u^{\ast}(\psi_{2}\circ\phi_{nm}\circ\xi)(f)u\|<7\varepsilon ~~~~~\forall f\in F.$$
(If we identify $AffT(WQW^{\ast})D^{j}(WQW^{\ast})$ with $AffTQD^{j}Q$,  and identify
$U(QD^{j}Q)/ \widetilde{SU}(QD^{j}Q)$
with $U((WQW^{\ast})D^{j}(WQW^{\ast}))/ \widetilde{S U}((WQW^{\ast})D^{j}(WQW^{\ast}))$
in the obvious way, we can write the above (2) and (3) as
$
 \|AffT\psi_{1}^{i,j}(h)-AffT\psi_{2}^{i,j}(h)\|<\sigma~~\forall h\in H^i,
$
~~and~~
$
 dist((\psi_{1}^{i,j})^{\natural}(g),(\psi_{2}^{i,j})^{\natural}(g))<\sigma~~ \forall g\in E^i.)
$
\begin{proof}
Let $\eta>0$ be the number  desired in Corollary 6.13,  Theorem 6.11  and Theorem 6.18 for each block $B^{i}$ with finite set $F^{i}$.  Further, we assume that $\eta$ is so small that $F^i$ is $(\frac{\varepsilon}3,\eta)$ continuous, as we will also apply Theorem 4.14.
Then consider
$$
B\xrightarrow{\xi} A_{n}\rightarrow A_{n+1}\rightarrow \cdots
$$
as an inductive limit with $B$ as the first term. Let $\kappa$ be a fixed simplicial decomposition of $\coprod Sp B^i$
such that every simplex of $(\coprod Sp B^i, \kappa)$  has diameter less than $\frac{\eta}{128}$.
Applying Theorem 4.4 to $B$ (in place of $A_n$) and $\frac{\eta}{32}$ (in place of $\eta$), there is an   $m_{1}> n$ (large enough), and there is  a factorization of $\phi_{n, m_{1}}\circ\xi$
as
$$
(\phi_{n, m_{1}}\circ\xi):~~B\stackrel{\pi_{1}}{\longrightarrow}B_{1}\stackrel{\widetilde{\phi}}{\longrightarrow} A_{m_{1}}
$$
such that each $B^{i}$ corresponds to possibly  several blocks $B_{1}^{i,j}=B^{i}|_{X_{i,j}}$, where $X_{i,j}\subset Sp(B^{i})$ are connected sub-complexes of $(\coprod Sp B^i, \kappa)$ and $\widetilde{\phi}^{(i,j),k}:B_{1}^{i,j}\rightarrow A_{m_{1}}^{k}$ is either zero map or has $sdp(\frac{\eta}{32},\delta)$ property for certain $\delta$, where $\pi_{1}:B=\oplus B^{i}\rightarrow B_{1}=\oplus B_{1}^{i,j}$ is given by the  restrictions $B^{i}\rightarrow\oplus_{j} B_{1}^{i,j}$.

For each block $B^{i}$, we can assume that for each $j>1$, $Sp(B_{1}^{i,j})= X_{i.j}\varsubsetneqq Sp(B^{i})$,  that is, only the first block $B_{1}^{i,1}$ may have the same spectrum as $B^{i}$.
For each $B_{1}^{i,1}$ with $Sp(B^{i,1})=Sp(B^{i})=S^{1}$, and $F^{i}\subset B^{i}$, $\eta> 0$, $\delta> 0$,  let $\widetilde{\eta}$ be as in Corollary 6.13 corresponding  to $\eta$, $\delta$ and $F^{i}$ (works for all such $F^{i}$).
By Theorem 4.4 again, there is $m_{2}>m_{1}$ such that the homomorphism $\phi_{m_{1},m_{2}}\circ \widetilde{\phi}|_{B_{1}^{i,1}}$ factors through
$$
B_{1}^{i,1}\stackrel{\pi_{2}}{\longrightarrow}\oplus_{j}B_{2}^{i,1,j}\stackrel{\widetilde{\widetilde{\phi}}}{\longrightarrow}A_{m_{2}}
$$
Such that ${\widetilde{\widetilde{\phi}}}^{(i,1,j),k}:~~B_{2}^{i,1,j}\rightarrow A_{m_{2}}^{k}$ is either zero map or has
$sdp(\frac{\widetilde{\eta}}{32},\widetilde{\delta})$ property for certain $\widetilde{\delta}$, where $\pi_{2}$ is defined in the  same way as that of $\pi_{1}$.

Let $L$ be larger than the numbers $L$  described in Corollary 6.13 and Theorem 6.18, namely,  the number $L$ in Corollary 6.13 for $F^{i}\subset M_{l_{i}}(C(S^{1}))$,
$\varepsilon$, $(\eta,\delta)$,$(\widetilde{\eta},\widetilde{\delta})$,   and the number $L$  in Theorem 6.18 for $F^{i}\subset M_{l'_{i}}(I_{_{n_{i}}})$ or $F^{i}\subset PM_{\bullet}(C(T_{\uppercase\expandafter{\romannumeral2}
,k}))P$, $\varepsilon$, $(\eta,\delta)$  and the above simplicial structure $\kappa$.

Denote those blocks of $B^{i}$ with $B^i=PM_{\bullet}(C(T_{\uppercase\expandafter{\romannumeral2},k})P$, $M_l(I_k)$ or $M_l(\C)$
%$Sp(B^{i})\neq S^{1}$
by $C_{4}^{i}$; ~blocks $B_{1}^{i,j}$ with $B^i=M_{\bullet}(C[0,1])$ or $B^i=M_{\bullet}(C(S^1))$ but $Sp(B_{1}^{i,j})\neq S^{1}$ by $C_{3}^{i}$; ~blocks $B_{2}^{i,1,j}$ with $Sp(B_{1}^{i,1})=S^{1}$ but $Sp(B_{2}^{i,1,j})\neq S^{1}$ by $C_{2}^{i}$, ~ and blocks $B_{2}^{i,1,1}$ with $Sp(B_{2}^{i,1,1})= S^{1}$ by $C_{1}^{i}$.

We can order the blocks  as $1\leq i\leq t_{1}$ for $C_{1}^{i}$,  and $t_{1}< i\leq t_{2}$ for $C_{2}^{i}$, and  $t_{2}< i\leq t_{3}$ for $C_{3}^{i}$,
 and finally $t_{3}< i\leq t_{4}$ for $C_{4}^{i}$ and let
$$
C'=\oplus C^{i}=(\oplus_{i=1}^{t_1}C_{1}^{i})\oplus(\oplus_{i=t_{1}+1}^{t_{2}}C_{2}^{i})\oplus(\oplus_{i=t_{2}+1}^{t_{3}}C_{3}^{i})\oplus(\oplus_{i=t_{3}+1}^{t_{4}}C_{4}^{i})~~.
$$
Let $\phi': C'\rightarrow A_{m_{2}}$ be defined as follows:

For $C_{1}^{i}$, $\phi'|_{C_{1}^{i}}=\widetilde{\widetilde{\phi}}|_{C_{1}^{i}}$ which has both properties $sdp(\frac{\eta}{32},\delta)$ and $sdp(\frac{\widetilde{\eta}}{32},\widetilde{\delta})$.

For $C_{2}^{i}$, $\phi'|_{C_{2}^{i}}=\widetilde{\widetilde{\phi}}|_{C_{2}^{i}}$, which has property $sdp(\widetilde{\frac{\eta}{32}},\widetilde{\delta})$. (Note that $Sp(C_{2}^{i})\subsetneqq S^{1}$  is either an interval or $\{pt\}$.)

For $C_{3}^{i}$, $\phi'|_{C_{3}^{i}}=\phi_{m_{1}n_{2}}\circ \widetilde{\phi}|_{C_{3}^{i}}$, which has property $sdp(\frac{\eta}{32},\delta)$.  (Also  $Sp(C_{3}^{i})\subsetneqq S^{1}$ or $Sp(C_3^i)\subseteq [0,1]$ is either an interval or $\{pt\}$.)

For $C_{4}^{i}$, $\phi'|_{C_{4}^{i}}=\phi_{n m_{2}}\circ \xi|_{C_{4}^i}$ which factors through $\oplus_{j}C_{4}^{i}|_{X_{i,j}}$ with $X_{i,j}\subset Sp(C_{4}^i)$.
%$(=[0,1])$, $T_{\uppercase\expandafter{\romannumeral2},k}$, or $\{pt\}$ but not $S^{1}$).
(That is,  we will consider this part of homomorphism as a whole --- not dividing them into pieces, but we will use the property that $\phi|_{C_{4}^{i}}$  factors through as $\oplus C_{4}^{i}|_{X_{i,j}}$ with each $X_{i,j}$ being sub-complex with respect to $\kappa$, and with each partial map having property $sdp(\frac{\eta}{32},\delta)$  as in Theorem 6.18.)

For the integer $L$, we can find $m>m_{2}$ such that $\phi_{m_{2}m}\circ \phi':C\rightarrow A_{m}$ satisfies the dichotomy of Proposition 3.2.

For each $C_{1}^{i}$ and $A_{m}^{j}$, if $(\phi_{m_{2}, m}\circ \phi')^{i,j}$ is not  $L$-large, then the homomorphism $(\phi_{m_{2}m}\circ \phi')^{i,j}$ satisfies (b) of Proposition 3.2 and factors as
 $$(\phi_{m_{2}m}\circ \phi')^{i,j}:~~C_{1}^{i}\rightarrow C_{1}^{i}|_{X_{i,j}}\xrightarrow{\widetilde{\widetilde{\widetilde{\phi}}}}A_{m}^{j},$$
where $X_{i,j}$ is a union of several disjoint small intervals centered at ${x_{1},x_{2},\cdots,x_{\bullet}}$ with length smaller than any pregiven number, and  there exists a homomorphism $$\psi:C_{1}^{i}\rightarrow C_{1}^{i}|_{\{x_{1},x_{2},\cdots,x_{\bullet}\}}\rightarrow A_{m}^{j}~~~~~~~~\mbox{with}~~~~
 \|\widetilde{\widetilde{\widetilde{\phi}}}(f)-\psi(f)\|<\varepsilon~~~~~\forall f\in F^{i'}\subset B^{i'}=C_{1}^{i}$$ (note that $C_{1}^{i}$ is one of the original $B^{i'}$ with $Sp(B^{i'})=S^{1}$).
Denote  these blocks $C_{1}^{i}|_{X_{i,j}}$ by $C_{5}^{k}$ for $t_{4}+1\leq k\leq t_{5}$. Then $Sp(C_{5}^{i})$ is the disjoint union of some small intervals.  Finally,  let $C=C'\oplus(\oplus_{i=t_{4}+1}^{t_{5}}C_{5}^{i})$,
and $\phi: C\rightarrow A_{m}$ be defined as below.
For $C_{2}^{i}$, $C_{3}^{i}$, $C_{4}^{i} $ and any block $A_{m}^{j}$, we define  $\phi^{i,j}=(\phi_{m_{2}m}\circ \phi')^{i,j}$.  For $C_{1}^{i}$ and $A_{m}^{j}$, if  the homomorphism $(\phi_{m_{2}m}\circ \phi')^{i,j}$ is at least $L$-large, define $\phi^{i,j}=(\phi_{m_{2}m}\circ \phi')^{i,j}$;  if $(\phi_{m_{2}m}\circ \phi')^{i,j}$ not $L$-large, then define $\phi^{i,j}$ to be zero (this part of the map will be defined to go  through $C_{5}^{i}$). For $C_{5}^{i}$, define $\phi|_{C_{5}^{i}}=\widetilde{\widetilde{\widetilde{\phi}}}|_{C_{5}^{i}}$. Now we obtain the following factorization
$$
\phi_{nm}\circ\xi:
B\stackrel{\pi}{\longrightarrow}C\stackrel{\phi}{\longrightarrow}A_{m}
$$
with the properties described as below. For each block $C^{i}$ of $C$, there is only one block $B^{i'}$
such that $\pi^{i',i}\neq 0$.
Furthermore for such block $B^{i'}$, $C^{i}=B^{i'}|_{X}$ for certain $X\subset Sp(B^{i'})$ and the map $\pi^{i',i}$ is the restriction. And there are four kinds of blocks $C^i$ of $C$ described as follows:\\
(i)~~ The first kind block  from $C_1^i$ with $Sp(C^{i})=S^{1}$:
For each such block $C^{i}$, there is a block $B^{i'}$ with $B^{i'}=C^{i}$,  the
map $\pi^{i',i}:B^{i'}\rightarrow C^{i}$ is the identity map, and furthermore for each $A_{m}^{j}$, $\phi^{i,j}$ is either zero map, or is an $L$-large map with $sdp(\frac{\eta}{32},\delta)$ and $sdp(\frac{\widetilde{\eta}}{32},\widetilde{\delta})$ properties.\\
(ii) ~~The second kind of block from $C_5^i$ with $Sp(C^{i})$ being finite union of intervals or points and with the following property:  for each $\phi^{i,j}$, there exists a $\psi^{i,j}$ defined by point evaluations such that
$$\|\phi^{i,j}(f)-\psi^{i,j}(f)\|<\varepsilon ~~~~~~\forall f\in \pi(F)^{(i)},$$
where $\pi(F)^{(i)}$ denote $\pi^{i',i}(F^{i'})$ for $F^{i'}\subset B^{i'}$ with $\pi^{i',i}\neq 0$.\\
(iii)~~ The third kind of block from $C_2^i$ or $C_3^i$ with  $Sp(C^{i})=$interval or $\{pt\}$: For each such block $C^i$ and any block $A_m^j$, if $\phi^{i,j}\neq 0 $, then $\phi^{i,j}$ has either property $sdp(\frac{\eta}{32},\delta)$ (if $C^i$ is from $C_3^i$) or property $sdp(\frac{\widetilde{\eta}}{32},\widetilde{\delta})$ (if $C^i$ is from $C_2^i$).\\
%which is image of $\pi^{{i',i}}:B^{i'}\rightarrow C^{i}$ for some $B^{i'}=M_{\bullet}(C[0,1])$ or $M_{\bullet}(C(S^{1}))$, and satisfy the condition that if $\phi^{i,j}\neq 0 $, then $\phi^{i,j}$ has either property $sdp(\frac{\eta}{32},\delta)$ (for those blocks coming from $B^{i'}=M_{\bullet}(C[0,1])$ and some coming from $B^{i'}=M_{\bullet}(C(S^{1}))$ which corresponding to factorization involving $\widetilde{\phi}$) or has $sdp(\frac{\widetilde{\eta}}{32},\widetilde{\delta})$ (for some blocks coming from $M_{\bullet}(C(S^{1}))$ corresponding to factorization involve $\widetilde{\widetilde{\phi}}$).\\
(iv) ~~The fourth kind of block from $C_4^i$ with  $C^{i}=B^{i'}=PM_{\bullet}C(T_{\uppercase\expandafter{\romannumeral2},
k})P$,  $M_{l}(I_{k})$ or $M_l(\C)$: For each such block and any $A_{m}^{j}$,  if $\phi^{i,j}\neq 0$, then $\phi^{i,j}$ is a homomorphism satisfying the conditions (a) and  (b') (for $\phi$) of Theorem 6.18.  That is,  $\phi^{i,j}=\oplus_k(\phi^{i,j})_{k}$, where  for each $k$
$$(\phi^{i,j})_{k}:~~ C^{i}|_{X_{k}}\rightarrow A_{m}^{j}$$ has $sdp(\frac{\eta}{32},\delta)$ property (this is (a) of Theorem 6.18), $X_k$ is a sub-complex with respect to $\kappa$,  and $(\phi^{i,j})_{k}$ is either $L$-large or close to a homomorphism defined by point evaluations to within $\varepsilon$ on $F$.

Let $H^{\prime i}\subset AffTB^{i}$ be defined to include all the subsets of the following forms.

For $B^{i}=M_{\bullet}(C(S^{1}))$, the set $H^{\prime i}$ should  include $H$   in Corollary 6.13 for $(\eta,\delta)$ and $(\widetilde{\eta},\widetilde{\delta})$;  include $H$  in 6.11 for $(\eta,\delta)$; and  include $H$  in 6.11 for $(\widetilde{\eta},\widetilde{\delta})$, play the role of $(\eta,\delta)$.

For $B^{i}=M_{\bullet}(C[0,1])$, the set $H^{\prime i}$ should include $H$  in 6.11 for $(\eta,\delta)$.

For $B^{i}=PM_{\bullet}(C(T_{\uppercase\expandafter{\romannumeral2}
k}))P$, the set $H^{\prime i}$ should include $H$ as in 6.18 for $(\eta,\delta)$ and simplicial structure $\kappa$.

For $B^{i}=M_{l}(I_{k})$,  the set $H^{\prime i}$ should include $H$  as in 6.18 for $(\eta,\delta)$ and simplicial structure $\kappa$; and $H$  in 4.14 for $(\eta,\delta)$ (see Remark 4.15  also).

Let $H^{\prime}=\bigoplus_{i}H^{\prime i}\subset AffTB=\bigoplus AffTB^{i}.$

For each block $B^{i}=M_{\bullet}(C(S^{1}))$, let $E^{\prime i}=\{z\}\subset U (B^{i})/ \widetilde{SU(B^{i})} $, where $z$ is defined in 6.13. For blocks $B^{i}$ with $Sp(B^{i})\neq S^{1}$, let $E^{\prime i}=\{0\}\subset U(B^{i})/\widetilde{SU}(B^{i})$ , where $0\in U(B^{i})/\widetilde{SU}(B^{i})$ is represented by ${\bf 1}_{B^i}$. Let $\delta'$ be as described in 6.13 (for all blocks of $B^{i}=M_{\bullet}(C(S^{1}))$ for $(\eta,\delta)$).

Let $E^{\prime}=\bigoplus E^{\prime i}\subset U(B)/\widetilde{SU}(B)=\bigoplus U(B^i)/\widetilde{SU}(B^{i}).$

Denote  $\phi_{nm}\circ \xi ({\bf 1}_B)=\bar{R}$. We recall that the unital positive map $AffT(\phi_{n,m}\circ\xi): AffTB \to AffT \bar{R}A_m\bar{R}$ and the  contractive group homomorphism $(\phi_{n,m}\circ\xi)^{\natural}: U(B)/ \widetilde{SU}(B) \to U(\bar{R}A_m\bar{R})/ \widetilde{SU}(\bar{R}A_m\bar{R}) $ are described in (2.4) and (2.8).
Also recall that $\imath_T: AffT \bar{R}A_m\bar{R} \to AffT A_m$ and $\imath_*: \widetilde{SU}(\bar{R}A_m\bar{R}) \to \widetilde{SU}(A_m) $ are induced by inclusion as in (2.4) and (2.8).

For each $A_{m}^{j}$, let $H^{j}=\pi_{j}((\imath_T\circ AffT(\phi_{n,m}\circ\xi))(H^{\prime}))=\pi_{j}((\imath_T\circ AffT(\phi\circ\pi))(H^{\prime}))$
and $E^{j}=\pi_{j}((\imath_*\circ(\phi_{n,m}\circ\xi)^{\natural})(E^{\prime}))$.
Set  $H=\bigoplus H^{j}$ and $E=\bigoplus E^{j}$.
Let $N=\max\limits_{j}(rank(\textbf{1}_{A_{m}^{j}}))$ and $\sigma=\frac{1}{N}min(\frac{\delta}{4},\frac{\widetilde{\delta}}{4},\delta')$.

Let $\psi_{1},\psi_{2}: A_{m}\rightarrow D=\bigoplus\limits_{i=1}\limits^{s}D^{i}(\in \mathcal{HD})$ be two homomorphisms satisfying the condition of the  theorem  for the set $H, E$  and positive number $\sigma$.

For any $C^{i}\subset C$, $A_{m}^{j}\subset A_{m}$, with $\phi^{i,j}\neq0$, let
$R=\phi^{i,j}(\textbf{1}_{C^{i}})\in A_{m}^{j}$. Then $\frac{rank(R)}{rank(1_{A_{m}^{j}})}\geq\frac{1}{N}$.

Up to unitary equivalence, we can assume $\psi_{1}(R)=\psi_{2}(R)$ and $\psi_{1}(\textbf{1}_{A_{m}^{j}})=\psi_{2}(\textbf{1}_{A_{m}^{j}})$.

For a block $D^{k}$, let $Q_{1}=\psi_{1}^{j,k}(R)=\psi_{2}^{j,k}(R)$ and $Q=\psi_{1}^{j,k}(\textbf{1}_{A_{m}^{j}})=\psi_{2}^{j,k}(\textbf{1}_{A_{m}^{j}})$.

Note that the inclusion $RA_{m}^{j}R\rightarrow A_{m}^{j}$ induces the map  $$AffTRA_{m}^{j}R=C_{\mathbb{R}}(Sp(A_{m}^{j})) \longrightarrow AffTA_{m}^{j}=C_{\mathbb{R}}(Sp(A_{m}^{j}))$$ by sending $f$ to $\frac{rank(R)}{rank(\textbf{1}_{A_{m}^{j}})} f$.

Any two  unital homomorphisms
$\psi_{1},\psi_{2}:A_{m}^{j}\rightarrow QD^{k}Q$ (with $\psi_{1}(R)=\psi_{2}(R)=Q_1$ ) and the corresponding homomorphisms
$\widetilde{\psi}_{1}=\psi_{1}|_{RA_{m}^{j}R},\;\widetilde{\psi}_{2}=\psi_{2}|_{RA_{m}^{j}R}:RA_{m}^{j}R\rightarrow Q_{1}D^{k}Q_{1}$ have the following relations:\\
 For  any $f\in AffTRA_{m}^{j}R$ with $\imath_T(f) \in AffTA_{m}^{j}$
 (see 2.4),
$$\|AffT\widetilde{\psi}_{1}(f)-AffT\widetilde{\psi}_{2}(f)\|=\frac{rank(\textbf{1}_{A_{m}^{j}})}{rank(R)}\|AffT\psi_{1}(\imath_T(f))-AffT\psi_{2}(\imath_T(f))\|;$$
and for any $g\in U(RA_{m}^{j}R)/ \widetilde{SU}(RA_{m}^{j}R)$ with $\imath_*(g)\in U(A_{m}^{j})/ \widetilde{SU}(A_{m}^{j})$ (see 2.8),
$$
dist(\widetilde{\psi}_{1}^{\natural}(g),\widetilde{\psi}_{2}^{\natural}(g))\leq \frac{rank(\textbf{1}_{A_{m}^{j}})}{rank(R)}dist(\psi_{1}^{\natural}(\imath_*(g)),\psi_{2}^{\natural}(\imath_*(g))).
$$
%as left hand side are considered as element in $AffTQ_{1}D^{k}Q_{1}$ or $U(Q_{1}D^{k}Q_{1})/ \widetilde{SU}(Q_{1}D^{k}Q_{1})$ and right hand side are considered as element in $AffTQD^{k}Q$ or $U (QD^{k}Q)/ \widetilde{SU} (QD^{k}Q)$.

Using the above facts, it is routine to verify the conclusion of the theorem.\\
\end{proof}

%\newpage

\vspace{2mm}

\noindent\textbf{\S7. Proof of the main theorem}
%\vspace{3mm}

In this section, we will prove the following main theorem of this article.

\noindent\textbf{Theorem 7.1.} Suppose that $A\!=\!\lim(A_{n},\phi_{nm})$ and $B\!=\!\lim(B_{n},\psi_{nm})$ are two (not necessarily unital) $A\mathcal{HD}$ inductive
limit algebras  with the ideal property.  Suppose that there is an isomorphism $$\alpha:(\underline{K}(A),\underline{K}(A)^{+},\sum A)\rightarrow(\underline{K}(B),\underline{K}(B)^{+},\sum B)$$ which is compatible with Bockstein operations. Suppose that for each projection $p\in A$ and
$\overline{p}\in B$ with $\alpha([p])=[\overline{p}]$, there exist a unital positive linear isomorphism $$\xi^{p,\overline{p}}:AffT pAp\rightarrow AffT \overline{p}B\overline{p}$$ and an isometric group isomorphism $$\gamma^{p,\overline{p}}:U(pAp)/\widetilde{SU(pAp)}\rightarrow U(pBp)/\widetilde{SU(pBp)}$$
 satisfying  the following compatibility conditions: \\
(1) For each pair of projections $p<q\in A$ and $\overline{p}<\overline{q}\in B$ with $\alpha([p])=[\overline{p}], \alpha([q])=[\overline{q}]$, the  diagrams

$$\CD
  AffT(pAp) @>\xi^{p,\overline{p}}>> AffT(\overline{p}B\overline{p}) \\
  @V \ VV @V \ VV  \\
  AffT(qAq) @>\xi^{q,\overline{q}}>> AffT(\overline{q}B\overline{q})
\endCD
~~~\mbox{and}~~~
\CD
  U(pAp)/\widetilde{SU(pAp)} @>\gamma^{p,\overline{p}}>> U(\overline{p}B\overline{p})/\widetilde{SU(\overline{p}B\overline{p})} \\
  @V \ VV @V \ VV  \\
  U(qAq)/\widetilde{SU(qAq)} @>\gamma^{q,\overline{q}} >> U(\overline{q}B\overline{q})/\widetilde{SU(\overline{q}B\overline{q})}
\endCD
$$
commute, where the vertical maps are induced by the  inclusion homomorphisms.\\
(2) The maps $\alpha_{0}$ and $\xi^{p,\overline{p}}$ are compatible, that is, the diagram
$$\CD
  K_{0}(pAp) @>\rho>> AffT(pAp) \\
  @V \alpha_{0} VV @V \xi^{p,\overline{p}} VV  \\
  K_{0}(\overline{p}B\overline{p}) @>\rho>> AffT(\overline{p}B\overline{p})
\endCD
$$
commutes (this is not an extra requirement, since it follows from the commutativity of the first diagram in (1) above by [Ji-Jiang]), and then we have the
map (still denoted by $\xi^{p,\overline{p}}$) from $AffT(pAp)/\widetilde{\rho K_{0}(pAp)}$ to $AffT(\overline{p}B\overline{p})/\widetilde{\rho K_{0}(\overline{p}B\overline{p})}$.\\
(3) The maps $\xi^{p,\overline{p}}$ and $\gamma^{p,\overline{p}}$ are compatible, that is,  the diagram
$$\CD
  AffT(pAp)/\widetilde{\rho K_{0}(pAp)} @>\widetilde{\lambda}_{A}>> U(pAp)/\widetilde{SU(pAp)} \\
  @V \xi^{p,\overline{p}} VV @V \gamma^{p,\overline{p}} VV  \\
  AffT(\overline{p}B\overline{p})/\widetilde{\rho K_{0}(\overline{p}B\overline{p})} @>\widetilde{\lambda}_{B}>> U(\overline{q}B\overline{q})/\widetilde{SU(\overline{q}B\overline{q})}
\endCD$$
commutes\\
(4) The map $\alpha_{1}:K_{1}(pAp)/torK_{1}(pAp)\rightarrow K_{1}(\overline{p}B\overline{p})/tor(\overline{p}B\overline{p})$
(note that $\alpha$ keeps the positive cone of $\underline{K}(A)^{+}$ and  therefore takes  $K_{1}(pAp)\subset K_{1}(A)$ to $K_{1}(\overline{p}B\overline{p})\subset K_{1}(B)$) is  compatible with $\gamma^{p,\overline{p}}$,  that is,  the diagram
$$\CD
  U(pAp)/\widetilde{SU(pAp)} @>\widetilde{\pi}_{pAp}>> K_{1}(pAp)/torK_{1}(pAp) \\
  @V \gamma^{p,\overline{p}} VV @V \alpha_{1} VV  \\
  U(\overline{p}B\overline{p})/\widetilde{SU(\overline{p}B\overline{p})} @>\widetilde{\pi}_{pBp}>> K_{1}(\overline{p}B\overline{p})/ torK_{1}(\overline{p}B\overline{p})
\endCD$$
commutes\\
Then there is an isomorphism $\Gamma:A\rightarrow B$ such that\\
(a) $\underline{K}(\Gamma)=\alpha$,  and \\
(b) if $\Gamma_{p}:pAp\rightarrow \Gamma(p)B\Gamma(p)$ is the restriction of $\Gamma$,  then $AffT(\Gamma_{p})=\xi^{p,\overline{p}}$ and $\Gamma^{\natural}_{p}=\gamma^{p,\overline{p}}$, where $[\overline{p}]=[\Gamma(p)]$.

\noindent\textbf{Lemma 7.2.} Let $C\in \mathcal{HD}$ be a single building block with $p\in C$ a projection. Let $a=\frac{rank(\textbf{1}_{C})}{rank(p)}$, where
$rank: K_{0}(C)_{+}\rightarrow \mathbb{Z}_{+}$ is defined in 1.12. Let
$\imath_{\ast}:U(pCp)/\widetilde{SU(pCp)}\rightarrow U(C)/\widetilde{SU(C)}$
be given by $\imath_{\ast}(u)=u+(\textbf{1}_{C}-p)$ (see 2.8). Then for any two elements $u,v\in U(pCp)/\widetilde{SU}(pCp)$, we have that if
$\widetilde{D}_{C}(\imath_{\ast}(u),\imath_{\ast}(v))<\varepsilon$ then $\widetilde{D}_{pCp}(u,v)<a\varepsilon$.

\begin{proof}
 If $a\varepsilon>2$, then $\widetilde{D}_{pCp}(u,v)\leq 2 <a\varepsilon$. Let us assume $a\varepsilon\leq2$, and therefore $\varepsilon<2$. Then the element $[u^{\ast}v]\in U(pCp)/\widetilde{SU}(pCp)$ satisfies $\widetilde{\pi}_{pCp}([u^{\ast}v])=0,$ where $\widetilde{\pi}_{A}: U(A)/\widetilde{SU(A)}\rightarrow K_1(A)/torK_1(A)$.

Then there is an $h\in AffT(pCp)/\widetilde{\rho K_{0}(pCp)}$ such that $\widetilde{\lambda}_{pCp}(h)=[u^{\ast}v]$. Hence $\widetilde{\lambda}_{C}(\imath_{\ast}(h))=\imath_{\ast}([u^{\ast}v]).$
~~On the other hand, $\|\imath_{\ast}(h)\|=\frac{1}{a}\|h\|.$
We have $\widetilde{d}'(\imath_{\ast}(h),0)=\frac{1}{a}\widetilde{d}'(h,0),$
where $\widetilde{d}'$ is defined
in 2.4.
Notice the following easy fact, if $|e^{i\theta}-1|<\varepsilon$, then $|e^{ia\theta}-1|<a\varepsilon$ for $a>1$. Using the definition of $\widetilde{d}_{A}$ in 2.4 and the isometry  $\widetilde{\lambda}_{A}:AffT A/\widetilde{\rho K_{0}(A)}\rightarrow U(A)/\widetilde{SU(A)}$ in 2.7,
 it is easy to get the lemma.\\
\end{proof}

\noindent\textbf{7.3.}~~We will construct an  intertwining between the inductive limit systems  $(A_{n},\phi_{nm})$ and $(B_{n},\psi_{nm})$
(of course, we need to pass to subsequences). By Theorem 3.5, there are homomorphisms $\nu_{n}:A_{n}\rightarrow B_{n}$ and $\mu_{n}:B_{n}\rightarrow A_{n+1}$ such that\\
(a)~ $\mu_{n}\circ \nu_{n}\sim_{h}\phi_{n,n+1}$ and $\nu_{n+1}\circ \mu_{n}\sim_{h}\phi_{n,n+1}$,  and\\
(b)~ $\alpha\circ \underline{K}(\phi_{n,\infty})\!=\!\underline{K}(\psi_{n,\infty})\circ\underline{K}(\nu_{n})\!\in\! KL(A_{n},B)$ and
~~$\alpha^{-1}\circ \underline{K}(\psi_{n,\infty})\!=\!\underline{K}(\phi_{n+1,\infty})\circ\underline{K}(\mu_{n})\!\in\! KL(B_{n},A)$.

Notice  that for any  building blocks $A_{n}^{i}$, if two projections $p,q\in A_{n}^{i}$ with $[p]=[q]\in K_{0}(A_{n}^{i})$,  then $p\sim_{u}q$.  Consequently, we can modify (using unitary equivalence) $\nu_{n}$ and $\mu_{n}$ such that for any  blocks $A_{n}^{i}$ and $B_{n}^{j}$, $$\mu_{n}\circ\nu_{n}(\textbf{1}_{A_{n}^{i}})=\phi_{n,n+1}(\textbf{1}_{A_{n}^{i}})~~\mbox{and}~~~~\nu_{n+1}\circ\mu_{n}(\textbf{1}_{B_{n}^{i}})=\psi_{n,n+1}(\textbf{1}_{B_{n}^{i}}).$$
(See Remark 3.6.)

\noindent\textbf{7.4.}~~Let $F_{1}\subset A_{1}, F_{2}\subset A_{2},\cdots,F_{n}\subset A_{n}$, $\cdots\cdots$ be a sequence of finite sets such that  $F_{n}\supset\phi_{k,n}(F_{k})$
for all $k\leq n-1$, and that  $\bigcup\limits_{n=1}\limits^{\infty}\phi_{n,\infty}(F_{n})$ is dense in unit ball of $A$. Let $G_{1}\subset B_{1}, G_{2}\subset B_{2},\cdots,G_{n}\subset B_{n}$,$\cdots\cdots$ be defined similarly.
%  have similar property as $F_{i}\subset A_{i}$.
Let $\varepsilon_{1}<\varepsilon_{2}<\cdots$ with $\sum\varepsilon_{n}<+\infty$.
Now we will apply Theorem 5.10 (the main existence theorem), Proposition 5.9 (to pull back the map of invariants), Theorem 4.13 (to get the desired decomposition with some sets being weakly approximately constant) and Theorem 6.19 (the main uniqueness theorem) to prove  our main theorem.

Let $k_{1}=1$. Apply Theorem 4.13 to $A_{1}$, there are  an $m'>k_{1}=1$ and projections $Q_{0},Q_{1}\in A_{m'}$ such that $Q_{0}+Q_{1}=\phi_{1,m'}(\textbf{1}_{A_{i}})$, and there are
 unital maps $$\psi_{0}\in Map(A_{1},Q_{0}A_{m'}Q_{0})_{1},~~~~\mbox{and}~~~~\psi_{1}\in Hom(A_{1},Q_{1}A_{m'}Q_{1})_{1}$$
with the following properties:\\
%such that the following statement are true\\
(1) $\|\phi_{1,m'}(f)-\psi_{0}(f)\oplus\psi_{1}(f)\|<\varepsilon_1$ ~~for all $ f\in F$; \\
(2) $\omega(\psi_{0}(F_{1}))<\varepsilon_1;$\\
(3) The homomorphism $\psi_{1}$ factors through $C$ --- a direct sum of matrix algebra over $C[a,b]$ (where $[a,b]\subset[0,1]$) or $\mathbb{C}$
as $$\psi_{1}:A_{1}\xrightarrow{\xi_{1}}C\xrightarrow{\xi_{2}}Q_{1}A_{m'}Q_{1},$$
where $\xi_{1}$ and $\xi_{2}$ are unital homomorphisms.

Set
%\begin{center}
$I=\{i~|~~\xi_{2}^{i,j}\neq0$ for $j$ with $A_{m'}^{j}$ of type $PM_{\bullet}((T_{\uppercase\expandafter{\romannumeral2},k})P$ or $M_{l}(I_{k}))\}$
%\end{center}
and set\\
%\begin{center}
$J=\{j~|~~A_{m'}^{j}$ is of type $PM_{\bullet}C(T_{\uppercase\expandafter{\romannumeral2},k})P$ or $M_{l}((I_{k})\}$.
%\end{center}
 Let
$$\widetilde{A}_{1}=\bigoplus\limits_{i\in I}C^{i}\oplus\bigoplus\limits_{j\in J}Q_{0}^{j}A_{m'}^{j}Q_{0}^{j}\oplus\bigoplus\limits_{j\notin J}A_{m'}^{j},~\mbox{with subset}~
\widetilde{F}_{1}=\bigoplus\limits_{j\in J}\psi_{0}^{-,j}(F_{1})\oplus\bigoplus\limits_{i\in I}\xi_{1}^{i}(F_{1})\oplus\bigoplus\limits_{j\notin J}\phi_{\textbf{1}, m'}^{-j}(F_{1})\subset\widetilde{A}_{1}.$$
(See 1.3 for the notation.)

Let $\xi: \widetilde{A}_{1}\longrightarrow A_{m'}$ be the direct sum of the following partial maps: the inclusion from
$Q_{0}^{j}A_{m'}^{j}Q_{0}^{j}$ to $A_{m'}^{j}$ (for $j\in J$), the  identity map from $A_{m'}^{j}$ to $A_{m'}^{j}$ (for $j\notin J$),  and
$\xi_{2}^{i,j}: C^{i}\longrightarrow Q_{1}^{j}A_{m'}^{j}Q_{1}^{j}$ (for $i\in I$ and $j\in J$).

Since we will reserve the notation $\xi$ for the positive linear map for $AffT$ spaces, denote the above map $\xi: \widetilde{A}_{1}\longrightarrow A_{m'}$ by $\varphi$. Let $\widetilde{\varphi}: A_{1}\longrightarrow\widetilde{A}_{1}$ be given by:

if $\widetilde{A}_{1}^{k}=Q_{0}^{j}A_{m'}^{j}Q_{0}^{j}$ for $j\in J$, then $\widetilde{\varphi}^{i,k}: A_{1}^{i}\longrightarrow\widetilde{A}_{1}^{k}$
is defined by $\widetilde{\varphi}^{i,k}=\psi_{0}^{i,k}$;

if $\widetilde{A}_{1}^{k}=A_{m'}^{j}$ for $j\notin J$, then $\widetilde{\varphi}^{i,k}: A_{1}^{i}\longrightarrow\widetilde{A}_{1}^{k}$
is defined by $\widetilde{\varphi}^{i,k}=\phi_{1,m'}^{i,j}$;

if $\widetilde{A}_{1}^{k}=C^{j}$ for $j\in I$, then $\widetilde{\varphi}^{i,k}: A_{1}^{i}\longrightarrow\widetilde{A}_{1}^{k}$
is defined by $\widetilde{\varphi}^{i,k}=\xi_{1}^{i,j}: A_{1}^{i}\longrightarrow C^{j}$.

Then $\|\phi_{1,m'}(f)-\varphi\circ\widetilde{\varphi}(f)\|<\varepsilon_{1}$ ~~for all $f\in F.$

Now the algebra $\widetilde{A}_{1}$ with $\widetilde{F}_{1}\subset\widetilde{A}_{1}$ and $\varphi: \widetilde{A}_{1}\longrightarrow A_{m'}$ satisfies the condition of Theorem 6.19 (in place of  $B$ with $\xi: B\longrightarrow A_{n}$). By Theorem 6.19, there is $m>m'$, $H=\oplus H^{i}\subseteq\oplus AffTA_{m}^{i}$, $E=\oplus E^{i}\subseteq
\bigoplus\limits_{i}U(A_{m}^{i})/\widetilde{SU(A_{m}^{i})}$, and $\sigma_{1}>0$ as desired in the theorem.

Apply Theorem 5.10 to $H=\oplus H^{i}$ and $E=\oplus E^{i}$ (of course with $m$ in place of $n$), and $\frac{\sigma_{1}}{8}$ (in places of both $\varepsilon_{1}$ and $\varepsilon_{2}$), there is $m''>m$ as in the theorem (in place of $m$).

\noindent\textbf{7.5.}~~For each $A_{m''}^{i_{0}}$, and each $s\in \mathbb{Z}_{+}$, let $P_{s}=\phi_{m'',m''+s}(\textbf{1}_{A_{m''}^{i_{0}}})\in A_{m''+s}$, with
$P_{0}=\textbf{1}_{A_{m''}^{i_{0}}}.$ Let $\nu_{m''}: A_{m''}\longrightarrow B_{m''}$ be as in 7.3,  let $Q_{0}=\nu_{m''}(\textbf{1}_{A_{m''}^{i_{0}}})$, and $Q_{s}=\psi_{m'',m''+s}(Q_{0})\in B_{m''+s}$. Let $P=\phi_{m'',\infty}(\textbf{1}_{A_{m''}^{i_{0}}})\in A, Q=\psi_{m'',\infty}(Q_{0})\in B$.

(Note that the projections $Q_0, Q_s$ here are different from the projections $Q_0, Q_1$ in 8.4, which appeared in the decomposition of a homomorphism. We believe that it will not cause any confusion.)

Then $\alpha\in KK(A,B)$ induces $\alpha_{[PAP,QBQ]}\in KK(PAP,QBQ)$.  Consider two inductive limits
$$C\triangleq PAP=lim(C_{s}\triangleq P_{s}A_{m''+s}P_{s},\widetilde{\phi}_{s,t}=\phi_{m''+s,m''+t}|_{C_{s}}),~~~~\mbox{and}$$
$$D\triangleq QBQ=lim(D_{s}\triangleq Q_{s}B_{m''+s}Q_{s},\widetilde{\psi}_{s,t}=\psi_{m''+s,m''+t}|_{D_{s}}),$$ as in Proposition 5.9, with
$\alpha_{(C,D)}, \xi^{P,Q}, \gamma^{P,Q}$ as in 5.8, and with
$$H'\subseteq AffTC_{0}=AffTA_{m''}^{i_{0}}=C_{\mathbb{R}}(SpA_{m''}^{i_{0}})~~\mbox{and}~~E'\subseteq U(C_0)/\widetilde{SU(C_0)} =U(A_{m''}^{i_{0}})/\widetilde{SU(A_{m''}^{i_{0}})}$$ to be
$$H'=\bigcup_{i} \imath_{T,i}(AffT\phi_{m,m''}^{i,i_{0}}(H^{i}))\subset AffTA_{m''}^{i_{0}}~~\mbox{and}~~E'=\bigcup_{i}\imath_{*,i}((\phi_{m,m''}^{i,i_{0}})^{\natural}(E^{i}))\subset U(A_{m''}^{i_{0}})/\widetilde{SU(A_{m''}^{i_{0}})}, $$
respectively, where
$\imath_{T,i}: AffT(\phi_{m,m''}^{i,i_{0}}({\bf 1}_{A_m^{i}}) A_{m''}^{i_{0}}\phi_{m,m''}^{i,i_{0}}({\bf 1}_{A_m^{i}})) \to Aff TA_{m''}^{i_{0}}~~~~~~\mbox{and}~~~~~~~~$
$$\imath_{*,i}: U (\phi_{m,m''}^{i,i_{0}}({\bf 1}_{A_m^{i}}) A_{m''}^{i_{0}}\phi_{m,m''}^{i,i_{0}}({\bf 1}_{A_m^{i}}))/\widetilde{SU}(\phi_{m,m''}^{i,i_{0}}({\bf 1}_{A_m^{i}}) A_{m''}^{i_{0}}\phi_{m,m''}^{i,i_{0}}({\bf 1}_{A_m^{i}})) \to U(A_{m''}^{i_{0}})/\widetilde{SU(A_{m''}^{i_{0}})},$$
are induced by the inclusion maps (see 2.4 and 2.8).
Let $\varepsilon=
\frac{\sigma_{1}}{8rank(\textbf{1}_{A_{m''}^{i_{0}}})}.$ By Proposition 5.9, there is an  $s>0,$ there is a KK-element $\alpha_{1}\in KK(C_{0},D_{s})$ which can be realized by a homomorphism
(in fact, $\alpha_{1}$ can be choose to be $KK(\psi_{m'',m''+s}\circ\nu_{m''}|_{A_{m''}^{i_0}})$), there is a positive contraction   $$\xi_{1}: AffTC_{0}\longrightarrow AffTD_{s}$$
compatible with the $K_{0}$ map induced by $\alpha_{1}$ (see (c) of 5.6), and there is a contractive group homomorphism
$$\gamma_{1}: U(C_{0})/\widetilde{SU(C_{0})}\longrightarrow U(D_{s})/\widetilde{SU(D_{s})}$$ compatible with $\xi_{1}$
(see (e) of 5.6) and compatible with the $K_{1}$ map induced by $\alpha_{1}$ (see (d) of 5.6), such that\\
(1)  $\alpha_{1}\times KK(\widetilde{\psi}_{s,\infty})=KK(\widetilde{\phi}_{0,\infty})\times\alpha_{(C,D)},$\\
(2)  $\|AffT\widetilde{\psi}_{s,\infty}\circ\xi_{1}(h)-\xi^{P,Q}\circ AffT\widetilde{\phi}_{0,\infty}(h)\|<\frac{\sigma_{1}}{8 rank(\textbf{1}_{A_{m''}^{i_{0}}})} ~~~~\forall h\in H'$\\
(3)  $dist(\widetilde{\psi}^{\natural}_{s,\infty}\circ \gamma_{1}(g),\gamma^{P,Q}\circ\widetilde{\phi}^{\natural}_{0,\infty}(g))<\frac{\sigma_{1}}{8rank(\textbf{1}_{A_{m''}^{i_{0}}})}~~~~\forall g\in E^{\prime}.$

\noindent\textbf{Lemma 7.6.}  Let $A,B,C$ be basic $\mathcal{HD}$ building blocks, let  $\phi: A\longrightarrow B$ be a (not necessarily unital) homomorphism,
and $F\subseteq AffTA$  a finite set. Let $F'_{1}=AffT\phi(F)\subseteq AffT\phi(1_{A})B\phi(1_{A})$ which can be regarded as subset $F_1=\imath_{*}(F'_{1})$ of $AffTB$ (see 2.8).
%(where $\imath$ is the inclusion $\phi(1_{A})B\phi(1_{A})\xrightarrow{\imath}B$).
Suppose that  $\alpha\in KK(B,C)$ can be realized by a homomorphism. Let $\xi_{1},\xi_{2}: AffTB\longrightarrow AffT\alpha(\textbf{1}_{B})C\alpha(\textbf{1}_{B}))$ be two unital positive linear maps. Then
% \begin{center}
 $$\max\limits_{f\in F}\|(\xi_{1}\circ AffT\phi)_{([\textbf{1}_{A}], \alpha[\phi(\textbf{1}_{A})])}(f)-(\xi_{2}\circ AffT\phi)_{([\textbf{1}_{A}], \alpha[\phi(\textbf{1}_{A})])}(f)\|
=(\max\limits_{f\in F_{1}}\|\xi_{1}(f)-\xi_{2}(f)\|)\cdot\frac{rank(\textbf{1}_{B})}{rank(\phi(\textbf{1}_{A}))}. $$
%\end{center}
In particular,  if $\|\xi_{1}(f)-\xi_{2}(f)\|<\frac{\sigma_{1}}{8rank(\textbf{1}_{B})}$ ~for all $f\in F_{1}$ then
$$\|(\xi_{1}\circ AffT\phi)_{([\textbf{1}_{A}],\alpha[\phi(\textbf{1}_{A})])}(f)-(\xi_{2}\circ AffT\phi)_{([\textbf{1}_{A}],\alpha[\phi(\textbf{1}_{A})])}(f)\|<\frac{\sigma_{1}}{8}~~~~\forall f\in F.$$
\begin{proof}
Note that the map  $$\imath_*:~~AffT(\phi(\textbf{1}_{A})B\phi(\textbf{1}_{A}))=C_{\mathbb{R}}(SpB)\longrightarrow AffTB=C_{\mathbb{R}}(SpB)$$
induced by inclusion, is given by $g\mapsto\frac{rank\phi(\textbf{1}_{A})}{rank(\textbf{1}_{B})}\cdot g$. The lemma is obviously true.\\
\end{proof}

\noindent\textbf{Lemma 7.7.}~~ Let $A,B,C$ be basic building blocks. Let $\phi: A\longrightarrow B$ be a (not necessarily  unital) homomorphism, let  $\alpha\in KK(B,C)$ be induced by a homomorphism, and let $R\in C$ satisfy $\alpha[\textbf{1}_{B}]=[R]\in K_{0}(C)$. Let $E\subset U(A)/\widetilde{SU(A)}$ be a finite set and $$E_{1}=\phi^{\natural}(E)\subseteq U(\phi(\textbf{1}_{A})B\phi(\textbf{1}_{A}))/\widetilde{SU}(\phi(\textbf{1}_{A})B\phi(\textbf{1}_{A}))\subseteq U(B)/\widetilde{SU(B)}$$ (by regarding
$u\in U(\phi(\textbf{1}_{A})B\phi(\textbf{1}_{A}))$ as $u\oplus(\textbf{1}_{B}-\phi(\textbf{1}_{A}))\in U(B)).$ Let
 $$\gamma_{1},\gamma_{2}: U(B)/\widetilde{SU(B)}\longrightarrow
U(RCR)/\widetilde{SU}(RCR)$$
be two contractive group homomorphisms compatible with $\alpha_{*}: K_{1}(B)\longrightarrow K_{1}(C)$. Then
$$\max\limits_{g\in E}dist((\gamma_{1}\circ\phi^{\natural})|_{([\textbf{1}_{A}],\alpha[\phi(\textbf{1}_{A})])}(g),(\gamma_{2}\circ\phi^{\natural})|_{([\textbf{1}_{A}],\alpha[\phi(\textbf{1}_{A})])}(g))\leq
\frac{rank\phi(\textbf{1}_{A})}{rank(\textbf{1}_B)}\max\limits_{g\in E_{1}}dist(\gamma_{1}(g),\gamma_{2}(g)).$$
\begin{proof}
By Lemma 7.2, this lemma is also obvious.\\
\end{proof}

\noindent\textbf{7.8.} Applying  7.5 to each block $A_{m''}^{i}$ of $A_{m''}$, we can find a positive integer $s$ working for all $i$.
 Let us denote   $\psi_{m'',m''+s}\circ \nu_{m''}: A_{m''}\longrightarrow B_{m''+s}$ by $\chi$, and denote $m''+s$ by $l_{1}$. Then from 7.5, there is
$$\xi_{1}^{i}: AffTA_{m''}^{i}\longrightarrow AffT(\chi(\textbf{1}_{A_{m''}^{i}})B_{l_{1}}\chi(\textbf{1}_{A_{m''}^{i}}))$$
which is compatible with $K_{0}$ map induced by $\chi$, and there is
$$\gamma_{1}^{i}: U(A_{m''}^{i})/\widetilde{SU}(A_{m''}^{i})\longrightarrow U(\chi(\textbf{1}_{A_{m''}^{i}})B_{l_{1}}\chi(\textbf{1}_{A_{m''}^{i}})/\widetilde{SU}(\chi(\textbf{1}_{A_{m''}^{i}})B_{l_{1}}\chi(\textbf{1}_{A_{m''}^{i}}).$$
which is compatible with $\xi_{1}^{i}$ and $K_{1}$ map induced by $\chi$, with condition (1), (2), (3) in 7.5.

\noindent\textbf{Lemma 7.9.}~~ Let $A=\oplus A^{i}\in\mathcal{HD}$, $B=\oplus B^{j}\in\mathcal{HD}$, and let  $\chi: A\longrightarrow B$ be a homomorphism. Suppose that for each $i$, there are a unital positive map
$$\xi^{i}: AffTA^{i}\longrightarrow AffT\chi(\textbf{1}_{A^{i}})B\chi(\textbf{1}_{A^{i}}),$$
and a contractive group homomorphism
$$\gamma^{i}: U(A^{i})/\widetilde{SU(A^{i})}\longrightarrow U(\chi(\textbf{1}_{A^{i}})B\chi(\textbf{1}_{A^{i}}))/\widetilde{SU}(\chi(\textbf{1}_{A^{i}})B\chi(\textbf{1}_{A^{i}}))$$
which are compatible with each other and compatible with the $K$-theory map induced by $\chi$. Then there are a unique unital positive linear map
$$\xi: AffTA\longrightarrow AffT\chi(\textbf{1}_{A})B\chi(\textbf{1}_{A})$$
which is compatible with the  $K_{0}$ map induced by $\chi$, a unique contractive group homomorphism $$\gamma: U (A)/\widetilde{SU(A)}\longrightarrow U(\chi(\textbf{1}_{A})B\chi(\textbf{1}_{A}))/\widetilde{SU}(\chi(\textbf{1}_{A})B\chi(\textbf{1}_{A}))$$ which is compatible with $\xi$ and $K_{1} $ map induced by $\chi$, such that for each $i$
$$\xi|_{(\textbf{1}_{A^{i}},\chi(\textbf{1}_{A^{i}}))}=\xi^{i} ~~~~~~~\mbox{and}~~~~\gamma|_{(\textbf{1}_{A^{i}},\chi(\textbf{1}_{A^{i}}))}=\gamma^{i}.$$
\begin{proof} For $f=(f_1,f_2,\cdots, f_{\bullet})\in AffTA=\oplus_i AffTA^i$ and \\$g=(g_1,g_2,\cdots, g_{\bullet})\in U(A)/\widetilde{SU(A)}=\oplus_i U(A^i)/\widetilde{SU(A^i)}$, define
$$\xi(f)=\sum_i \imath_{T,i}(\xi^i(f_i))~~~~~~~ \mbox{and}~~~~~~ \gamma (g)=\sum_i\imath_{*,i}(\gamma^i(g_i)), $$
where $ \imath_{T,i}: AffT\chi(\textbf{1}_{A^{i}})B\chi(\textbf{1}_{A^{i}}) \to AffT\chi(\textbf{1}_{A})B\chi(\textbf{1}_{A})$
and \\$\imath_{*,i}:
U(\chi(\textbf{1}_{A^{i}})B\chi(\textbf{1}_{A^{i}}))/\widetilde{SU}(\chi(\textbf{1}_{A^{i}})B\chi(\textbf{1}_{A^{i}}))
\to U(\chi(\textbf{1}_{A})B\chi(\textbf{1}_{A}))/\widetilde{SU}(\chi(\textbf{1}_{A})B\chi(\textbf{1}_{A}))$ are induced by the inclusions. Evidently, $\xi$ and $\gamma$ are as desired.  \\
\end{proof}

\noindent\textbf{7.10.} From 7.8,  and applying  7.9, we get  a $KK$ element $\alpha_{1}\in KK(A_{m''},\chi(\textbf{1}_{A_{m''}})B_{l_{1}}\chi(\textbf{1}_{A_{m''}}))$, a  unital positive linear map
$\xi_{1}: AffTA_{m''}\longrightarrow AffT(\chi(\textbf{1}_{A_{m''}})B_{l_{1}}\chi(\textbf{1}_{A_{m''}}))$, and a contractive group homomorphism
$$\gamma_{1}: U(A_{m''})/\widetilde{SU}(A_{m''})\longrightarrow U(\chi(\textbf{1}_{A_{m''}})B_{l_{1}}\chi(\textbf{1}_{A_{m''}})/\widetilde{SU}(\chi(\textbf{1}_{A_{m''}})B_{l_{1}}\chi(\textbf{1}_{A_{m''}})),$$ which are compatible with
each other and satisfying  the following conditions:\\
(1) $\alpha_{1}\times KK(\psi_{l_{1},\infty})=KK(\psi_{m'',\infty})\times\alpha\in KK(A_{m''},B);$\\
~~~~\\
(2) $\|AffT\psi_{l_{1},\infty}\circ\xi_{1}|_{(\textbf{1}_{A_{m''}^{i_{0}}},\psi_{l_{1},\infty}[\chi(\textbf{1}_{A_{m''}^{i_{0}}})])}(h)\\~~~~~~~~~~~-\xi^{[\phi_{m'',\infty}(\textbf{1}_{A_{m''}^{i_{0}}}),\alpha(\phi_{m'',\infty}(\textbf{1}_{A_{m''}^{i_{0}}}))]}
\circ AffT\phi_{m'',\infty}|_{(\textbf{1}_{A_{m''}^{i_{0}}}, \phi_{m''\infty}(\textbf{1}_{A_{m''}^{i_{0}}}))}(h)\|
<\frac{\sigma_{1}}{8rank(\textbf{1}_{A_{m''}^{i_{0}}})},\\~~~~~~~~~~~~~~~\forall h\in H^{\prime}\subset AffT(A_{m''}^{i_{0}})$ (as described 7.5);\\
(3) $dist(((\psi_{l_{1},\infty})^{\natural}\circ\gamma_{1})_{(\textbf{1}_{A_{m''}^{i_{0}}},\psi_{l_{1},\infty}\chi[\textbf{1}_{A_{m''}^{i_{0}}}])}(g),\\
~~~~~~~~~~~~~~~~~\gamma^{[\phi_{m'',\infty}(\textbf{1}_{A_{m''}^{i_{0}}}),\alpha(\phi_{m'',\infty}(\textbf{1}_{A_{m''}^{i_{0}}}))]}\circ (\phi_{m'',\infty})^{\natural}|_{(\textbf{1}_{A_{m''}^{i_{0}}},\psi_{l_{1},\infty}\chi[\textbf{1}_{A_{m''}^{i_{0}}}])}(g))<\frac{\sigma_{1}}{8rank(\textbf{1}_{A_{m''}^{i_{0}}})},\\~~~~~~~~ ~~~~~~~~~\forall
g\in E^{\prime}\subset U(A_{m''}^{i_{0}})/\widetilde{SU(A_{m''}^{i_{0}})}$(as described in 7.5).

Applying  Theorem 5.10 and notice the choice of $m''$, and $\varepsilon_{1}=\varepsilon_{2}=\frac{\sigma_{1}}{8}$, there is a homomorphism
$$\wedge': A_{m}\longrightarrow B_{l_{1}},$$ such that\\
(4) $KK(\wedge')=KK(\phi_{m,m''})\times\alpha_{1}\in KK(A_{m},B_{l_{1}}).$\\
(5) $\|(AffT\wedge')|_{\textbf{1}_{A_{m}}^{i},\wedge(\textbf{1}_{A_{m}}^{i})}(h)-(\xi_{1}\circ AffT\phi_{m,m''})|_{\textbf{1}_{A_{m}}^{i},\wedge(\textbf{1}_{A_{m}}^{i})}(h)\|<\frac{\sigma_{1}}{8} ~~\forall h\in H^i$\\
(6) $dist(\wedge^{'\natural}|_{(\textbf{1}_{A_{m}}^{i},\wedge(\textbf{1}_{A_{m}}^{i}))}(g), (\gamma_{1}\circ\phi_{m,m''}^{\natural})_{(\textbf{1}_{A_{m}}^{i},\wedge(\textbf{1}_{A_{m}}^{i}))}(g))<\frac{\sigma_{1}}{8} ~~~~~~~\forall g\in E^{i}.$
\\Define $\wedge_{1}:~~A_{k_{1}}=A_{1}\longrightarrow B_{l_{1}}$ by
$\wedge_{1}=\wedge'\circ\phi_{1,m}=\wedge'\circ\phi_{m',m}\circ\phi_{1,m'}$. Note that\\
(7) $\|\wedge_{1}(f)-\wedge'\circ\phi_{m',m}\circ\varphi\circ\widetilde{\varphi}(f)\|<\varepsilon, ~~~~~\forall f\in F_{1}\subset A_{1},$
\\where $\widetilde{\varphi}: A_{1}\longrightarrow \widetilde{A}_{1}$ and $\varphi: \widetilde{A}_{1}\longrightarrow A_{m'}$ are as in 7.4. Thus  we get the first
map $\wedge_{1}: A_{k_{1}}\longrightarrow B_{l_{1}},$ which also satisfies\\
(8) $KK(\wedge_{1})=KK(\phi_{k_{1},l_{1}})\times KK(\nu_{l_{1}})=KK(\nu_{k_{1}})\times KK(\psi_{k_{1},l_{1}}).$

\noindent\textbf{7.11.} From 7.10, combining (1) and (4), we obtain\\
(a) $KK(\psi_{l_{1},\infty}\circ\wedge')=KK(\phi_{m,\infty})\times\alpha.$

Combining (2) and (5) and applying 7.6,  we obtain\\
(b)\! $\|AffT(\psi_{l_{1},\infty}\circ\wedge')|_{(\textbf{1}_{A_{m}}^{i},\psi_{l_{1},\infty}\circ\wedge'(\textbf{1}_{A_{m}}^{i}))}(h)-\xi^{\phi_{m,\infty}(\textbf{1}_{A_{m}}^{i}),\alpha(\phi_{m,\infty}(\textbf{1}_{A_{m}}^{i}))}\circ
AffT(\phi_{m,\infty}|_{A_{m}^{i}})(h)\|<\frac{\sigma_{1}}{4}~~\forall h\in H^{i}.$

 Combining (3) and (6) and applying 7.7, we obtain\\
(c) $dist((\psi_{l_{1},\infty}\circ\wedge')^{\natural}|_{(\textbf{1}_{A_{m}}^{i},\psi_{l_{1},\infty}\circ\wedge'(\textbf{1}_{A_{m}}^{i}))}(g),
\gamma^{\phi_{m,\infty}(\textbf{1}_{A_{m}}^{i}),\alpha\circ\phi_{m,\infty}(\textbf{1}_{A_{m}}^{i})}\circ(\phi_{m,\infty}|_{A_{m}^{i}})^{\natural}(g))<\frac{\sigma_{1}}{4}~~\forall g\in E^{i}.$

\noindent\textbf{7.12.} Go the above procedure 7.4---7.11  for $B_{l_{1}}$ (in place of $A_{1}$) $$G_{1}'=G_{l_{1}}\cup\wedge_{1}(F_{1})\cup\wedge'\circ\phi_{m',m}\circ\widetilde{\varphi}(\widetilde{F}_{1})\cup\wedge'(F_{m})\subset B_{l_{1}}$$ (in place of $F_{1}$) and $\varepsilon$. Applying Theorem 4.13, we obtain $B_{l_{1}}\xrightarrow{\widetilde{\varphi'}}\widetilde{B}_{1}\xrightarrow{\varphi'}B_{n'}$ with $\widetilde{G}_{1}=\oplus\widetilde{G}_{1}^{i}\subset\widetilde{B}_{1}=\oplus\widetilde{B}_{1}^{i}$ such that $\widetilde{G}_{1}\supset\widetilde{\varphi'}(G_{1}')$,  such that $\omega(\widetilde{G}_{1}^{i})<\varepsilon_{1}$ if $\widetilde{B}_{1}^{i}$ is of form $PM(C(T_{\uppercase\expandafter{\romannumeral2},k}))P$ or $M_{l}(I_{k})$, and such that
%\begin{center}
%(**)~~~~~~~~~~~~~~~~~~~~~~~~~~~~~~~~~~~
%$\|\varphi'\circ\widetilde{\varphi'}(g)-\psi_{l_{1},n'}(g)\|<\varepsilon_{1}, \forall g\in G_{1}^{\prime}.$
%\end{center}
\begin{equation*}
(**)~~~~~~~~~~~~~~~~~~~~~~~~~~~
\|\varphi'\circ\widetilde{\varphi'}(g)-\psi_{l_{1},n'}(g)\|<\varepsilon_{1},~~~~~~ \forall g\in G_{1}^{\prime}.
\end{equation*}
We will repeat the procedure from  7.4 to 7.11, briefly. Applying Theorem 6.19 to $\varphi': \widetilde{B}_{1}\longrightarrow B_{n'}$ and $\varepsilon_{1}$, we obtain $n>n'$, $\sigma_{2}>0$, $$\widetilde{H}=\oplus\widetilde{H}^{i}\subseteq AffTB_{n}=\bigoplus\limits_{i}AffTB_{n}^{i},~~\mbox{and}~~\widetilde{E}=\oplus\widetilde{E}^{i}\subseteq U(B_{n})/\widetilde{SU(B_{n})}=\bigoplus\limits_{i}U(B_{n}^{i})/\widetilde{SU(B_{n}^{i})}$$
as desired in Theorem 6.19. Furthemore, we can assume
$$\widetilde{H}\supseteq AffT(\psi_{l_{1},n}\circ\wedge')(H)~~~~\mbox{and}~~~~\widetilde{E}\supset(\psi_{l_{1},n}\circ\wedge')^{\natural}(E),$$
where $H\subseteq AffTA_{m}$ and $E\subseteq U(A_{m})/\widetilde{SU(A_{m})}$ are as in 7.4.  (Here, for any nonzero projection $P\in B$, we regard $AffT(PBP)$ and $U(PBP)/\widetilde{SU}(PBP)$ as the subsets of $AffTB$ and $U(B)/\widetilde{SU}(B)$ by the injective
maps $\imath_T$ and $\imath_*$ in 2.4 and 2.8, respectively.)

Let $M=\max\limits_{i}rank(\textbf{1}_{B_{n}^{i}}).$ Going the procedure  from 7.4 to 7.11, and applying Theorem 5.10 and Proposition 5.9, we can obtain $k'_{2}>n$,
$\mu': B_{n}\longrightarrow A_{k'_{2}}$ such that\\
(a$'$) $KK(\phi_{k'_{2}, \infty}\circ\mu')=KK(\psi_{n,\infty})\times\alpha^{-1};$ \\
(b$'$)$\|AffT(\phi_{k'_{2},\infty}\circ\mu')|_{(\textbf{1}_{B_{n}}^{i},\phi_{k'_{2},\infty}\circ\mu'(\textbf{1}_{B_{n}^{i}}))}(h)\\~~~~~~~~~~~~~~~~ -(\xi^{{\alpha^{-1}(\psi_{n,\infty}(\textbf{1}_{B_{n}}^{i})),\psi_{n,\infty}(\textbf{1}_{B_{n}}^{i})}})^{-1}\circ
AffT(\psi_{n,\infty})|_{B_{n}^{i}}(h)\|<\frac{min(\sigma_{2},\frac{\sigma_{1}}{M})}{4},~~~~~~ \forall h\in \widetilde{H}^i;$  \\
(c$'$) $dist(((\phi_{k'_{2},\infty}\circ \mu')^{\natural})|_{(\textbf{1}_{B_{n}^{i}},\phi_{k'_{2},\infty}\circ\mu'(\textbf{1}_{B_{n}^{i}}))}(g),\\~~~~~~~~~~~~~~~~~~~~~~~~~~~~~
(\gamma^{\alpha^{-1}(\psi_{n,\infty}(\textbf{1}_{B_{n}}^{i}),\psi_{n,\infty}(\textbf{1}_{B_{n}}^{i}))})^{-1}\circ(\psi_{n,\infty}|_{B_{n}^{i}})^{\natural}(g))<\frac{min(\sigma_{2},\frac{\sigma_{1}}{M})}{4}, ~~~~~~~\forall g\in\widetilde{E}^i.$

\noindent\textbf{7.13.}~~Applying  7.2 and  7.6 (see 7.7 and 7.9 also) and notice  that $M=\max\limits_{i}rank(\textbf{1}_{B_{n}^{i}}),$ we get the following estimations.
(Let $P_{i}=\phi_{m,k'_{2}}(\textbf{1}_{A_{m}^{i}})$ which is equivalent to $(\mu' \circ\psi_{l_{1},n}\circ\wedge')(\textbf{1}_{A_{m}^{i}})$.)\\
Combining (a) and (a$^{'}$), we get\\
($a''$) $KK(\wedge')\times KK(\psi_{l_{1},n})\times KK(\mu')\times KK(\phi_{k'_2,\infty})=KK(\phi_{m,k'_{2}})\times KK(\phi_{k'_2,\infty}).$\\
Combining (b) and ($b^{'}$), we get\\
($b''$) $\|AffT\phi_{k'_{2},\infty}|_{P_{i}A_{k_{i}}P_{i}}\Big{\{}AffT\mu'\circ\psi_{l_{1},n}\circ\wedge'|_{(\textbf{1}_{A_{m}^{i}},P_{i})}(h)
-AffT\phi_{m,k'_{2}}|_{(\textbf{1}_{A_{m}^{i}}, P_{i})}(h)\Big{\}}\|<\frac{\sigma_{1}}{2},~~ \forall h\in H^{i}.$\\
Combining (c) and ($c'$), we get\\
($c''$)
$dist((\phi_{k'_{2},\infty}|_{P_{i}A_{k_{i}}P_{i}})^{\natural}(\mu'\circ\psi_{l_{1},n}\circ\wedge')^{\natural}|_{(\textbf{1}_{A_{m}^{i}},P_{i})}(g),
(\phi_{k'_{2},\infty}|_{P_{i}A_{k_{2}'}P_{i}})^{\natural}(\phi_{m,k'_{2}})^{\natural}|_{(\textbf{1}_{A_{m}^{i}},P_{i})}(g))<\frac{\sigma_{1}}{2},~~ \forall g\in E^{i}.$

Evidently,  passing  to $k_{2}>k'_{2}$,
%using a routine  argument about inductive limit procedure,
and replacing
 $\mu'$ by $\mathcal{M}'=\phi_{k'_{2},k_{2}}\circ \mu'$, we get (now  denote $\phi_{m,k_{2}}(\textbf{1}_{A_{m}^{i}})$ by $P_{i}$,  which is equivalent to $\mathcal{M}'\circ\psi_{l_{1},n}\circ\wedge'(\textbf{1}_{A_{m}^{i}})$)\\
$(A)$ $KK(\wedge')\times KK(\psi_{l_{1},n})\times KK (\mathcal{M}')=KK(\phi_{m,k_{2}}).$\\
$(B)$ $\|AffT(\mathcal{M}' \circ\psi_{l_{1},n}\circ\wedge')|_{(\textbf{1}_{A_{m}^{i}},P_{i})}(h)-
AffT\phi_{m,k_{2}}|_{(\textbf{1}_{A_{m}^{i}},P_{i})}(h)\|<\sigma_{1},~~~~\forall h\in H^{i}.$\\
$(C)$ $dist\big((\mathcal{M}' \circ\psi_{l_{1},n}\circ\wedge')^{\natural}|_{(\textbf{1}_{A_{m}^{i}},P_{i})}(g),
(\phi_{m,k_{2}})^{\natural}|_{(\textbf{1}_{A_{m}^{i}},P_{i})}(g)\big)<\sigma_{1},~~~~ \forall g\in E^i.$

\noindent\textbf{7.14.}~~Since $m,\sigma_{1},H,E$ are chosen for $\widetilde{F}_{1}\subset\widetilde{A}_{1}$ and $\varphi: \widetilde{A}_{1}\longrightarrow A_{m'}$ as in
Theorem 6.19,  there is a unitary $u\in A_{k_{2}}$ such that
%\begin{center}
$\|Adu\circ\mathcal{M}'\circ\psi_{l_{1},n}\circ\wedge'\circ\phi_{m',m}\circ\varphi(f)-\phi_{m',k_{2}}\varphi(f)\|<7\varepsilon, ~~\forall\;f\in\widetilde{F}_{1}.$\\
%\end{center}
Combining the above with (7) of 7.10, we get
%\begin{center}
$$\parallel Adu\circ\mathcal{M}'\circ\psi_{l_{1},n}\circ\wedge_{1}(f)-\phi_{k_{1}k_{2}}(f)\|<8\varepsilon_{1},~~~~ \forall f\in F_{1}\subset A_{1}=A_{k_1}.$$
%\end{center}
Let $\mathcal{M}_{1}=Adu\circ\mathcal{M}'\circ\psi_{l_{1},n}: B_{l_{1}}\longrightarrow A_{k_{2}}$, we get the following almost commuting diagram
\vspace{-4mm}
\begin{center}$
\xymatrix{
  A_{1} \ar[rr]^{\phi_{k_{1}k_{2}}} \ar[dr]^{~~(F_{1},8\varepsilon)}_{\wedge_{1}}
                &  &   A_{k_2}   \\
                  & B_{l_1}    \ar[ru]_{\mathcal{M}_{1}}           }
$
\end{center}
\vspace{-4mm}
That is,  $\|\mathcal{M}_{1}\circ\wedge_{1}(f)-\phi_{k_{1},k_{2}}(f)\|<8\varepsilon,~~~~~~ \forall f\in F_{1}.$
%\begin{figure*}[!t]
%\centering
%\includegraphics[width=1.5in]{fig1.eps}
%\end{figure*}

\noindent\textbf{7.15.}~~By the same procedure, we can construct the next almost commuting diagram. We need to use $\widetilde{B}_{1}$, the factorization $B_{l_{1}}\xrightarrow{\widetilde{\varphi'}}\widetilde{B}_{1}\xrightarrow{\varphi'}B_{n'}$, $(**)$ from 7.12,  and finally,  the choice of $n>n'$, to get almost commutative  diagram
$$
\xymatrix{
                &  A_{k_2} \ar[dr]^{\wedge_{2}}             \\
 B_{l_1} \ar[ur]^{\mathcal{M}_{1}}_{~~(G'_{1},8\varepsilon_{1})} \ar[rr] & &     B_{l_2}        }
$$
When we work on $\wedge_{2}$ we should use $\varepsilon_{2}$ (in place of $\varepsilon_{1}$) in order to make the next two diagrams almost commutes up to $8\varepsilon_{2}$ on the given sets (e.g., $F_{k_2}$ and $G_{l_2}$). Note that $G'_{1}\supset G_{l_{1}}(\subset B_{l_{1}})$ and
$G'_{1}\supset \wedge_{1}(F_{1})$. Finally we will get the following almost commuting diagram:
$$
\xymatrix{
A_{k_1}\ar[rr]\ar[dd]  &  &  A_{k_2}\ar[rr]\ar[dd]  &  &  A_{k_3}\ar[rr]\ar[dd]  &  &  \\
     &  &  &  &  &  &  \\
     B_{l_1}\ar[rr]\ar[ruru]^{(F_{k_{1}},8\varepsilon_{1})}_{(G_{l_{1}},8\varepsilon_{1})} &  &  B_{l_2}\ar[rr]
     \ar[ruru]^{(F_{k_{2}},8\varepsilon_{2})}_{(G_{l_{2}},8\varepsilon_{2})}  &  & B_{l_3} \ar[rr] &  &  \\
 }
$$
Since $\sum_{i=1}^{\infty}8\varepsilon_{i}<+\infty$, we get approximately intertwining between $A\!=\!lim(A_{n},~\phi_{nm})$ and $B\!=\!lim(B_{n},~\psi_{nm}).$

Hence our main result Theorem 7.1 is proved.\\

$\\$ Guihua Gong, Department of Mathematics, University of Puerto Rico at Rio Piedras, PR 00936, USA\\
%email address:  guihua.gong@upr.edu
Chunlan Jiang, College of Mathematics and Information Science, Hebei Normal University,  Shijiazhuang, Hebei, 050024,    China\\
%email address:  cljiang@hebtu.edu.cn
Liangqing Li, Department of Mathematics, University of Puerto Rico at Rio Piedras, PR 00936, USA
email address:  liangqing.li@upr.edu\\

\clearpage

\begin{tiny}

\begin{small}

\end{small}

\end{tiny}

\end{document}